\documentclass[a4paper,12pt]{report}

\newif\ifpdf
\ifx\pdfoutput\undefined
\pdffalse % we are not running PDFLaTeX
\else
\pdfoutput=1 % we are running PDFLaTeX
\pdftrue
\fi

\ifpdf
\usepackage[pdftex]{graphicx}
\else
\usepackage{graphicx}
\fi

\usepackage{calc}   % permet de faire des opérations arithmétiques sur
\fontencoding{T1}
\usepackage[latin1]{inputenc}

\usepackage[margin=2.53cm]{geometry}

\usepackage{theorem}

\theoremstyle{break}\theorembodyfont{\itshape}
\newtheorem{theo}{Théorème}[section]

\newtheorem{prop}[theo]{Proposition}
\newtheorem{lemm}[theo]{Lemme}
\newtheorem{coro}[theo]{Corollaire}

\theoremstyle{break}\theorembodyfont{\upshape}
\newtheorem{exem}{Exemple}[section]
\newtheorem{rema}{Remarque}[section]

\usepackage[english,french]{babel}

\newcommand{\Dire}[1]{\guillemotleft~#1~\guillemotright{}}

\usepackage{mathrsfs}  % police RSFS : \mathscr
\usepackage{amsmath,  % \text, \boldsymbol , \pmb, 
            amssymb} % \mathbb, \mathfrak, msam et msbm
\usepackage[all]{xy} % diagrammes commutatifs
\usepackage{eucal} % redéfinit \mathcal
\usepackage{stmaryrd}
\newcommand{\BBA}{\mathbb{A}} 
\newcommand{\BBB}{\mathbb{B}} 
\newcommand{\BBC}{\mathbb{C}} \newcommand{\BBP}{\mathbb{P}}
 
\newcommand{\BBE}{\mathbb{E}} \newcommand{\BBR}{\mathbb{R}}
 \newcommand{\BBS}{\mathbb{S}}

\newcommand{\BBM}{\mathbb{M}} 

 \newcommand{\MFU}{\mathfrak{U}}

 \newcommand{\MFs}{\mathfrak{s}}
 \newcommand{\MFt}{\mathfrak{t}}

\newcommand{\Dem}{\textbf{Démonstration}\ {}}
\newcommand{\QED}{\hbox{}\nobreak\hfill
    \quad\hbox{\ding{111}}}
\newcommand{\QEDb}{\hbox{}\nobreak\hfill
    \quad\hbox{\ding{111}}\linebreak}

\usepackage{pifont}
\usepackage{txfonts}

\DeclareMathOperator{\Id}{\mathrm{Id}}
\DeclareMathOperator{\PB}{\mathrm{PB}}
\DeclareMathOperator{\FONC}{\mathscr{F}\!\!\mathfrak{onc}}
\DeclareMathOperator{\GFONC}{\mathscr{G\!F}\!\!\mathfrak{onc}}
\DeclareMathOperator{\GMOR}{\mathscr{G\!M}\!\!\mathfrak{or}}
\DeclareMathOperator{\GSQUARE}{\mathscr{G\!C}\!\mathfrak{arr\acute{e}e}}
\DeclareMathOperator{\GMORALEX}{\mathscr{G\!M}\!\!\mathfrak{or}
                                                       \mathscr{A}\!\mathfrak{lex}}
\DeclareMathOperator{\NAT}{\mathscr{N}\!\!\mathfrak{at}}
\DeclareMathOperator{\Ima}{\mathrm{Im}}
\DeclareMathOperator{\FIX}{\mathrm{Fix}}
\DeclareMathOperator{\Ob}{\mathrm{Ob}}
\DeclareMathOperator{\Mor}{\mathrm{Mor}}
\DeclareMathOperator{\Ouv}{\mathscr{O}\!\mathfrak{uv}}

\DeclareMathOperator{\GTrans}{\mathscr{G\!T}\!\mathfrak{rans}}
\DeclareMathOperator{\GTRANSALEX}{\mathscr{G\!T}\!\!\mathfrak{rans}
                                                       \mathscr{A}\!\mathfrak{lex}}
\DeclareMathOperator{\OP}{\#}
\DeclareMathOperator{\PART}{\mathscr{P}}
\DeclareMathOperator{\Card}{\mathrm{Card}}
\DeclareMathOperator{\MorTo}{\xrightarrow{\ \star\ \ }}
\usepackage{makeidx}
\makeindex
\title{Les groupements}
\author{Freddy Bonnin}
\date{8 avril 2004}
\begin{document}
\setlength{\parindent}{0cm}
%%%%%%
%%%%%%
\maketitle
%
%\begin{abstract}
%
%\end{abstract}
%
\tableofcontents
%%%%%%%%%%%%%%%%%%%%%%%%%%%%%%%%%%%%%%%%%%%%%%%%%%
%%%%%%%%%%%%%%%%%%%%%%%%%%%%%%%%%%%%%%%%%%%%%%%%%%
\setcounter{chapter}{-1}
\chapter{Introduction}
Je me dois de prévenir le lecteur que ce qui est exposé dans ce mémoire n'est 
qu'une tentative de généralisation de la notion de catégorie.  Une tentative et rien de plus. 
Sur bien des points, cet essai n'est pas tout à fait finalisé mais les principaux objectifs 
que nous nous étions donnée ont été à peu près atteints.

Depuis très longtemps, différentes extensions de la notion de catégorie ont été 
envisagées. Une des premières a l'avoir été fut la notion bien connue de 
graphe. Cette structure fut et est encore très étudier pour ses propriétés combinatoires et
son utilisation possible dès que l'on souhaite étudier des choses qui ont une origine et une 
fin. En effet un graphe est tout simplement un ensemble muni de deux applications qui
indiquent la source et le but de chaque élément. De part leur simplicité ont les rencontre
souvent dans des applications pratiques (informatique et économie). Mais c'est une 
extension qui présente le défaut de ne pas avoir de composition. Pour résoudre ce problème 
on s'intéresse souvent à la catégorie libre engendrée par le graphe.

D'autres extensions ont été étudiée. Leur point commun est un affaiblissement de la 
composition. Pour être plus précis, on s'intéresse ici à un graphe $(G;s;t)$ muni
d'une opération binaire partiellement définie (\textit{i.e.} $\circ:G_{1}\to G$, où
$G_{1}$ est un sous-ensemble de $G\times G$) soumise aux conditions suivantes :
si $g\circ f$ existe alors $s(g)=t(f)$, $s(g\circ f)=s(f)$ et $t(g\circ f)=t(g)$. Si on 
suppose de plus que les compositions $t(f)\circ f$, $f\circ s(f)$ existent toujours et 
sont égales à $f$ alors on obtient un graphe multiplicatif fortement 
identitaire que Charles Ehresmann \cite{Ehr65} a plus simplement appelé graphe multiplicatif ou 
néocatégorie. Lurtz Schöder et Horst Herrlich \cite{ScHe00} ont quand à eux affaibli la condition 
précédente en n'imposant pas l'existence de $t(f)\circ f$ et $f\circ s(f)$ mais en les 
obligeant à être égales à $f$ si elles existent. Ils parlent alors de graphe multiplicatif 
faiblement identitaire. En supposant que la composition est aussi associative, ils aboutissent 
à la notion de semicatégorie. Paulo Mateus, Amilcar Sernadas et Cristina Sernadas \cite{MSS99} ont 
utilisé les précatégories (en même temps néocatégorie et semicatégorie) pour étudier 
la combinaison des automates probabilistes (informatique).

Aussi intéressantes soient-ils, tous ces cas présentent le défaut d'imposer le fait que 
les éléments sources et buts sont nécessairement des identités à partir du moment où ils 
sont composables. Cette condition s'avère un obstacle très difficile (voir impossible à lever)
quand on souhaite généraliser certaine construction géométrique. Par exemple, 
Mikhail Kapranov et Vladimir Voevodsky \cite{KaVo91} ont échoué dans leur tentative de généralisation 
des espaces de lacets de Moore en dimension supérieur car leur construction ne satisfait 
pas les conditions sur les identités. Il semblait par conséquent intéressant de voir ce qui 
pouvait se passer si on supprimait ces conditions identitaires.

Le premier chapitre est consacré à quelques rappels simples sur les pullbacks  et les 
univers de Grothendieck. Bien que ces derniers ne fournissent pas une base ensembliste très
satisfaisante, nous nous en contenterons.

Le second chapitre est juste une mise en forme 
la propre possible d'une manière classique de définir les catégories.
On peut par exemple ce référer au le livre de Saunders MacLane \cite{ML71}.
Ordinairement une catégorie est composé d'un ensemble d'objets, d'une famille d'ensembles
indexée par les couples d'objets dont les éléments forment les morphismes de la catégorie et 
d'une composition des morphismes vérifiant les axiomes d'identité et l'associativité. En fait 
tout cela peut être simplifier en ne s'intéressant qu'à l'ensemble des morphismes muni 
de deux applications, source et but, et d'une composition satisfaisant certains axiomes.
Cette construction classique pour les catégories est beaucoup plus difficile à trouver dans la 
littérature pour les foncteurs et les transformations naturelles. Bien que n'apportant rien de bien 
nouveau, ce travail est intéressant car il fournit des idées importantes pour toute la suite.

Le chapitre trois est le coeur du mémoire. Ici on définit la notion de groupement qui est 
tout simplement la généralisation souhaitée. Nous y avons aussi défini les morphismes de 
groupements que nous avons appelés g-morphismes. L'absence d'identités rend très difficiles 
la définition adéquate d'une notion de transformations entre g-morphismes. L'un des principaux
obstacles est la difficulté de construire des applications source et but sur l'ensemble 
de telles transformations et donc de parler de composition. Nous avons abordé cette question à la 
fin du chapitre dans le but de montrer où se situait les obstacles principaux.

Le quatrième chapitre est né de la volonté de trouver un exemple non trivial de groupement.
Or l'espace des chemins de Moore d'un espace topologique est trivialement muni 
d'une composition. Il suffit de juxtaposer les chemins. Cette dernière est bien 
associative mais ne vérifie pas les conditions identitaires habituellement imposé à une 
catégorie. Puisque dès le départ nous les avons supprimé de notre théorie, les groupements
fournissent un cadre naturel pour l'étude des chemins de Moore. Bien sûr s'il n'y avait 
que cela, l'intérêt serait limité. Mais il se trouve qu'il est extrêmement facile de généraliser
cette construction en dimension supérieure. C'est ce que nous avons fait avec les surfaces 
et plus généralement les $I$-espaces de Moore.

Le cinquième chapitre est une tentative pour fournir un cadre assez général pour définir 
des transformations entre g-morphismes. Ce chapitre n'est pas tout à fait satisfaisant mais 
je l'es\-père peut fournir quelques idées intéressantes.

Le dernier chapitre a pour unique but de donner quelques exemples d'ensembles 
possèdant deux structures différentes  de groupements et d'en tirer une définition 
de 2-groupements.

Je tiens à remercier Bertrand Toen pour m'avoir proposé un sujet portant sur les bicatégories. 
Bien que ce mémoire semble en être éloigné, il en est la conséquence directe. En effet 
c'est en essayant de donner une définition dans laquelle on ne ferait plus référence à un ensemble de
 0-cellules et à un ensemble de 1-cellules, mais simplement à un unique ensemble de cellules que nous 
 avons été amené à définir catégories, foncteurs et transformations naturelles sans utiliser d'objets. 
 En se faisant, nous avons constater qu'une grande partie de la théorie pouvait être faite sans 
 les conditions identitaires. Et l'exemple des surfaces de Moore et leurs généralisations possibles en 
 dimensions supérieures nous ont convaincus que ce point de vue méritait d'être exposé.
 
Les mots me manquent pour remercier Carlos Simpson pour le soutien et l'aide qu'il a su m'apporter.
Rien n'aurait pu être fait sans ses remarques pertinentes, sa gentillesse et sa patience.
%%%%%%%%%%%%%%%%%%%%%%%%%%%%%%%%%%%%%%%%%%%%%%%%%%
\chapter{Préliminaires ensemblistes}

Le seul but de ce chapitre est de rappeler quelques propriétés des pullbacks et surtout de
préciser nos notations et le cadre ensembliste dans lequel nous travaillerons.
%%%%%%%%%%%%%%%%%
%%% section 1 %%%
%%%%%%%%%%%%%%%%%
\section{Une propriété intéressante des ensembles : le pullback}\label{sec:preli:1}
\subsection{Définitions et premières propriétés}
\label{subsec:preli:1.1} 
Les ensembles ont la propriété que si deux applications $f_1:A_1\to B$ et
$f_2:A_2\to B$ ont le même but  
$B$ alors il existe deux applications $g_1:C\to A_1$ et $g_2:C\to A_2$ telles que
\begin{enumerate}
\item $f_1g_1=f_2g_2$
\item et si $h_1:D\to A_1$ et $h_2:D\to A_2$ sont deux applications vérifiant
$f_1h_1=f_2h_2$
alors il existe une unique application $h:D\to C$ satisfaisant $h_1=g_1h$ et
$h_2=g_2h$.
\end{enumerate}
\begin{rema}\label{rema:preli:1.1}
Nous utiliserons toujours la notation $fg$ pour parler de la composition des 
applications $f$ et $g$. Lorsque nous rencontrerons des compositions différentes 
de la composition naturelle d'applications d'ensembles nous utiliserons toujours 
d'autres notations.
\end{rema}
\begin{rema}\label{rema:preli:1.2} 
Une construction peut être donnée en prenant $C=\left\lbrace\left(x_1;x_2\right)\in
A_1\times A_2\mid f\left(x_1\right)=f_2\left(x_2\right)\right\rbrace$ et $g_1$, $g_2$ les restrictions des projections canoniques.
\end{rema}
Dans le langage des catégories cela signifie que la \Dire{catégorie} des ensembles
possède des \emph{pullbacks}\index{pullback} ou encore que, ci-dessous, le diagramme de gauche peut être
compléter pour que le carré de droite soit cartésien.
$$\xymatrix{
 & A_1 \ar[d]^{f_1}\\ A_2 \ar[r]_{f_2} & B}
\qquad
\xymatrix{
C \ar[r]^{g_1} \ar[d]_{g_2} & A_1 \ar[d]^{f_1}\\
A_2 \ar[r]_{f_2} & B}$$
En fait le pullback n'est pas $C$ comme on le dit bien souvent mais le 
couple d'application $\left(g_1;g_2\right)$. On notera
$\left(g_1;g_2\right)=\PB\left(f_1;f_2\right)$ et $C=A_1\times_B A_2$.
Pour l'instant cette notation veut simplement dire que
$\left(g_1;g_2\right)$ est \textbf{un} pullback de $\left(f_1;f_2\right)$. 
En effet il n'est pas unique au sens fort du terme, mais unique à une 
bijection près comme cela va être expliqué dans le lemme ci-dessous.
%%%%%%% 
\begin{prop}\label{prop:preli:1.1} 
Soient $f_1:A_1\to B,\ f_2:A_2\to B$ et $f_3: A_3\to B$ trois applications
d'ensembles ayant le même but.
\begin{enumerate}
\item Si $\left(g_1;g_2\right)=\PB\left(f_1;f_2\right)$ alors
$\left(g_2;g_1\right)=\PB\left(f_2;f_1\right)$;
\item Si $\left(g_1;g_2\right)=\PB\left(f_1;f_2\right)$ et 
$\left(h_1;h_2\right)=\PB\left(f_1;f_2\right)$ alors il existe une unique
bijection $l$ telle que $h_1=g_1l$ et $h_2=g_2l$.
\item Si $\left(g_1;g_2\right)=\PB\left(f_1;f_2\right)$,
$\left(g_3;g_4\right)=\PB\left(f_2;f_3\right)$,
$\left(h_1;h_2\right)=\PB\left(f_1g_1;f_3\right)$ et\\
$\left(h_3;h_4\right)=\PB\left(f_1;f_2g_3\right)$, alors il existe une unique bijection
$$\varphi:A_1\times_B\left(A_2\times_B A_3\right)\to\left(A_1\times_B
A_2\right)\times_B A_3$$ 
telle que $$\left(g_1h_1\right)\varphi=h_3,\ 
\left(g_2h_1\right)\varphi=g_3h_4\ \mathrm{et}\ h_2\varphi=g_4h_4$$
\end{enumerate} 
\end{prop}
\Dem 
\begin{enumerate}
\item \'Evident d'après la définition.
\item D'après la première condition sur les pullbacks, $f_1h_1=f_2h_2$. Donc la
seconde, implique l'existence de l'application $l$. De la même façon, il existe une
unique application $l'$ telle que $g_1=h_1l'$ et $g_2=h_2l'$. On en déduit que $ll'$ 
est une application satisfaisant $h_1=h_1(ll')$ et $h_2=h_2(ll')$. D'après la 
deuxième propriété des pullbacks, une seule application peut satisfaire ces égalités.
C'est l'identité. Donc $ll'=\Id$. De même $l'l=\Id$. Ce qui signifie que $l$ est une
bijection.
\item Puisque $\left(g_1;g_2\right)=\PB\left(f_1;f_2\right)$ et 
$f_1h_3=f_2\left(g_3h_4\right)$, il existe une unique application $\varphi_1$ 
pour laquelle $g_1\varphi_1=h_3$ et $g_2\varphi_1=g_3h_4$. D'où
$$\left(f_1g_1\right)\varphi_1=f_1\left(g_1\varphi_1\right)=f_1h_3
=f_2\left(g_3h_4\right)=\left( f_2g_3\right)h_4=\left( f_3g_4\right)h_4
=f_3\left(g_4h_4\right)$$
Or, par hypothèse, $\left(h_1;h_2\right)=\PB\left(f_1g_1;f_3\right)$. Ainsi il 
existe une unique application $\varphi$ satisfaisant
$\varphi_1=h_1\varphi$ et $g_4h_4=h_2\varphi$. Si nous remontons le raisonnement 
nous voyons que $\varphi$ est la seule application vérifiant les égalités
$\left(g_1h_1\right)\varphi=h_3$, $\left(g_2h_1\right)\varphi=g_3h_4$ et
$h_2\varphi=g_4h_4$.

De même il existe une unique application $\varphi'$ qui vérifie
$h_3\varphi'=g_1h_1$, $\left(g_3h_4\right)\varphi'=g_2h_1$ et
$g_4h_4\varphi'=h_2$. Des différentes égalités que nous venons de trouver on déduit
$$h_3\left(\varphi'\varphi\right)=\left(g_1h_1\right)\varphi=h_3,\ 
g_3\left(h_4\left(\varphi'\varphi\right)\right)=\left(g_2h_1\right)\varphi=g_3h_4\ 
\mathrm{et}\  
g_4\left(h_4\left(\varphi'\varphi\right)\right)=h_2\varphi=g_4h_4$$
Comme $\left(g_3;g_4\right)=\PB\left(f_2;f_3\right)$, on trouve 
$$h_3\left(\varphi'\varphi\right)=h_3\ 
\mathrm{et}\  
h_4\left(\varphi'\varphi\right)=h_4$$
Et finalement $\varphi'\varphi=\Id$ car $\left(h_3;h_4\right)
=\PB\left(f_1;f_2g_3\right)$. De la même façon $\varphi\varphi'=\Id$. Ce qui montre 
bien que $\varphi$ est une bijection.\QED
\end{enumerate}
%%%%%
Par la suite nous considérerons toujours que dans la notation
$\left(g_1;g_2\right)=\PB\left(f_1;f_2\right)$
les applications $g_1$ et $g_2$ sont construites comme décrites dans la remarque
\ref{rema:preli:1.2}.
\subsection{Quelques résultats techniques}\label{subsec:preli:1.2} 
Considérons deux applications d'ensembles $s,t:A\to B$.
Posons $\left(\pi_s^t;\pi_t^s\right)=\PB\left(s;t\right)$ et notons 
$A\times_BA$ 
le domaine commun à $\pi_s^t$ et $\pi_t^s$.
\begin{rema}\label{rema:preli:1.3} 
Il peut parfois être commode de voir $A\times_BA$ comme l'ensemble des 
couples $(x;y)\in A$ vérifiant $s(x)=t(y)$. Et $\pi_s^t\left(x;y\right)=x$, 
$\pi_t^s\left(x;y\right)=y$
\end{rema}
\'Etant donné $\left(s_{1},t_{1}:A_1\to B_1\right)$ et 
$\left(s_{2},t_{2}:A_2\to B_2\right) $ deux couples d'applications 
d'ensembles,
si $f,g:A_1\to A_2$ sont deux applications d'ensembles telles que
\begin{itemize}
\item $s_2f\pi_{s_1}^{t_1}=t_2g\pi_{t_1}^{s_1}$ alors il est évident qu'il 
existe une unique 
application $$f\times_c g:A_1\times_{B_1}A_1\to A_2\times_{B_2}A_2$$
satisfaisant
$$\pi_{s_2}^{t_2}\left(f\times_c g\right)=f\pi_{s_1}^{t_1}\ \mathrm{et}\ 
\pi_{t_2}^{s_2}\left(f\times_c g\right)=g\pi_{t_1}^{s_1}$$
\item $s_2f=t_2g$ alors il existe une unique application 
$$\left\langle f;g\right\rangle:A_1\to A_2\times_{B_2}A_2$$ vérifiant
$$f=\pi_{s_2}^{t_2}\left\langle f;g\right\rangle\ \mathrm{et}\ 
g=\pi_{t_2}^{s_2}\left\langle f;g\right\rangle$$
\end{itemize}
\begin{rema}\label{rema:preli:1.4} 
Avec la construction du pullback vue dans la section précédente, on a
$$f\times_c g\left(x;y\right)=\left(f(x);g(y)\right)\qquad\mathrm{et}
\qquad\left\langle f;g\right\rangle(x)=\left(f(x);g(y)\right)$$
\end{rema}
\begin{prop}\label{prop:preli:1.2} 
On considère $\left(s_1,t_1:A_1\to B_1\right)$, 
$\left(s_2,t_2:A_2\to B_2\right)$ et $\left(s_3,t_3:A_3\to B_3\right)$
trois couples d'applications d'ensembles.
\begin{enumerate}
\item Soient $f,g:A_1\to A_2$ et
$f',g':A_2\to A_3$ quatre applications d'ensembles telles que
$$s_2f\pi_{s_1}^{t_1}=t_2g\pi_{t_1}^{s_1}\qquad\mathrm{et}\qquad
s_3f'\pi_{s_2}^{t_2}=t_3g'\pi_{2_1}^{s_2}$$
Alors $f\times_c g$, $f'\times_c g'$ et $\left(f'f\right)\times_c\left(
g'g\right)$ existent et on a de plus
$$\left(f'f\right)\times_c\left(g'g\right)=\left(f'\times_c g'\right)
\left(f\times_c g\right)$$
\item Soient $f,g:A_1\to A_2$ et
$f',g':A_2\to A_3$ quatre applications d'ensembles telles que
$$s_2f=t_2g\qquad\mathrm{et}\qquad
s_3f'\pi_{s_2}^{t_2}=t_3g'\pi_{2_1}^{s_2}$$
Alors $\left\langle f;g\right\rangle$, $f'\times_c g'$ et
$\left\langle f'f;g'g\right\rangle$ existent et
$$\left\langle f'f;g'g\right\rangle=\left(f'\times_c g'\right)
\left\langle f;g\right\rangle$$
\item Soient $f:A_1\to A_2$ et
$g,h:A_2\to A_3$ trois applications d'ensembles telles que
$$s_3g=t_3h$$
Alors $\left\langle g;h\right\rangle$ et $\left\langle gf;hf\right\rangle$
existent avec
$$\left\langle gf;hf\right\rangle=\left\langle g;h\right\rangle f$$
\end{enumerate} 
\end{prop}
\Dem Dans chacun des cas, les hypothèses impliquent l'existence 
de toutes les applications sauf de la dernière.
\begin{enumerate}
\item 
$s_3\left(f'f\right)\pi_{s_1}^{t_1}=s_3f'\left(f\pi_{s_1}^{t_1}\right)
=s_3f'\pi_{s_2}^{t_2}\left(f\times_c g\right)
=t_3g'\pi_{t_2}^{s_2}\left(f\times_c g\right)
=t_3g'\left(g\pi_{t_1}^{s_1}\right)=t_3\left(g'g\right)\pi_{t_1}^{s_1}$
D'où l'existence de $\left(f'f\right)\times_c\left(g'g\right)$. Par 
définition c'est la seule application qui vérifie
$$f'f\pi_{s_1}^{t_1}=\pi_{s_3}^{t_3}\left(\left(f'f\right)
\times_c\left(g'g\right)\right)\quad\mathrm{et}\quad
g'g\pi_{t_1}^{s_1}=\pi_{t_3}^{s_3}\left(\left(f'f\right)
\times_c\left(g'g\right)\right)$$
Or$$\pi_{s_3}^{t_3}\left(f'\times_c g'\right)\left(f\times_c g\right)
=f'\pi_{s_2}^{t_2}\left(f\times_c g\right)=f'f\pi_{s_1}^{t_1}$$
et
$$\pi_{t_3}^{s_3}\left(f'\times_c g'\right)\left(f\times_c g\right)
=g'\pi_{t_2}^{s_2}\left(f\times_c g\right)=g'g\pi_{t_1}^{s_1}$$
D'où l'égalité recherchée.
\item L'existence de $\left\langle f'f;g'g\right\rangle $ est donnée par 
le calcul
$$s_3f'f=s_3f'\pi_{s_2}^{t_2}\left\langle f;g\right\rangle
=t_3g'\pi_{t_2}^{s_2}\left\langle f;g\right\rangle
=t_3g'g$$ Comme
$$\pi_{s_3}^{t_3}\left(f'\times_c g'\right)\left\langle f;g\right\rangle
=f'\pi_{s_2}^{t_2}\left\langle f;g\right\rangle
=f'f$$
$$\pi_{t_3}^{s_3}\left(f'\times_c g'\right)\left\langle f;g\right\rangle
=g'\pi_{t_2}^{s_2}\left\langle f;g\right\rangle
=g'g$$
et que, par définition, $\left\langle f'f;g'g\right\rangle$ est la seule 
application à vérifier ces égalités, on trouve bien le résultat annoncé.
\item \'Evident.\QED
\end{enumerate} 
\begin{prop}\label{prop:preli:1.3} 
Considérons $s,t:A\to B$ deux applications ayant même source et même but.
\begin{enumerate}
\item A fin de simplifier nos notations, posons 
\begin{itemize}
\item $\left(g_s;g_t\right)=\PB\left(s;t\right)$,
\item $\left(f_s;f_t\right)=\PB\left(sg_t;t\right)$ et
\item $\left(h_s;h_t\right)=\PB\left(s;tg_s\right)$
\end{itemize} 
Alors il existe une unique bijection 
$$\omega:\left(A\times_B A\right)\times_B A\to A\times_B\left(A\times_B
A\right)$$ qui vérifie les égalités ci-dessous
$$h_s\omega=g_sf_s,\qquad g_sh_t\omega=g_tf_s\qquad\mathrm{et}\qquad
g_th_t\omega=f_t$$
\item De plus, si $h_1,h_2,h_3:C\to A$ sont trois applications 
d'ensembles telles que
$$sh_1=th_2\qquad\mathrm{et}\qquad sh_2=th_3$$
alors on peut bien évidemment construire les applications 
$$\left\langle h_1;h_2\right\rangle:C\to A\times_B A,\qquad
\left\langle h_2;h_3\right\rangle:C\to A\times_B A$$
et il existe une unique application 
$$\left\langle\left\langle h_1;h_2\right\rangle;h_3\right\rangle:
C\to\left(A\times_B A\right)\times_B A$$ telle que 
$$g_sf_s\left\langle\left\langle h_1;h_2\right\rangle;h_3\right\rangle
=h_1,\quad 
g_tf_s\left\langle\left\langle h_1;h_2\right\rangle;h_3\right\rangle
=h_2,\quad 
f_t\left\langle\left\langle h_1;h_2\right\rangle;h_3\right\rangle
=h_3$$
ainsi qu'une unique application
$$\left\langle h_1;\left\langle h_2;h_3\right\rangle\right\rangle
:C\to A\times_B\left(A\times_B A\right)$$ 
vérifiant
$$h_s\left\langle\left\langle h_1;h_2\right\rangle;h_3\right\rangle
=h_1,\quad 
g_sh_t\left\langle\left\langle h_1;h_2\right\rangle;h_3\right\rangle
=h_2,\quad 
g_th_t\left\langle\left\langle h_1;h_2\right\rangle;h_3\right\rangle
=h_3$$
De plus on a la relation
$$\omega\left\langle\left\langle h_1;h_2\right\rangle;h_3\right\rangle
=\left\langle h_1;\left\langle h_2;h_3\right\rangle\right\rangle$$
\item Si $A=B$, alors les applications $\left\langle \Id;s\right\rangle
\times_c\Id$ et $\Id\times_c\left\langle t;\Id\right\rangle$ sont bien 
définies et $$\omega\left(\left\langle \Id;s\right\rangle
\times_c\Id\right)=\Id\times_c\left\langle t;\Id\right\rangle$$
\end{enumerate}
\end{prop}
\Dem
\begin{enumerate}
\item Puisque $\left(g_s;g_t\right)=\PB\left(s;t\right)$ et que
$s\left(g_tf_s\right)=tf_t$ (ou encore $\left(sg_t\right)f_s=tf_t$),
il existe une unique application $\omega_1$ pour laquelle
$$g_s\omega_1=g_tf_s\qquad\mathrm{et}\qquad g_t\omega_1=f_t$$
On a alors
$$\left(tg_s\right)\omega_1=t\left(g_s\omega_1\right)=t\left(g_tf_s\right)
=\left(tg_t\right)f_s=\left(sg_s\right)f_s=s\left(g_sf_s\right)$$
Comme $\left(h_s;h_t\right)=\PB\left(s;tg_s\right)$, on trouve bien l'unique 
application $\omega$ vérifiant
$$h_s\omega=g_sf_s\qquad\mathrm{et}\qquad h_t\omega=\omega_1$$
c'est-à-dire
$$h_s\omega=g_sf_s,\qquad g_sh_t\omega=g_tf_s\qquad\mathrm{et}\qquad
g_th_t\omega=f_t$$
De la même manière, on trouve une unique application $\omega'$ telle que
$$f_t\omega'=g_th_t,\qquad g_tf_s\omega'=g_sh_t\qquad\mathrm{et}\qquad
g_sf_s\omega'=h_s$$
Par conséquent $\omega'\omega$ est la seule application à satisfaire les 
égalités
$$f_t\left(\omega'\omega\right)=f_t,\qquad g_tf_s\left(\omega'\omega\right)=
g_tf_s,\qquad g_sf_s\left(\omega'\omega\right)=g_sf_s$$
et $\omega\omega'$ la seule à satisfaire
$$h_s\left(\omega\omega'\right)=h_s,\qquad g_sh_t\left(\omega\omega'\right)=
g_sh_t,\qquad g_th_t\left(\omega\omega'\right)=g_th_t$$
Ce qui implique $\omega'\omega=\Id$ et $\omega\omega'=\Id$.
\item Rappelons que $\left\langle h_1;h_2\right\rangle$ est l'unique application 
telle que
$$g_s\left\langle h_1;h_2\right\rangle=h_1\quad\mathrm{et}\quad
g_t\left\langle h_1;h_2\right\rangle=h_2$$
Comme
$sg_t\left\langle h_1;h_2\right\rangle=sh_2=th_3$,
il existe bien une unique application 
$$\left\langle\left\langle h_1;h_2\right\rangle;h_3\right\rangle:
C\to\left(A\times_B A\right)\times_B A$$
vérifiant
$$f_s\left\langle\left\langle h_1;h_2\right\rangle;h_3\right\rangle
=\left\langle h_1;h_2\right\rangle\quad\mathrm{et}\quad
f_t\left\langle\left\langle h_1;h_2\right\rangle;h_3\right\rangle
=h_3$$
D'où l'existence et l'unicité annoncées. De même pour
$$\left\langle h_1;\left\langle h_2;h_3\right\rangle\right\rangle$$
Or on a
$$h_s\omega\left\langle\left\langle h_1;h_2\right\rangle;h_3\right\rangle
=g_sf_s\left\langle\left\langle h_1;h_2\right\rangle;h_3\right\rangle
=h_1$$
$$g_sh_t\omega\left\langle\left\langle h_1;h_2\right\rangle;h_3\right\rangle
=g_tf_s\left\langle\left\langle h_1;h_2\right\rangle;h_3\right\rangle
=h_2$$
$$g_th_t\omega\left\langle\left\langle h_1;h_2\right\rangle;h_3\right\rangle
=f_t\left\langle\left\langle h_1;h_2\right\rangle;h_3\right\rangle
=h_3$$
Par unicité, on obtient la relation voulue.
\item On a
$$s\Id=ts,\quad sg_t\left\langle\Id;s\right\rangle g_s=ssg_s=sg_s=tg_t=
t\Id g_t$$
D'où l'existence de $\left\langle\Id;s\right\rangle\times_c\Id$. On montre 
aussi facilement l'existence de $\Id\times_c\left\langle t;\Id\right\rangle$.\\
Le reste est immédiat d'après les calculs ci-dessous et l'unicité.
$$h_s\omega\left(\left\langle\Id;s\right\rangle\times_c\Id\right)
=g_sf_s\left(\left\langle\Id;s\right\rangle\times_c\Id\right)
=g_s\left\langle\Id;s\right\rangle g_s
=\Id g_s=g_s$$
$$g_sh_t\omega\left(\left\langle\Id;s\right\rangle\times_c\Id\right)
=g_tf_t\left(\left\langle\Id;s\right\rangle\times_c\Id\right)
=g_t\left\langle\Id;s\right\rangle g_s
=sg_s=tg_t$$
$$g_th_t\omega\left(\left\langle\Id;s\right\rangle\times_c\Id\right)
=f_t\left(\left\langle\Id;s\right\rangle\times_c\Id\right)
=\Id g_t=g_t$$
et
$$h_s\left(\Id\times_c\left\langle t;\Id\right\rangle\right)
=\Id g_s=g_s$$
$$g_sh_t\left(\Id\times_c\left\langle t;\Id\right\rangle\right)
=g_s\left\langle t;\Id\right\rangle g_t
=t g_t$$
$$g_th_t\left(\Id\times_c\left\langle t;\Id\right\rangle\right)
=g_t\left\langle t;\Id\right\rangle g_t
=\Id g_t=g_t$$
\QED
\end{enumerate}
\begin{rema}\label{rema:preli:1.5}
Cette proposition est en fait évidente si on pense à 
$\left(A\times_c A\right)\times_c A$ comme étant l'ensemble
$$\left\lbrace\left((x;y);z\right)\mid s(x)=t(y),\ s(y)=t(z)\right\rbrace$$
et à $A\times_c\left(A\times_c A\right)$ comme étant
$$\left\lbrace\left(x;(y;z)\right)\mid s(x)=t(y),\ s(y)=t(z)\right\rbrace$$
On a alors $\omega\left(\left(x;y\right);z\right)=\left(x;\left(y;z\right)\right)$.
\end{rema}
\begin{rema}\label{rema:preli:1.6} 
Bien que l'ensemble $\left(A_1\times_B A_2\right)\times_B A_3$
de la proposition \ref{prop:preli:1.1} et l'ensemble 
$\left(A\times_B A\right)\times_B A$ de la proposition \ref{prop:preli:1.3}
semblent avoir des notations qui peuvent porter à confusion, cela n'est 
pas le cas car elles ne s'appliquent pas dans les mêmes cas.
\end{rema}
%%%%%%%%%%%%%%%%%
%%% section 2 %%%
%%%%%%%%%%%%%%%%%
\section{Les univers de Grothendieck}
\label{sec:preli:2} 
Afin de pouvoir parler de catégorie des ensembles, de catégories des catégories, 
de catégorie de foncteurs, nous sommes obligé de travailler avec des ensembles 
\Dire{suffisamment petits} mais qui se comportent bien vis-à-vis des opérations 
usuelles de la théorie des ensembles. C'est ainsi que la notion d'univers a été 
introduite par Grothendieck.
%%%%%%%%%%%%%%%%%
\subsection{Définition et premières propriétés}\label{subsec:preli:2.1} 
Pour être plus précis, un \emph{univers}\index{univers} est un ensemble
$\MFU$ vérifiant les axiomes suivants
\begin{description}
\item[(U 1)] si $x\in\MFU$ et $y\in x$, alors $y\in\MFU$ ;
\item[(U 2)] si $x\in\MFU$ et $y\in\MFU$, alors $\left\lbrace x;y\right\rbrace\in\MFU$ ;
\item[(U 3)] si $I\in\MFU$ et si, pour chaque $i\in I$, $x_i\in\MFU$, alors 
$\bigcup_{i\in I}x_i\in\MFU$ ;
\item[(U 4)] si $x\in\MFU$, alors $\PART(x)\in\MFU$ où $\PART(x)$ est l'ensemble des 
        sous-ensembles de $x$ ;
\item[(U 5)] $\omega\in\MFU$ où $\omega$ n'est autre que l'ensemble des ordinaux finis.
\end{description}
Ces axiomes impliquent les propriétés de stabilité ci-dessous.
\begin{prop}\label{prop:preli:2.1} 
\begin{enumerate}
\item Si $x\in\MFU$, alors $\{x\}\in\MFU$.
\item Si $x\in\MFU$ et $y\in\MFU$, alors $\left(x;y\right)\in\MFU$.
\item Si $x\in\MFU$ et $y\in\MFU$, alors $x\times y\in\MFU$.
\item Si $x\in\MFU$ et $y\subseteq x$, alors $y\in\MFU$.
\item Si $x\in\MFU$ et $y\in\MFU$, alors $x^{y}\in\MFU$ où $x^{y}$ est l'ensemble
        des applications de $y$ dans $x$.
\item Si $I\in\MFU$ et si, pour chaque $i\in I$, $x_i\in\MFU$, alors 
        $\bigsqcup_{i\in I}x_i\in\MFU$.
\item Si $I\in\MFU$ et si, pour chaque $i\in I$, $x_i\in\MFU$, alors 
        $\prod_{i\in I}x_i\in\MFU$.
\item Si $x\in\MFU$, alors $x\subseteq\MFU$.
\end{enumerate}
\end{prop}
\Dem
\begin{enumerate}
\item C'est un cas particulier de l'axiome (U 1).
\item On sait que $\left(x;y\right)=\left\lbrace x;\left\lbrace x;y\right\rbrace
        \right\rbrace$. Il suffit d'appliquer deux fois l'axiome (U 2). 
\item Soit $i\in x$. Pour tout $j\in y$, on vient de voir que $\left(i;j\right)\in\MFU$.
        Il s'en suit que $\left\lbrace\left(i;j\right)\right\rbrace\in\MFU$. Comme $y\in\MFU$, 
        l'axiome (U 3) nous donne $\bigcup_{j\in y}\left\lbrace\left(i;j\right)\right
        \rbrace\in\MFU$. De même, on obtient $\bigcup_{i\in x}\bigcup_{j\in y}\left\lbrace
        \left(i;j\right)\right\rbrace\in\MFU$. Or par définition $\bigcup_{i\in x}
        \bigcup_{j\in y}\left\lbrace\left(i;j\right)\right\rbrace=x\times y$.
\item Il suffit de remarquer que $y\in\PART(x)\in\MFU$ et d'utiliser l'axiome (U 1).
\item Une application de $y$ dans $x$ peut être vue comme un sous-ensemble particulier 
        de $y\times x$. Par conséquent $x^y\subseteq\PART\left(x\times\right)$. D'après ce que 
        l'on vient de voir, l'ensemble de droite appartient à $\MFU$. Le résultat se déduit 
        alors de la propriété précédente.
\item Posons $\bar{x}_i=x\times\{i\}$. Les propriétés ci-dessus nous permettent 
        d'affirmer que $\bar{x}_i\in\MFU$. D'après l'axiome (U 3), $\bigcup_{i\in I}\bar{x}_i
        \in\MFU$. Or par définition, $\bigcup_{i\in I}\bar{x}_i=\bigsqcup_{i\in I}x_i$.
\item On peut définir $\left(u_i\right)_{i\in I}$ avec $u_i\in x_i$ comme une 
        application de $I$ dans $\bigsqcup_{i\in I}x_i$ telle que l'image de $i$ soit dans 
        $x_i$. Ainsi $\prod_{i\in I}x_i$ peut être vu comme un sous-ensemble de $\left(
        \bigsqcup_{i\in I}x_i\right)^I$. Les propriétés 4, 5 et 6, nous donne le résultat.
\item D'après l'axiome (U 1), tout élément de $x$ est un élément de $\MFU$. Par 
        conséquent, $x$ est une partie de l'univers, $x\subseteq\MFU$.\QED
\end{enumerate}
D'après l'axiome (U 5) et toutes les propriétés de stabilité vérifiées par un 
univers, toutes les constructions mathématiques usuelles peuvent être faites dans un 
univers. En particuliers les notions d'ensembles quotients, de nombres réels, de limites 
inductives et projectives peuvent être définies à l'intérieur d'un univers.

\subsection{Un nouvel axiome}\label{subsec:preli:2.2} 
C'est un euphémisme de dire qu'il est très difficile d'exhiber un univers 
en utilisant uniquement les axiomes de la théorie des ensembles ZFC. 
Donc afin de pouvoir travailler 
avec les univers, nous sommes amener à faire l'hypothèse que
\begin{center}
$(\star)$\quad\Dire{tout ensemble appartient à un univers.}
\end{center}  
Ceci fixé, nous pouvons maintenant donner quelques propriétés supplémentaires.
\begin{prop}\label{prop:preli:2.2} 
Toute intersection d'univers est un univers.
\end{prop}
\Dem
Soit $\left(\MFU_i\right)_{i\in I}$ une famille d'univers. Posons $\MFU=\bigcap_{i\in I}
\MFU_i$. Montrons l'axiome (U~1). Soit $i\in I$. Comme $x\in\MFU$, on a $x\in\MFU_i$ et
$y\in x$. $\MFU_i$ étant un univers, on en déduit $y\in\MFU_i$. Puisque c'est vrai pour 
tout $i\in I$, on a $y\in\MFU$. Les démonstrations des quatre autres axiomes sont 
semblables.\QED
\begin{prop}\label{prop:preli:2.3} 
Si $\MFU$ est un univers alors il existe un plus petit univers, que nous 
noterons $\MFU^{+}$, contenant $\MFU$. 
\end{prop}
\Dem
Conséquence immédiate de la proposition précédente et de l'axiome $(\star)$.
\QED

Si $\MFU$ est un univers, nous appellerons \emph{$\MFU$-ensembles}
\index{ensemble!$\MFU$-ensemble} 
les sous-ensembles de $\MFU$ et \emph{petits $\MFU$-ens\-embles}
\index{ensemble!petit $\MFU$-ensemble} ses éléments.

Pour la suite, choisissons un univers $\MFU$ et appellons respectivement ensembles
\index{ensemble} et petits ensembles les $\MFU$-ensembles et petits $\MFU$-ensembles.

%%%%%%%%%%%%%%%%%%%%%%%%%%%%%%%%%%%%%%%%%%%%%%%%%%
\chapter{Les catégories sans objets}\label{chap:cat}
%%%%%%%%%%%%%%%%%
%%% section 0 %%%
%%%%%%%%%%%%%%%%%
\section{Idées directrices}\label{sec:cat:0}
Classiquement, une catégorie $\BBC$ est composée d'un ensemble d'objets $\Ob(\BBC)$, d'un ensemble de morphismes $\Mor(\BBC)$, d'une application source $s:\Mor(\BBC)\to\Ob(\BBC)$, d'une application
but $t:\Mor(\BBC)\to\Ob(\BBC)$ (le $t$ vient de l'anglais \textit{target}), d'une application
identité $\Id:\Ob(\BBC)\to\Mor(\BBC)$ et d'une composition $\OP:\Mor(\BBC)\times_O\Mor(\BBC)\to
\Mor(\BBC)$ qui à tous morphismes $f:X\to Y$ et $f':X'\to Y'$ tels que $X'=s(f')=t(f)=Y'$ associe un
nouveau morphisme $f'\OP f:X\to Y$. Ici nous avons utilisé la représentation usuelle d'un morphisme
$f$ par une flèche $f:X\to Y$ quand $s(f)=X$ et $t(f)=Y$. De plus, pour que ces ensembles et
applications forment bien une catégorie, ils doivent satisfaire les axiomes suivants :
\begin{description}
\item[(associativité)] Pour tous morphismes $f$, $g$ et $h$ vérifiant $s(h)=t(g)$ et $s(g)=t(f)$, on a
$$h\OP(g\OP h)=(h\OP g)\OP f$$
\item[(identité)] Pour tout objet $X$, $s(\Id_X)=X$ et $t(\Id_X)=X$. De plus pour tout morphisme $f$, on a 
$$f\OP\Id_{s(f)}=f\quad\text{et}\quad\Id_{t(f)}\OP f=f$$
\end{description}
Il est a remarqué qu'à chaque objet est associé une unique identité et réciproquement qu'à chaque identité 
est associée un unique objet. Cette remarque nous amènera, dans la prochaine section, à définir les catégories 
en supprimant l'ensemble des objets.
 
Il est évident que les applications source et but se comportent les une par rapport aux autres de la manière 
suivante :
$$s(\Id s)=s\quad;\quad s(\Id t)=t\quad;\quad t(\Id t)=t\quad;\quad t(\Id s)=s$$
Ces égalités et l'identification des objets et identités nous donnerons l'axiome (CAT 1).

En regardant la définition de la composition, on voit que si $f$ et $f'$ sont des morphismes composables alors
$$s(f'\OP f)=s(f)\qquad\text{et}\qquad t(f'\OP f)=t(f')$$
D'où l'axiome (CAT 2).

Les deux derniers axiomes (CAT 3) et (CAT 4) ne sont respectivement que les axiomes d'identité et 
d'associativité.

Usuellement un foncteur $F$ entre la catégorie $\BBC_1$ et la catégorie $\BBC_2$ est une paire
d'applications $F^{\Ob}:\Ob(\BBC_1)\to\Ob(\BBC_2)$ et $F^{\Mor}:\Mor(\BBC_1)\to\Mor(\BBC_2)$ telles
que $F^{\Mor}(\Id_X)=\Id_{F^{\Ob}(X)}$ et $F^{\Mor}(f'\OP_1f)=F^{\Mor}(f')\OP_2F^{\Mor}(f)$
pour tout objet $X$ et tout couple $(f';f)\in\BBC_1\times_O\BBC_1$ de morphismes composables.
Puisque nous allons identifier objets et identités, nous ne conserverons que l'application des
morphismes et les deux conditions précédentes nous donnerons respectivement les axiomes
(FONC 1) et (FONC 2).

Le cas des transformations naturelles est plus complexe. C'est même en fait le point de départ de notre
travail. Habituellement une transformation naturelle $\eta$ entre deux foncteurs $F_1,F_2:\BBC_1\to
\BBC_2$ est une famille de morphismes $\eta_X:F_1^{\Ob}(X)\to F_2^{\Ob}(X)$ indexée par les objets
$\Ob(\BBC_1)$ de la catégorie source $\BBC_1$ pour laquelle le carré
$$\xymatrix{
F_1^{\Ob}(X)\ar[r]^{\eta_X}\ar[d]_{F_1^{\Mor}(f)}  & F_2^{\Ob}(X)\ar[d]^{F_2^{\Mor}(f)}\\
F_1^{\Ob}(Y)\ar[r]_{\eta_Y} & F_2^{\Ob}(Y)}$$
est commutatif pour tout morphisme $f$ de $\BBC_1$.

Nos catégories n'ayant pas d'objets, nos transformations naturelles ne peuvent pas être définies
de cette façon. Nous serons obligé de les définir comme des applications de l'ensemble des morphismes
de la catégorie source dans l'ensemble des morphismes de la catégorie but. Mais puisque nos objets
sont identifiés avec les identités et que, nous le montrerons dans la dernière section de ce chapitre,
l'ensemble de ces dernières n'est autre que l'image de l'application source $s$, nous leur imposerons
l'axiome (NAT 2). De plus en regardant les sources et buts des morphismes de la famille composant la
transformation naturelle, nous somme amené à définir l'axiome (NAT 1). Le dernier axiome (NAT 3)
n'étant rien d'autre que la traduction de la commutativité du carré ci-dessus.

Cette section n'était qu'une introduction sommaire à ce chapitre. Une étude plus approfondie de la
correspondance entre le point de vue classique et le point de vue sans objets sera faite dans la
dernière section.
%%%%%%%%%%%%%%%%%
%%% section 1 %%%
%%%%%%%%%%%%%%%%%
\section{Les catégories}\label{sec:cat:1}

Partant du principe que les objets peuvent être identifiés avec leurs identités, nous sommes
amené à définir une \emph{catégorie}\index{catégorie} (sans objets) $\left(\BBB;s,t,\#\right)$
comme la donnée
\begin{itemize}
\item d'un ensemble $\BBB$
\item de deux applications $s:\BBB\to\BBB$, la \emph{source}\index{source}
et $t:\BBB\to\BBB$, le \emph{but}\index{but},
\item et d'une application $\#:\BBB\times_c\BBB\to\BBB$, la \emph{composition}
\index{composition},
\end{itemize} 
vérifiant les axiomes suivants
\begin{description}
\item[(CAT 1)] $ss=s,\qquad st=t,\qquad tt=t\qquad\text{et}\qquad ts=s$
\item[(CAT 2)] $s\#=s\pi_t^s:\BBB\times_c\BBB\to\BBB\qquad\mathrm{et}\qquad
t\#=t\pi_s^t:\BBB\times_c\BBB\to\BBB$
\item[(CAT 3)] $\#\left\langle\Id_{\BBB};s\right\rangle=\Id_{\BBB},\qquad 
\#\left\langle t;\Id_{\BBB}\right\rangle=\Id_{\BBB}$
\item[(CAT 4)] Avec $\omega:\left(\BBB\times_c\BBB\right)\times_c\BBB\to
\BBB\times_c\left(\BBB\times_c\BBB\right)$
la bijection de la proposition \ref{prop:preli:1.3}, on a
$$\#\left(\#\times_c\Id_{\BBB}\right)=
\#\left(\Id_{\BBB}\times_c\#\right)\omega:
\left(\BBB\times_c\BBB\right)\times_c\BBB\to\BBB$$ 
\end{description} 
Dans cette définition nous parlons des applications
$\left\langle\Id_{\BBB};s\right\rangle$, 
$\left\langle t;\Id_{\BBB}\right\rangle$,  $\Id_{\BBB}\times_c\#$ et
$\#\times_c\Id_{\BBB}$ sans en avoir vérifié l'existence. 
Ce qui est immédiat d'après les calculs suivants
$$ts=s=s\Id_{\BBB},\qquad st=t=t\Id_{\BBB}$$
$$t\left(\#h_t\right)=tg_sh_t=sh_s=s\left(\Id h_s\right),\qquad
s\left(\#f_s\right)=sg_t f_s=tf_t=t\left(\Id f_t\right)$$
où nous avons utilisé les notations de la proposition \ref{prop:preli:1.3} 
($g_s=\pi_s^t$ et $g_t=\pi_t^s$).
\begin{rema}\label{rema:cat:1.1} 
Tout ensemble peut être muni d'une structure de catégorie. En effet si $X$ est 
un ensemble, il suffit de prendre $s=t=\#=\Id_X$. Cette dernière étant possible 
car $\left(\Id;\Id\right)=\PB\left(\Id;\Id\right)$.
\end{rema}
\begin{rema}\label{rema:cat:1.2} 
Pour être plus rigoureux, nous aurions dû parler de $\MFU$-catégorie
\index{catégorie!$\MFU$-catégorie} plutôt que simplement de catégorie.
\end{rema}
Une \emph{petite catégorie}\index{catégorie!petite catégorie} 
$\left(\BBB;s;t;\OP\right)$ est une 
catégorie dont l'ensemble de base $\BBB$ est un petit ensemble.
\begin{lemm}\label{lemm:cat:1.1}
Une petite catégorie est un élément de l'univers.
\end{lemm}
\Dem
La démonstration est basée sur la proposition \ref{prop:preli:2.1}.
Supposons que $\left(\BBB;s;t;\OP\right)$ soit une petit $\MFU$-catégorie.
 Comme $\BBB$ 
est $\MFU$-petit, $\BBB\times\BBB$ est petit (propriété 3)
et donc aussi $s$, $t$ et $\OP$ (propriété 5). Il s'ensuit que le 
quadrulet $\left(\BBB;s;t;\OP\right)$ est un élément de $\MFU$ (propriété 7 
et axiomes (U 1) et (U 5)).
\QEDb
Le lemme suivant est une conséquence immédiate de la définition.
% lemme
\begin{lemm}\label{lemm:cat:1.2}
Si $(\BBB;s;t;\OP)$ est une catégorie alors si $x$ est un élément de $\BBB$ 
tel que $s(x)=x$ ou $t(x)=x$ alors $s(x)=t(x)=x$ et pour tous $y$ et $z$ dans $\BBB$ 
vérifiant $s(y)=t(x)$ et $s(x)=t(z)$, on a 
$$y\OP x=y\quad\text{et}\quad x\OP z=z$$
On dit que l'élément $x$ est une \emph{identité}\index{identité}.
\end{lemm}
\Dem 
Si $s(x)=x$, alors $ts(x)=t(x)$. Or $ts(x)=s(x)$ d'après l'axiome (CAT 1). Donc 
$s(x)=t(x)=x$. De même si $t(x)=x$.\\
De plus comme $s(y)=t(x)$ et $s(x)=t(z)$, on a, d'après l'axiome (CAT 3),
$$y\OP x=y\OP t(x)=y\OP s(y)=y\quad\text{et}\quad
x\OP z=s(x)\OP z=t(z)\OP z=z$$
\QED
%%%%%%%%%%%%%%%%%
%%% section 2 %%%
%%%%%%%%%%%%%%%%%
\section{Les foncteurs}\label{sec:cat:2} 
Un \emph{foncteur}\index{foncteur} $f:\left(\BBB_1;s_1;t_1;\OP_1\right)\to
\left(\BBB_2;s_2;t_2;\OP_2\right)$, souvent
noté $f:\BBB_1\to\BBB_2$, est un triplet $\left(f;(\BBB_1;s_1;t_1;\OP_1);
(\BBB_2;s_2;t_2;\OP_2)\right)$ où $f:\BBB_1\to\BBB_2$ est une application 
d'ensembles  vérifiant les axiomes ci-dessous
\begin{description}
\item[(FONC 1)] $fs_1=s_2f\qquad\mathrm{et}\qquad ft_1=t_2f$
\item[(FONC 2)] $f\#_1=\#_2\left(f\times_c f\right):\BBB_1\times_c\BBB_1\to\BBB_2$
\end{description}
L'application $f\times_c f$ existe car\quad
$s_2\left(f\pi_{s_1}^{t_1}\right)=fs_1\pi_{s_1}^{t_1}=ft_1\pi_{t_1}^{s_1}
=t_2\left(f\pi_{t_1}^{s_1}\right)$
\begin{rema}\label{rema:cat:2.1} 
Il est immédiat que toute application d'ensembles $f:X\to Y$ est un foncteur 
$$f:\left(X;\Id_X;\Id_X;\Id_X\right)\to\left(Y;\Id_Y;\Id_Y;\Id_Y\right)$$
entre les ensembles $X$ et $Y$ munis des structures de catégories vues à la 
remarque \ref{rema:cat:1.1}.
\end{rema}
Donnons une première conséquence, très classique ,de cette définition.
% lemme
\begin{lemm}\label{lemm:cat:2.0}
Si $f:\BBB_{1}\to\BBB_{2}$ est un foncteur de catégories, alors $f(s_{1}(x))$ et 
$f(t_{1}(x))$ sont des identités pour tout $x\in\BBB_{1}$.
\end{lemm}
\Dem
L'axiome (FONC 1) implique
$$f(s_{1}(x))=s_{2}f(x)\quad\text{et}\quad f(t_{1}(x))=t_{2}f(x)$$
Or $s_{2}s_{2}=s_{2}$ et $t_{2}t_{2}=t_{2}$. Le résultat découle alors directement 
du lemme \ref{lemm:cat:1.2}.\QED
% proposition
\begin{prop}\label{prop:cat:2.1} 
\begin{itemize}
\item Si $\left(\BBB;s;t\right)$ est un catégorie, l'application identité 
$\Id_{\BBB}$ définit un foncteur aussi noté $\Id_{\BBB}$.
\item Si $f_1:\left(\BBB_1;s_1;t_1\right)\to\left(\BBB_2;s_2;t_2\right)$ et 
$f_2:\left(\BBB_2;s_2;t_2\right)\to\left(\BBB_3;s_3;t_3\right)$ sont deux 
foncteurs alors la fonction d'ensembles $f_2f_1$ définit elle aussi un 
foncteur 
$$f_2\bullet f_1:\left(\BBB_1;s_1;t_1\right)\to\left(\BBB_3;s_3;t_3\right)$$
\end{itemize}
\end{prop}
\Dem
\begin{itemize}
\item \'Evident.
\item En effet on a
$$\left(f_2f_1\right)s_1=f_2s_2f_1=s_3\left(f_2f_1\right),\qquad
\left(f_2f_1\right)t_1=f_2t_2f_1=t_3\left(f_2f_1\right)$$
et
$$\left(f_2f_1\right)\#_1=f_2\#_2\left(f_1\times_c f_1\right)=
\#_3\left(f_2\times_c f_2\right)\left(f_1\times_c f_1\right)=
\#_3\left(\left(f_2f_1\right)\times_c\left(f_2f_1\right)\right)$$
d'après la proposition \ref{prop:preli:1.2}.\QED
\end{itemize}
\begin{lemm}\label{lemm:cat:2.1}
Si $f:\BBB_1\to\BBB_2$ est un $\MFU$-foncteur entre les petites $\MFU$-catégories 
$\BBB_1$ et $\BBB_2$ alors c'est un élément de l'univers $\MFU$.
\end{lemm}
\Dem
La démonstration est semblable à celle du lemme \ref{lemm:cat:1.1}.  
\QEDb
Posons 
\begin{itemize}
\item $\FONC$ l'ensemble des \emph{petits foncteurs}\index{foncteur!petit foncteur}
(foncteurs entre petites catégories),
\item $s:\FONC\to\FONC,\ s\left(f:\BBB_1\to\BBB_2\right)=\Id_{\BBB_1}$\\
$t:\FONC\to\FONC,\ t\left(f:\BBB_1\to\BBB_2\right)=\Id_{\BBB_2}$
\item $\bullet:\FONC\times_c\FONC\to\FONC,\ \bullet\left(f_2;f_1\right)=
f_2\bullet f_1$
\end{itemize}
\begin{theo}\label{theo:cat:2.1} 
$\left(\FONC;s;t;\bullet\right)$ est une catégorie.
\end{theo}
\Dem Les trois premiers axiomes sont évidents et (CAT 4) se déduit immédiatement de 
l'associativité des applications d'ensembles.\QED
%%%%%%%%%%%%%%%%%
%%% section 3 %%%
%%%%%%%%%%%%%%%%%
\section{Les transformations naturelles}\label{sec:cat:3} 
Une \emph{transformation naturelle}\index{transformation naturelle}
$\eta:\left(f_1:\BBB_1\to\BBB_2\right)\leadsto\left(f_2:\BBB_1\to\BBB_2\right)$ 
est un triplet $$\left(\eta:\BBB_1\to\BBB_2;(f_1:\BBB_1\to\BBB_2);
(f_2:\BBB_1\to\BBB_2)\right)$$
où $\eta:\BBB_1\to\BBB_2$ est une application d'ensembles telle que
\begin{description}
\item[(NAT 1)] $s_2\eta s_1=s_2f_1\qquad\mathrm{et}\qquad t_2\eta s_1=s_2f_2$
\item[(NAT 2)] $\eta=\eta s_1$ 
\item[(NAT 3)] $\#_2\left\langle f_2;\eta s_1\right\rangle=
\#_2\left\langle\eta t_1;f_1\right\rangle$
\end{description}
Pour simplifier, nous noterons souvent $\eta:f_1\leadsto f_2$.
\begin{lemm}\label{lemm:cat:3.0}
L'axiome (NAT 1) est équivalent à l'axiome suivant
\begin{description}
\item[(NAT 1')] $s_2\eta t_1=t_2f_1\qquad\mathrm{et}\qquad t_2\eta t_1=t_2f_2$
\end{description}
\end{lemm}
\Dem Supposons que l'axiome (NAT 1) soit vrai. D'après les axiomes (CAT 1) et (FONC 1), on a
$$s_2\eta t_1=s_2\eta s_1t_1=s_2f_1t_1=s_2t_2f_1=t_2f_1\qquad\mathrm{et}\qquad 
t_2\eta t_1=t_2\eta s_1t_1=s_2f_2t_1=s_2t_2f_2=t_2f_2$$
La même démonstration marche dans l'autre sens.\QEDb
Ce lemme implique que les applications $\left\langle f_2;\eta s_1\right\rangle$ et 
$\left\langle\eta t_1;f_1\right\rangle$ existent puisque c'est équivalent à dire
$$s_2f_2=t_2(\eta s_1)\qquad\mathrm{et}\qquad s_2(\eta t_1)=t_2f_2$$
\begin{rema}\label{rema:cat:3.1} 
Si nous observons cette preuve, nous avons juste utilisé la condition (NAT 1).  
En fait (NAT 2) est utile pour assurer une certaine forme d'unicité qui sera 
nécessaire plus tard pour définir une 
composition des transformations naturelles (voir la proposition \ref{prop:cat:3.2}).
\end{rema}
\begin{prop}\label{prop:cat:3.1} 
Si $f:\BBB_1\to\BBB_2$ est un foncteur, alors $fs_1$ est une transformation 
naturelle de $f$ dans lui-même.
\end{prop}
\Dem
\begin{description}
\item[(NAT 1)] $s_2\left(fs_1\right)s_1=s_2s_2fs_1=s_2\left(fs_1\right)
\qquad\mathrm{et}\qquad t_2\left(fs_1\right)s_1=t_2s_2fs_1=s_2\left(fs_1\right)$
\item[(NAT 2)] $\left(fs_1\right)s_1=fs_1$
\item[(NAT 3)] $\#_2\left\langle f;\left(fs_1\right)s_1\right\rangle
=\#_2\left\langle f;fs_1\right\rangle
=\#_2\left\langle\Id_{\BBB_2}f;s_2f\right\rangle
=\#_2\left\langle \Id_{\BBB_2};s_2\right\rangle f
=\Id_{\BBB_2}f=f$\\ et\quad
$\#_2\left\langle\left(fs_1\right)t_1;f\right\rangle
=\#_2\left\langle ft_1;f\right\rangle
=\#_2\left\langle t_2f;\Id_{\BBB_2}f\right\rangle
=\#_2\left\langle t_2;\Id_{\BBB_2}\right\rangle f
=\Id_{\BBB_2}f=f$ \QED
\end{description}
\begin{prop}\label{prop:cat:3.2} 
Si
$$\eta_1:\left(f_1:\left(\BBB_1;s_1;t_1;\#_1\right)\to
\left(\BBB_2;s_2;t_2;\#_2\right)\right)
\leadsto\left(f_2:\left(\BBB_1;s_1;t_1;\#_1\right)
\to\left(\BBB_2;s_2;t_2;\#_2\right)\right)$$
et
$$\eta_2:\left(f_2:\left(\BBB_1;s_1;t_1;\#_1\right)\to
\left(\BBB_2;s_2;t_2;\#_2\right)\right)
\leadsto\left(f_3:\left(\BBB_1;s_1;t_1;\#_1\right)\to
\left(\BBB_2;s_2;t_2;\#_2\right)\right)$$
sont deux transformations naturelles, alors l'application d'ensembles 
$\left\langle \eta_2;\eta_1\right\rangle$ existe et l'application 
$\#_2\left\langle \eta_2;\eta_1\right\rangle$ définit une transformation
naturelle de $f_1$ dans $f_3$.
\end{prop}
\Dem
\begin{itemize}
\item Pour démontrer l'existence de $\left\langle \eta_2;\eta_1\right\rangle$, 
il suffit de prouver que $s_2\eta_2=t_2\eta_1$. Ce qui est fait ci-
dessous en utilisant l'axiome (NAT 2).
$$s_2\eta_2=s_2\eta_2s_1=s_2f_2 \qquad\mathrm{et}
\qquad t_2\eta_1=t_2\eta_1s_1=s_2f_2$$
\item Il nous reste à prouver que $\#_2\left\langle \eta_2;\eta_1\right\rangle$ est une application
naturelle de $f_1$ vers 
\begin{description}
\item[(NAT 1)] $s_2\#_2\left\langle \eta_2;\eta_1\right\rangle s_1
=s_2\pi_{t_2}^{s_2}\left\langle \eta_2;\eta_1\right\rangle s_1
=s_2\eta_1s_1=s_2f_1$\\
et $t_2\#_2\left\langle \eta_2;\eta_1\right\rangle s_1
=t_2\pi_{s_2}^{t_2}\left\langle \eta_2;\eta_1\right\rangle s_1
=t_2\eta_2s_1=s_2f_3$
\item[(NAT 2)] $\#_2\left\langle \eta_2;\eta_1\right\rangle s_1
=\#_2\left\langle \eta_2s_1;\eta_1s_1\right\rangle
=\#_2\left\langle \eta_2;\eta_1\right\rangle$
\item[(NAT 3)] Conséquence des calculs suivants
\begin{align*}
\#_2\left\langle f_3;\#_2\left\langle \eta_2;\eta_1
\right\rangle s_1\right\rangle
&=\#_2\left\langle f_3;\#_2\left\langle\eta_2s_1;\eta_1s_1
\right\rangle\right\rangle\displaybreak[0]\\
&=\#_2\left(\Id_{\BBB_2}\times_c\#_2\right)\left\langle 
f_3;\left\langle\eta_2s_1;\eta_1s_1\right\rangle\right\rangle
\displaybreak[0]\\
&=\#_2\left(\#_2\times_c\Id_{\BBB_2}\right)\omega^{-1}
\left\langle f_3;
\left\langle\eta_2s_1;\eta_1s_1\right\rangle\right\rangle
\displaybreak[0]\\
&=\#_2\left(\#_2\times_c\Id_{\BBB_2}\right)
\left\langle\left\langle f_3;\eta_2s_1\right\rangle;
\eta_1s_1\right\rangle
\displaybreak[0]\\
&=\#_2\left\langle\#_2\left\langle f_3;\eta_2s_1\right\rangle;
\eta_1s_1\right\rangle
\displaybreak[0]\\
&=\#_2\left\langle\#_2\left\langle\eta_2t_1;f_2\right\rangle;
\eta_1s_1\right\rangle\displaybreak[0]\\
&=\#_2\left(\#_2\times_c\Id_{\BBB_2}\right)
\left\langle\left\langle\eta_2t_1;f_2\right\rangle;
\eta_1s_1\right\rangle\displaybreak[0]\\
&=\#_2\left(\Id_{\BBB_2}\times_c\#_2\right)\omega
\left\langle\left\langle\eta_2t_1;f_2\right\rangle;
\eta_1s_1\right\rangle\displaybreak[0]\\
&=\#_2\left(\Id_{\BBB_2}\times_c\#_2\right)
\left\langle\eta_2t_1;\left\langle f_2;
\eta_1s_1\right\rangle\right\rangle\displaybreak[0]\\
&=\#_2\left\langle\eta_2t_1;\#_2\left\langle f_2;
\eta_1s_1\right\rangle\right\rangle\displaybreak[0]\\
&=\#_2\left\langle\eta_2t_1;\#_2\left\langle\eta_1t_1;
f_1\right\rangle\right\rangle\displaybreak[0]\\
&=\#_2\left(\Id_{\BBB_2}\times_c\#_2\right)
\left\langle\eta_2t_1;\left\langle\eta_1t_1;
f_1\right\rangle\right\rangle\displaybreak[0]\\
&=\#_2\left(\#_2\times_c\Id_{\BBB_2}\right)\omega^{-1}
\left\langle\eta_2t_1;\left\langle\eta_1t_1;
f_1\right\rangle\right\rangle\displaybreak[0]\\
&=\#_2\left(\#_2\times_c\Id_{\BBB_2}\right)
\left\langle\left\langle\eta_2t_1;\eta_1t_1\right\rangle;
f_1\right\rangle\displaybreak[0]\\
&=\#_2\left\langle\#_2\left\langle\eta_2t_1;\eta_1t_1
\right\rangle; f_1\right\rangle\displaybreak[0]\\
&=\#_2\left\langle\#_2\left\langle\eta_2;\eta_1
\right\rangle t_1; f_1\right\rangle
\end{align*}\QED
\end{description}
\end{itemize} 
\begin{lemm}\label{lemm:cat:3.1}
Une \emph{petite transformation naturelle}\index{transformation naturelle!
petite transformation naturelle} (transformation naturelle entre petits foncteurs)
est un élément de l'univers.
\end{lemm}
\Dem
La démonstration est semblable à celle du lemme \ref{lemm:cat:1.1}.  
\QEDb
On peut donc considérer
\begin{itemize}
\item $\NAT$ l'ensemble des  \emph{petites transformations naturelles},
\item $s:\NAT\to\NAT,\ s\left(\eta:f_1\leadsto f_2\right)=f_1s_1$\\
$t:\NAT\to\NAT,\ t\left(\eta:f_1\leadsto f_2\right)=f_2s_1$
\item $\star:\NAT\times_c\NAT\to\NAT,\ 
\star\left(\eta_2,\eta_1\right)=\eta_2\star\eta_1
=\#_2\left\langle\eta_2;\eta_1\right\rangle$\\
où $\NAT\times_c\NAT=\left\lbrace\left(\eta_2;\eta_1\right)\in\NAT\times\NAT\mid
s\eta_2=t\eta_1\right\rbrace $.
\end{itemize} 
\begin{theo}\label{theo:cat:3.1} 
$\left(\NAT;s;t;\star\right)$ est une catégorie.
\end{theo}
\Dem
\begin{description}
\item[(CAT 1)] $ss\left(\eta\right)=s\left(f_1s_1\right)=f_1s_1=s\left(\eta\right),\qquad
st\left(\eta\right)=s\left(f_2s_1\right)=f_2s_1=t\left(\eta\right)$\\
$tt\left(\eta\right)=t\left(f_2s_1\right)=f_2s_1=t\left(\eta\right),\qquad
ts\left(\eta\right)=s\left(f_1s_1\right)=f_1s_1=s\left(\eta\right)$
\item[(CAT 2)] $s\star\left(\eta_2;\eta_1\right)=s\#_2 \left\langle\eta_2;\eta_1\right\rangle
=f_1s_1=s\left(\eta_1\right)=s\pi_{t}^{s}\left(\eta_2;\eta_1\right)$\\
$t\star\left(\eta_2;\eta_1\right)=t\#_2 \left\langle\eta_2;\eta_1\right\rangle
=f_3s_1=t\left(\eta_2\right)=t\pi_{s}^{t}\left(\eta_2;\eta_1\right)$
\item[(CAT 3)] On a
\begin{align*}
\star\left\langle\Id_{\BBB_2};s\right\rangle\left(\eta\right)
=\star\left(\eta;s\eta\right)
&=\#_2\left\langle\eta;s\eta\right\rangle\displaybreak[0]\\
&=\#_2\left\langle\eta;f_1s_1\right\rangle\displaybreak[0]\\
&=\#_2\left\langle\eta;s_2f_1\right\rangle\displaybreak[0]\\
&=\#_2\left\langle\eta;s_2\eta s_1\right\rangle
=\#_2\left\langle\eta;s_2\eta\right\rangle
=\#_2\left\langle\Id_{\BBB_2};s_2\right\rangle\left(\eta\right)
=\eta
\end{align*}
et de même $\star\left\langle t;\Id_{\BBB_2}\right\rangle
\left(\eta\right)=\eta$.
\item[(CAT 4)] Finalement
\begin{align*}
\star\left(\star\times_c\Id_{\NAT}\right)
\left(\left(\eta_3;\eta_2\right);\eta_1\right)
&=\star\left(\star\left(\eta_3;\eta_2\right);\eta_1\right)
\displaybreak[0]\\
&=\#_2\left\langle\#_2\left\langle\eta_3;\eta_2\right\rangle;
\eta_1\right\rangle\displaybreak[0]\\
&=\#_2\left(\#_2\times_c\Id_{\BBB_2}\right)
\left\langle\left\langle\eta_3;\eta_2\right\rangle;\eta_1\right\rangle
\displaybreak[0]\\
&=\#_2\left(\Id_{\BBB_2}\times_c\#_2\right)\omega
\left\langle\left\langle\eta_3;\eta_2\right\rangle;\eta_1\right\rangle
\displaybreak[0]\\
&=\#_2\left(\Id_{\BBB_2}\times_c\#_2\right)
\left\langle\eta_3;\left\langle\eta_2;\eta_1\right\rangle\right\rangle
\displaybreak[0]\\
&=\#_2\left\langle\eta_3;\#_2\left\langle\eta_2;\eta_1\right\rangle
\right\rangle\displaybreak[0]\\
&=\star\left(\eta_3;\star\left(\eta_2;\eta_1\right)\right)
\displaybreak[0]\\
&=\star\left(\Id_{\NAT}\times_c\star\right)
\left(\eta_3;\left(\eta_2;\eta_1\right)\right)\displaybreak[0]\\
&=\star\left(\Id_{\NAT}\times_c\star\right)\omega
\left(\left(\eta_3;\eta_2\right);\eta_1\right)
\end{align*}\QED
\end{description}
La proposition suivante montre l'un des intérêts du point de vue que nous avons 
adopté. 
\begin{prop}\label{prop:cat:3.3} 
Soient $\eta:\left(f_1:\BBB_1\to\BBB_2\right)\leadsto\left(f_2:\BBB_1
\to\BBB_2\right)$ une application naturelle et $g:\left(\BBA;s;t;\#\right)
\to\left(\BBB_1;s_1;t_1;\#_1\right)$, $h:\left(\BBB_2;s_2;t_2;\#_2\right) 
\to\left(\BBC;s';t';\#'\right) $ deux foncteurs. L'application d'ensembles 
$h\eta g:\BBA\to\BBC$ est en fait une transformation naturelle du foncteur 
$h\bullet f_1\bullet g$ vers le foncteur $h\bullet f_2\bullet g$.
\end{prop}
\Dem Vérifions les trois axiomes
\begin{description}
\item[(NAT 1)] On a
$$s'\left(h\eta g\right)s=hs_2\eta s_1g=hs_2f_1g=s'\left(h
\bullet f_1\bullet g\right)$$
$$t'\left(h\eta g\right)s=ht_2\eta s_1g=hs_2f_2g=s'\left(h\bullet f_2
\bullet g\right)$$
\item[(NAT 2)] $\left(h\eta g\right)s=h\eta s_1g=h\eta g$
\item[(NAT 3)] On a
\begin{align*}
\#'\left\langle h\bullet f_2\bullet g;h\eta g s\right\rangle
&=\#'\left\langle h f_2 g;h\eta s_1g\right\rangle\displaybreak[0]\\
&=\#'\left(h\times_c h\right)\left\langle f_2;\eta s_1\right\rangle
g\displaybreak[0]\\
&=h\#_2\left\langle f_2;\eta s_1\right\rangle g\displaybreak[0]\\
&=h\#_2\left\langle \eta t_1;f_1\right\rangle g\displaybreak[0]\\
&=\#'\left(h\times_c h\right)\left\langle \eta t_1;f_1\right\rangle
g\displaybreak[0]\\
&=\#'\left\langle h\eta t_1g;hf_1g\right\rangle\displaybreak[0]\\
&=\#'\left\langle h\eta gt;hf_1g\right\rangle\displaybreak[0]\\
&=\#'\left\langle h\eta gt;h\bullet f_1\bullet g\right\rangle
\end{align*}\QED
\end{description} 
%%%%%%%%%%%%%%%%%
%%% section 4 %%%
%%%%%%%%%%%%%%%%%
\section{Relation avec les notions habituelles}\label{sec:cat:4} 
\subsection{Catégorie \Dire{classique} associée à une catégorie 
\Dire{sans objets}}\label{subsec:cat:4.1} 
Soit $\left(\BBB;s;t\right)$ une catégorie \Dire{sans objet}. Parmi les 
éléments de $\BBB$, certains jouent un rôle particulier : les identités.

Par définition, $x\in\BBB$ est une identité si pour tous éléments $y$ et 
$z$ de $\BBB$ tels que $s\left(y\right)=t\left(x\right)$ et 
$s\left(x\right)=t\left(z\right)$, on a $y\#x=y$ et $x\#z=z$. Notons 
$\Id\left(\BBB\right)$ l'ensemble des identités de $\BBB$.

Bien que nous n'en ayons besoin que plus tard, nous pouvons maintenant 
définir la notion d'éléments inversible. On dit que $x\in\BBB$ est inversible 
quand il existe $y\in\BBB$ tel que $x\#y\in\Id\left(\BBB\right)$ et 
$y\#x\in\Id\left(\BBB\right)$.
\begin{lemm}\label{lemm:cat:4.1} 
Si pour toute application $f:\BBB\to\BBB$, on pose
$$\Ima\left(f\right)=\left\{x\in\BBB\mid\exists y\in\BBB,f\left(y\right)=x\right\}
\qquad\mathrm{et}\qquad
\FIX\left(f\right)=\left\{x\in\BBB\mid f\left(x\right)=x\right\}$$
alors on a
$$\Id\left(\BBB\right)=\Ima\left(s\right)=\Ima\left(t\right)=
\FIX\left(s\right)=\FIX\left(t\right)$$
\end{lemm}
\Dem Comme $st=t$ et $ts=s$, on a $\Ima\left(t\right)\subseteq\Ima\left(s\right)$ et
$\Ima\left(s\right)\subseteq\Ima\left(t\right)$.\\
Il est clair que $\FIX\left(s\right)\subseteq\Ima\left(s\right)$ et $\FIX\left(t\right)\subseteq\Ima\left(t\right)$.\\
Soit $x\in\Ima\left(s\right)$. On peut écrire $x=s\left(y\right)$ avec $y\in\BBB$. On a
ainsi $s\left(x\right)=ss\left(y\right)=s\left(y\right)=x$.
Donc $\Ima\left(s\right)\subseteq\FIX\left(s\right)$ et de même 
$\Ima\left(t\right)\subseteq\FIX\left(t\right)$.\\
Il ne nous reste plus qu'à prouver $\Id\left(\BBB\right)=\FIX\left(s\right)$.\\
Si $x\in\Id\left(\BBB\right)$, alors $x\#s\left(x\right)=s\left(x\right)$. Or d'après 
l'axiome (CAT 3), on a $x\#s\left(x\right)=x$. D'où $s\left(x\right)=x$. C'est-à-dire 
$\Id\left(\BBB\right)\subseteq\FIX\left(s\right)$.\\
De plus on vient de voir que si $x$ appartient à $\FIX\left(s\right)$, alors
$x=s\left(x\right)=t\left(x\right)$. Donc d'après l'axiome (CAT 3), $x$ est une identité.
\QEDb
Nous pouvons maintenant définir une catégorie au sens classique du terme en 
prenant les données suivantes :
\begin{itemize}
\item $\Ob=\Id\left(\BBB\right)$,
\item Pour tout $\left(x;y\right)\in\Ob^2$, 
$$\Mor_{\BBB}\left(x;y\right)=\left\lbrace f\in\BBB\mid s\left(f\right)
=x,\ t\left(f\right)=y\right\rbrace$$
\item Pour chaque $x\in\Ob$, 
$$\Id_x=x\in\Mor_{\BBB}\left(x;x\right)$$
\item Pour tout triplet $\left(x;y;z\right)\in\Ob^3$,
\begin{align*}
\circ :\Mor_{\BBB}\left(y;z\right)\times\Mor_{\BBB}
\left(x;y\right)&\longrightarrow\Mor_{\BBB}\left(x;z\right)\\
\left(f;g\right)&\longmapsto f\circ g=f\# g
\end{align*}
Cette dernière est bien définie car $s\left(f\right)=x=
t\left(g\right)$
\end{itemize}
En effet les axiomes sont satisfaits
\begin{description}
\item[(identité)] Si $f:x\to y$ pour tout $\left(x;y\right)\in\Ob^2$,
alors $$f\circ\Id_x=f\#x=f\quad\mathrm{et}\quad\Id_y\circ f=y\# f=f$$
d'après la définition même de $\Id\left(\BBB\right)$.
\item[(associativité)] Si $f:z\to y$, $g:y\to x$ et $h:x\to v$ sont 
trois morphismes quelconques, alors
\begin{align*}
\left(f\circ g\right)\circ h
&=\#\left(\#\times_c\Id_{\BBB}\right)\left(\left(f;g\right);h\right)
\displaybreak[0]\\ 
&=\#\left(\Id_{\BBB}\times_c\#\right)\omega\left(\left(f;g\right);
h\right)\displaybreak[0]\\ 
&=\#\left(\Id_{\BBB}\times_c\#\right)\left(f;\left(g;h\right)\right)
\displaybreak[0]\\ 
&=f\circ\left(g\circ h\right)
\end{align*}
\end{description}
\subsection{Catégorie \Dire{sans objets} associée à une catégorie \Dire{classique}}
\label{subsec:cat:4.2} 
Soit $\BBC$ une catégorie classique dont la composition est notée $\circ$.
On définit une catégorie \Dire{sans objets} $\left(\BBB;s;t;\#\right)$ en prenant
\begin{itemize}
\item $\BBB=\bigcup_{(x;y)\in\Ob(\BBC)^2}\Mor_{\BBC}\left(x;y\right)$
\item $s:\BBB\to\BBB,\ s\left(f:x\to y\right)=\Id_x$\\
$t:\BBB\to\BBB,\ t\left(f:x\to y\right)=\Id_y$
\item $\#:\BBB\times_c\BBB\to\BBB,\ \#\left(f;g\right)=f\circ g$\\
où $f\circ g$ est bien définie car
$\Id_y=s\left(f:y\to z\right)=t\left(g:x\to y'\right)=\Id_{y'}$, i.e.
$y=y'$.
\end{itemize}
Les axiomes sont satisfaits
\begin{description}
\item[(CAT 1)] \'Evident.
\item[(CAT 2)] Si $\left(f:z\to y;g:x\to y\right)\in\BBB\times_c\BBB$ alors
$$s\#\left(f;g\right)=s\left(f\circ g\right)=x\quad\mathrm{et}\quad 
s\pi_{s}^{t}\left(f;g\right)=s\left(g\right)=x$$
$$t\#\left(f;g\right)=t\left(f\circ g\right)=z\quad\mathrm{et}\quad 
t\pi_{t}^{s}\left(f;g\right)=t\left(f\right)=z$$
\item[(CAT 3)] Soit $f:x\to y\in\BBB$.
$$\OP\left\langle\Id_{\BBB};s\right\rangle(f)=\OP(f;\Id_x)
=f\circ\Id_x=f=\Id_{\BBB}(f)$$
$$\OP\left\langle t;\Id_{\BBB}\right\rangle(f)=\OP(\Id_y;f)
=\Id_y\circ f=f=\Id_{\BBB}(f)$$
\item[(CAT 4)] Soit $\left(\left(f;g\right);h\right)\in\left( 
\BBB\times_c\BBB\right)\times_c\BBB$. On a
$$\#\left(\#\times_c\Id\right)\left(\left(f;g\right);h\right)
=\left(f\circ g\right)\circ h$$
et
$$\#\left(\Id\times_c\#\right)\omega\left(\left(f;g\right);h\right)
=\#\left(\Id\times_c\#\right)\left(f;\left(g;h\right)\right)
=f\circ\left(g\circ h\right)$$
or $\left(f\circ g\right)\circ h=f\circ\left(g\circ h\right)$.
\end{description}
\subsection{Les foncteurs}
Soit $f:\left(\BBB_1;s_1;t_1;\#_1\right)\to\left(\BBB_2;s_2;t_2;\#_2\right)$ 
un foncteur de catégories \Dire{sans objets}.
Les données suivantes
\begin{itemize}
\item pour chaque $x\in\Ob_1=\Id\left(\BBB_1\right)$, 
$$F\left(x\right)=f\left(x\right)\in\Ob_2=\Id\left(\BBB_2
\right)$$
C'est bien défini car $f\left(x\right)=f\left(s_1\left(s\left(
x\right)\right)\right)=s_2f\left(x\right)\in\Ima\left(s_2\right)=\Ob_2$.
\item pour chaque $u\in\Mor_{\BBB_1}\left(x;y\right)$,
$$F\left(u\right)=f\left(u\right)\in\Mor_{\BBB_2}
\left(F(x);F(y)\right)$$
car $s_2\left(f\left(u\right)\right)=f\left(s_1\left(u\right)\right)
=f\left(x\right)$ et $t_2\left(f\left(u\right)\right)=
f\left(t_1\left(u\right)\right)=f\left(y\right)$.
\end{itemize}
définissent bien un foncteur car
\begin{itemize}
\item pour tout $x\in\Ob_1$,
$$F\left(\Id_x\right)=f\left(x\right)=\Id_{f(x)}=
\Id_{F(x)}$$
\item et pour tous morphismes $u:y\to z$, $v:x\to y$,
$$F\left(u\circ_1 v\right)=f\left(u\#_1 v\right)
=f\left(u\right)\#_2 f\left(v\right)=F\left(u\right)
\circ_2 F\left(v\right)$$
\end{itemize}

Soit $F:\BBC_1\to\BBC_2$ un foncteur de catégories \Dire{classiques}.
La fonction d'ensembles $f:\BBB_1\to\BBB_2$ donnée par
$$f\left(u\right)=F\left(u\right)\in\Mor_{\BBC_1}\left(F(x);F(y)\right)
\subseteq\BBB_2$$pour tout $u\in\Mor_{\BBC_1}\left(x;y\right)$, est un 
foncteur de catégories \Dire{sans objets}.
\begin{description}
\item[(FONC 1)]pour $u\in\Mor_{\BBC_1}\left(x;y\right)$, $$f\left(s_1\left(u\right)\right)=f\left(\Id_x\right)=F\left(\Id_x\right)
=\Id_{F(x)}=s_2\left(F\left(u\right)\right)=s_2f\left(u\right)$$
De même $f\left(t_1\left(u\right)\right)=t_2\left(f\left(u\right)\right)$.
\item[(FONC 2)]pour $\left(u;v\right)\in\Mor_{\BBC_1}\left(y;z\right)
\times\Mor_{\BBC_1}\left(x;y\right)$,
$$\#_2\left(f\times_c f\right)\left(u;v\right)=
\#_2\left(f(u);f(v)\right)=F(u)\circ_2F(v)=F\left(u\circ_1v\right)
=f\#_1\left(u;v\right)$$
\end{description} 

\subsection{Les transformations naturelles}\label{subsec:cat:4.3} 
Les calculs étant toujours les mêmes, nous nous contenterons ici de donner
les constructions.

Soit $\eta:\left(f_1:\BBB_1\to\BBB_2\right)\leadsto
\left(f_1:\BBB_1\to\BBB_2\right)$ une transformation naturelle (de catégories 
\Dire{sans objets}). On obtient une transformation naturelle $\Xi:F_1\leadsto F_2$ 
en prenant pour tout objet $x\in\Ob_1$
$$\Xi_{x}=\eta\left(x\right)\in\Mor_{\BBC_2}\left(F_1(x);F_2(x)\right)$$

Et si $\Xi:\left(F_1:\BBC_1\to\BBC_2\right)\leadsto
\left(F_2:\BBC_1\to\BBC_2\right)$ est une transformation naturelle (de 
catégories \Dire{classiques}), alors la fonction d'ensembles
$\eta:f_1\to f_2$ donnée par
$$\eta\left(u\right)=\Xi\left(x\right)\in\Mor_{\BBC_2}\left(F_1(x);F_2(x)
\right)\subseteq\BBB_2$$
pour $u\in\Mor_{\BBC_1}\left(x;y\right)\subseteq\BBB_1$, est une transformation
naturelle.

%%%%%%%%%%%%%%%%%%%%%%%%%%%%%%%%%%%%%%%%%%%%%%%%%%
\chapter{Les groupements}\label{chap:gr}
%%%%%%%%%%%%
%%%  section 1 %%%
%%%%%%%%%%%%
\section{Les origines}\label{sec:gr:1}
En lisant le chapitre précédent, une première évidence s'impose : L'utilisation 
du pull-back complique la théorie. Quitte à perdre sur l'unicité, rien ne nous empêche 
d'étendre la composition $\OP$ à l'ensemble produit $\BBB\times\BBB$.
Ce qui nous donne les définitions suivantes pour les catégories, foncteurs et 
transformations naturelles :
\begin{itemize}
\item Une catégorie $\left(\BBB;s;t;\OP\right)$ est un ensemble $\BBB$ muni de trois applications
\begin{itemize}
\item $s:\BBB\to\BBB$,
\item $t:\BBB\to\BBB$,
\item et $\OP:\BBB\times\BBB\to\BBB$,
\end{itemize}    
telles que les axiomes suivants soient satisfaits :
\begin{description}
\item[(CAT 1)] $ss=s$,\quad$st=t$,\quad$tt=t$,\quad$ts=s$ ;
\item[(CAT 2)] Si $x$ et $y$ sont deux éléments de $\BBB$ satisfaisant 
        $s(x)=t(y)$ alors $$s\left(x\OP y\right)=s(y)\quad \text{et}\quad 
        t\left(x\OP y\right)=t(x)$$
\item[(CAT 3)] pour tout $x\in\BBB$ alors $x\OP s(x)=x$ et $t(x)\OP x=x$ ;
\item[(CAT 4)] si $x,y,z$ sont trois éléments de $\BBB$ tels que $s(x)=t(y)$ et
        $s(y)=t(z)$, alors $$\left(x\OP y\right)\OP z=x\OP\left(y\OP z\right)$$
\end{description}
\item Un foncteur $f$ entre la catégorie $\left(\BBB_1;s_1;t_1;\OP_1\right)$ et la 
catégorie $\left(\BBB_2;s_2;t_2;\OP_2\right)$ est un triplet $\left((\BBB_1;s_1;t_1;
\OP_1);(\BBB_2;s_2;t_2;\OP_2);f\right)$, souvent noté $f:\left(\BBB_1;s_1;t_1;\OP_1\right)
\longrightarrow\left(\BBB_2;s_2;t_2;\OP_2\right)$, où $f$ est 
une application d'ensembles $f:\BBB_1\to\BBB_2$ vérifiant les axiomes
\begin{description}
\item[(FONC 1)] $fs_1=s_2f$,\quad $ft_1=t_2f$ ;
\item[(FONC 2)] si $x$ et $y$ sont deux éléments de $\BBB_1$ tels que 
        $s(x)=t(y)$, alors $$f\left(x\OP_1y\right)=f(x)\OP_2f(y)$$
\end{description}
\item Une application naturelle $\eta:f_1\leadsto f_2$ entre le foncteur 
$f_1$ et le foncteur $f_2$ est un triplet $\left(\eta;f_1;f_2\right)$ où 
$\eta:\BBB_1\to\BBB_2$ une application d'ensembles satisfaisant les trois axiomes 
ci-dessous
\begin{description}
\item[(NAT 1)] $s_2\eta s_1=f_1s_1$,\quad $t_2\eta s_1=f_2s_1$ ;
\item[(NAT 2)] $\eta=\eta s_1$.
\item[(NAT 3)] pour tout $x$ dans $\BBB_1$, on a\quad
        $f_2(x)\OP_2\eta\left(s_1(x)\right)=\eta\left(t_1(x)\right)\OP_2f_1(x)$ ;
\end{description}
\end{itemize}
D'après la remarque \ref{rema:cat:3.1}, l'axiome (NAT 2) ne semble par très naturelle. 
Pour être plus précis il n'a été utilisé que pour assurer l'existence de la composition 
$\OP\left\langle\eta_1;\eta_2\right\rangle$ et assurer la véracité de l'axiome (CAT 3) 
pour la catégorie $\NAT$.

En fait comme nous l'avons déjà dit, l'axiome (NAT 2) permet de réduire le nombre de 
transformations naturelles en identifiant celles qui nous semblent avoir les mêmes 
propriétés. Un moment de réflexion, nous amène à penser que l'axiome (CAT 3) des catégories joue un
rôle très semblable. Il nous sert à identifier les objets et les identités.
C'est-à-dire à identifier les catégories qui ont les mêmes identités. C'est 
la supression de ces deux axiomes qui nous pousse à introduire la notion de groupement. 
%%%%%%%%%%%%
%%%  section 2  %%%
%%%%%%%%%%%%
\section{Définition d'un groupement}\label{sec:gr:2}
On définit un \emph{groupement}\index{groupement} $\left(\BBB;s;t;\OP\right)$
comme étant un ensemble $\BBB$ muni de trois applications
\begin{itemize}
\item $s:\BBB\to\BBB$, la \emph{source}\index{source},
\item $t:\BBB\to\BBB$, le \emph{but}\index{but},
\item et $\OP:\BBB\times\BBB\to\BBB$, la \emph{composition}\index{composition},
\end{itemize}    
telles que les axiomes suivants soient satisfaits :
\begin{description}
\item[(GR 1)] $ss=s$,\quad$st=t$,\quad$tt=t$,\quad$ts=s$ ;
\item[(GR 2)] Si $x$ et $y$ sont deux éléments de $\BBB$ satisfaisant 
        $s(x)=t(y)$ alors $$s\left(x\OP y\right)=s(y)\quad \text{et}\quad 
        t\left(x\OP y\right)=t(x)$$
\item[(GR 3)] si $x,y,z$ sont trois éléments de $\BBB$ tels que $s(x)=t(y)$ et
        $s(y)=t(z)$, alors $$\left(x\OP y\right)\OP z=x\left(y\OP z\right)$$
\end{description}
\begin{exem}\label{exem:gr:2.1}
\'Etant donné que notre objectif est de généraliser légèrement la théorie de catégories, il
est normal que celle-ci soient des groupements. Les axiomes (CAT 1), (CAT 2) et (CAT 4) sont
mot pour mot les trois axiomes (GR1), (GR2) et (GR3).
\end{exem}
\begin{exem}\label{exem:gr:2.2}
Soit $(M;\bullet)$ un monoïde non vide. Choisissons un élément $c$ quelconque dans 
$M$ et notons $\breve{c}$ l'application constante de $M$ dans lui-même qui a tout 
élément associe $c$. Il est immédiat que $\left(M;\breve{c};\breve{c};\bullet\right)$
est un groupement. Il est aussi à remarquer que l'on peut ainsi associer au monoïde 
$(M;\bullet)$ un grand nombre de groupements.   
\end{exem}
\begin{rema}\label{rema:gr:2.1}
Dans l'exemple précédent, si le monoïde est un groupe $(G;\bullet;e)$ d'élément 
neutre $e$, alors on peut prendre $e$ à la place de $c$. Clairement 
$(G;\breve{e};\breve{e};\bullet)$ 
est un groupement. De plus on vérifie aisément que c'est une catégorie. 
C'est d'ailleurs de cette manière que l'on montre habituellement que les groupes 
sont des exemples de catégorie.
\end{rema}
Un groupement est dit \emph{petit}\index{groupement!petit groupement} si $\BBB$ 
est un petit ensemble.
\begin{rema}\label{rema:gr:2.2} 
Pour être plus précis nous devrions parler de \emph{$\MFU$-groupements}
\index{groupement!$\MFU$-groupement} et de 
\emph{petits $\MFU$-groupements}\index{groupement!petit $\MFU$-groupement}.
\end{rema} 
\begin{lemm}\label{lemm:gr:2.1}
Un petit $\MFU$-groupement est un élément de l'univers $\MFU$.
\end{lemm}
\Dem 
La démonstration est semblable à celle du lemme \ref{lemm:cat:1.1}. 
\QEDb
Il est à remarquer dans la définition des groupements, qu'il y a une certaine symétrie 
entre le rôle de l'application source et celui de l'application but. La seule différence
se situe au niveau de la composition. Soyons plus précis en prenant un groupement 
$\left(\BBB;s;t;\OP\right)$. Il est aisé de montrer que le quadruplet $\left(\BBB^\ast;
s^\ast;t^\ast;\OP^\ast\right)$, où $\BBB^\ast=\BBB$, $s^\ast=t$, $t^\ast=s$ et 
$x\OP^\ast y=y\OP x$ pour tous $x$ et $y$ dans $\BBB^\ast$, est un groupement.
\begin{description}
\item[(GR 1)] On a $s^\ast s^\ast=tt=t=s^\ast $. De même \quad 
$s^\ast t^\ast=t^\ast $,\quad$t^\ast t^\ast=t^\ast $,\quad$t^\ast s^\ast=s^\ast $.
\item[(GR 2)] Si $x$ et $y$ sont deux éléments de $\BBB^\ast=\BBB$ satisfaisant 
        $s^\ast(x)=t^\ast(y)$ alors $$s^\ast\left(x\OP^\ast y\right)=t(y\OP x)=t(y)=
        s^\ast(y)\quad\text{et}\quad t^\ast\left(x\OP^\ast y\right)=s(y\OP x)=s(x)=t^\ast(x)$$
        car $t(x)=s(y)$.
\item[(GR 3)] si $x,y,z$ sont trois éléments de $\BBB^\ast=\BBB$ tels que 
$s^\ast(x)=t^\ast(y)$ et
        $s^\ast(y)=t^\ast(z)$, alors $$\left(x\OP^\ast y\right)\OP^\ast z=
        z\OP\left(y\OP x\right)=\left(z\OP y\right)\OP x=x\OP^\ast\left(
        y\OP^\ast \right)$$
        car $t(x)=s(y)$ et $t(y)=s(z)$.
\end{description}
Le groupement $\left(\BBB^\ast;s^\ast;t^\ast;\OP^\ast\right)$ est appelé \emph{dual}
\index{groupement!groupement dual} du groupement $\left(\BBB;s;t;\OP\right)$.
%%%%%%%%%%%%%
%%%  section 3  %%%%
%%%%%%%%%%%%%
\section{Les g-foncteurs et g-morphismes}\label{sec:gr:3}
Un \emph{g-foncteur}\index{g-foncteur} $f$ entre le groupement
$\left(\BBB_1;s_1;t_1;\OP_1\right)$
et le groupement $\left(\BBB_2;s_2;t_2;\OP_2\right)$ est un triplet
$\left((\BBB_1;s_1;t_1;
\OP_1);(\BBB_2;s_2;t_2;\OP_2);f\right)$, souvent noté
$f:\left(\BBB_1;s_1;t_1;\OP_1\right)
\longrightarrow\left(\BBB_2;s_2;t_2;\OP_2\right)$, où $f$ est
une application d'ensembles $f:\BBB_1\to\BBB_2$ vérifiant les axiomes suivants
\begin{description}
\item[(GFONC 1)] $fs_1=s_2f$,\quad $ft_1=t_2f$ ;
\item[(GFONC 2)] si $x$ et $y$ sont deux éléments de $\BBB_1$ tels que
        $s_1(x)=t_1(y)$, alors $$f\left(x\OP_1y\right)=f(x)\OP_2f(y)$$
\end{description}
Par abus de notation, nous noterons souvent $$f:\BBB_1\to\BBB_2$$ à la place de
$$f:\left(\BBB_1;s_1;t_1;\OP_1\right)\longrightarrow\left(\BBB_2;s_2;t_2;
\OP_2\right)$$
\begin{rema}\label{rema:gr:3.1}
Comme pour les groupements, on devrait parler de \emph{$\MFU$-g-foncteurs}
\index{g-foncteur!$\MFU$-g-foncteur}, plutôt 
que de g-foncteurs.
\end{rema}
\begin{exem}\label{exem:gr:3.1}
Il est clair que tout foncteur entre catégories est un g-foncteur.
\end{exem}
\begin{rema}\label{rema:gr:3.2}
Quand nous avons associé des groupements à un monoïde, nous avons choisi des éléments
de celui-ci, ce qui fait que les morphismes de monoïdes ne deviennent pas nécessairement
des g-foncteurs. Nous sommes ici dans la même situation que celle où se sont trouvés les
topologues quand ils ont commencé à étudier le groupe fondamental. Comme eux, nous pourrions
résoudre notre problème en ne considérant  pas simplement des monoïdes $M$ mais des monoïdes
pointés $(M;c)$ où $c\in M$.
\end{rema}
\begin{exem}\label{exem:gr:3.2}
Il est usuel de voir les groupes comme des catégories et les homomorphismes de groupes comme des
foncteurs. Puisque toute catégorie est un groupement et tout foncteur est un g-groupement.
\`A la différence des monoïdes, le choix de l'élément peut-être fait de façon canonique.
Il suffit de prendre l'élément neutre qui est respecté par les homomorphismes de groupes.
\end{exem}
\begin{exem}\label{exem:gr:3.3}
Si $f:\BBB_{1}\to\BBB_{2}$ est un g-foncteur alors le triplet $\left(f^\ast;\BBB^\ast_{1};\BBB^\ast_{2}\right)$, 
où $f^\ast=f$ et $\BBB^\ast_{1}$, $\BBB^\ast_{2}$ sont respectivement les groupements duaux de 
$\BBB_{1}$ et $\BBB_{2}$, est lui aussi un g-foncteur. 
\begin{description}
\item[(GFONC 1)] $f^\ast s^\ast_1=ft_{1}=t_2f=s^\ast_{2}f^\ast $,\quad $f^\ast t^\ast_1=fs_{1}=s_2f=
t^\ast_{2}f^\ast$.
\item[(GFONC 2)] si $x$ et $y$ sont deux éléments de $\BBB^\ast_1$ tels que
        $s^\ast_1(x)=t^\ast_1(y)$, alors $t_{1}(x)=s_{1}(y)$ et 
        $$f^\ast\left(x\OP^\ast_1y\right)=f\left(y\OP x\right)=f(y)\OP_2f(x)=f^\ast(x)\OP^\ast f^\ast(y)$$
\end{description}
$f^\ast:\BBB^\ast_{1}\to\BBB^\ast_{2}$ est le \emph{g-foncteur dual}\index{g-foncteur!g-foncteur dual} 
du g-foncteur $f:\BBB_{1}\to\BBB_{2}$.\\
Il est amusant de noter que l'application d'ensembles à la base de chacun de ces g-foncteurs est la même. 
\end{exem}
\begin{prop}\label{prop:gr:3.1}
Soit $\left(\BBB;s;t;\OP\right)$ un groupement. L'application identité
$\Id_\BBB:\BBB\to\BBB$,
donnée par $\Id(x)=x$ pour chaque $x\in\BBB$, définit un g-foncteur de
$\left(\BBB;s;t;\OP\right)$ dans lui-même.
\end{prop}
\Dem\
\begin{description}
\item[(GFONC 1)] $\Id s(x)=s(x)=s\Id(x)$,\quad $\Id t(x)=t(x)=t\Id(x)$.
\item[(GFONC 2)] Soient $x$ et $y$ sont deux éléments de $\BBB$ tels que
        $s(x)=t(y)$.\\ On a $\Id\left(x\OP y\right)=\Id(x)\OP\Id(y)$.\QED
\end{description}
% lemme
\begin{lemm}\label{lemm:gr:3.1}
Si $f:\BBB_1\to\BBB_2$ un g-foncteur entre deux petits $\MFU$-groupements
$\BBB_1$ et $\BBB_2$,
alors $f$ est un élément de l'univers $\MFU$ \textit{i.e.} un petit ensemble.
\end{lemm}
\Dem
La démonstration est semblable à celle du lemme \ref{lemm:cat:1.1}.
\QEDb
Les g-foncteurs qui vérifient les conditions de ce lemme sont appelés \emph{petits
g-foncteurs}\index{g-foncteur!petit g-foncteur}.
Ainsi nous pouvons parler de l'ensemble $\GFONC$ des petits g-foncteurs.

Considérons les applications
\begin{itemize}
\item $s:\GFONC\to\GFONC$ définie par $s\left(f:\BBB_1\to\BBB_2\right)=\Id_{\BBB_1}$,
\item $t:\GFONC\to\GFONC$ définie par $t\left(f:\BBB_1\to\BBB_2\right)=\Id_{\BBB_2}$,
\item $\OP:\GFONC\times\GFONC\to\GFONC$ définie par
$$f_2\OP f_1=
\begin{cases}
f_2f_1\qquad\text{si }s\left(f_2\right)=t\left(f_1\right)\\
f_1\quad\qquad\text{sinon}
\end{cases}$$
\end{itemize}
Nous devons tout de même vérifier que $f_2\OP f_1$ est bien un g-foncteur quand
$$s(f_2)=t(f_1)=\Id_{\BBB'}$$ Premièrement la composition d'applications d'ensembles $f_2f_1$ est
possible vu la condition imposée à $f_1$ et $f_2$. De plus on a
\begin{description}
\item[(GFONC 1)] \'Evident d'après les définitions de $s$ et $t$ car $f_2f_1$ est une application
ayant même domaine $\BBB_1$ que $f_1$ et même but $\BBB_2$ que $f_2$.
\item[(GFONC 2)] Soient $x$ et $y$ deux éléments de $\BBB_1$ pour lesquels $s_1(x)=t_1(y)$.
Comme $$s'(f_1(x))=f_1s_1(x)=f_1t_1(y)=t'(f_1(y))$$ on a
\begin{align*}
(f_2\OP f_1)(x\OP_1y)=f_2(f_1(x\OP_1 y))&=f_2(f_1(x)\OP' f_1(y))\\
&=(f_2f_1(x))\OP_2(f_2f_1(y))=(f_2\OP f_1)(x)\OP_2(f_2\OP f_1)(y)
\end{align*}
\end{description}
% théorème
\begin{theo}\label{theo:gr:3.1}
$\left(\GFONC;s;t;\OP\right)$ est un groupement. C'est même une catégorie.
\end{theo}
\Dem Nous n'avons qu'à vérifier les axiomes (GR 1) à (GR 3) pour montrer que c'est
un groupement.
\begin{description}
\item[(GR 1)] Pour tout g-foncteur $f$, on a les égalités suivantes :
        $$ss\left(f:\BBB_1\to\BBB_2\right)=s\left(\Id_{\BBB_1}\right)=\Id_{\BBB_1}
        =s\left(f:\BBB_1\to\BBB_2\right)$$
        $$st\left(f:\BBB_1\to\BBB_2\right)=s\left(\Id_{\BBB_2}\right)=\Id_{\BBB_2}
        =t\left(f:\BBB_1\to\BBB_2\right)$$
        $$tt\left(f:\BBB_1\to\BBB_2\right)=t\left(\Id_{\BBB_2}\right)=\Id_{\BBB_2}
        =t\left(f:\BBB_1\to\BBB_2\right)$$
        $$ts\left(f:\BBB_1\to\BBB_2\right)=t\left(\Id_{\BBB_1}\right)=\Id_{\BBB_1}
        =s\left(f:\BBB_1\to\BBB_2\right)$$
\item[(GR 2)] Supposons que $f_1:\BBB_1\to\BBB'_1$ et $f_2:\BBB_2\to\BBB'_2$ 
        sont deux 
        g-foncteurs tels que $$s\left(f_2\right)=t\left(f_1\right)$$ On a
        $$s\left(f_2\OP f_1\right)=
        s\left(f_2f_1\right)=\Id_{\BBB_1}=
        s\left(f_1\right)$$
        $$t\left(f_2\OP f_1\right)=
        s\left(f_2f_1\right)=
        \Id_{\BBB'_2}=
        t\left(f_2\right)$$
\item[(GR 3)] C'est une conséquence directe de l'associativité de la 
        composition des applications d'ensembles.
\end{description}
C'est bien une catégorie car l'axiome (CAT 3) est vérifié par les calculs ci-dessous.
$$f\OP s(f)=f\OP\Id_{\BBB_1}=f\Id_{\BBB_1}=f\quad\text{et}\quad
t(f)\OP f=\Id_{\BBB_2}\OP f=\Id_{\BBB_2}f=f$$
où $f:\BBB_1\to\BBB_2$ est un foncteur.
\QED 
\begin{rema}\label{rema:gr:3.3}
Dans ce théorème, nous avons défini $f_2\OP f_1$ en prenant $f_1$ quand $s(f_2)
\neq t(f_1)$. En regardant la démonstration, on s'aperçoit immédiatement que ce choix
est arbitraire. On aurait très bien pu choisir $f_2$ ou autre chose. Ces ainsi que nous 
sommes amenés à dire que 
deux groupements $\left(\BBB_1;s_1;t_1;\OP_1\right)$ et 
$\left(\BBB_2;s_2;t_2;\OP_2\right)$ sont 
\emph{presque égaux}\index{groupement!groupements presque égaux} si 
$\BBB_1=\BBB_2$, $s_1=s_2$, $t_1=t_2$ et $x\OP_1y=x\OP_2y$ lorsque 
$s_1(x)=s_2(x)=t_2(y)=t_1(y)$. Il est évident que l'application 
identité $\Id_{\BBB_1}=
\Id_{BBB_2}$ définit un g-foncteur 
$\Id_{\BBB_1}^{\BBB_2}:\BBB_1\to\BBB_2$ 
et un g-foncteur $\Id_{\BBB_2}^{\BBB_1}:\BBB_2\to\BBB_1$ tels que 
$$\Id_{\BBB_1}^{\BBB_2}\Id_{\BBB_2}^{\BBB_1}=\Id_{\BBB_2}\qquad\text{et}\qquad 
\Id_{\BBB_2}^{\BBB_1}\Id_{\BBB_1}^{\BBB_2}=\Id_{\BBB_1}$$
En d'autres termes, deux groupements presque égaux sont isomorphes.
Il est clair que cela définit une relation d'équivalence sur l'ensemble des 
petits groupements.
\end{rema}
L'idée centrale dans notre travail est que les objets et donc les identités ne sont pas des notions
très \Dire{naturelles} et doivent être, autant que possible, éliminés de notre théorie. Ainsi
l'axiome (GFONC 1) cadre mal avec cette idée. Il a été introduit, en théorie des catégories, afin de 
rendre les identités \Dire{stables} par les foncteurs. L'axiome (GFONC 2) intéressant mais imparfait.
Il faut que les compositions aient un sens. Ce dernier étant, précédemment, donné par l'axiome (GFONC 1),
nous sommes amenés à imposer quelques conditions sur les éléments.
D'où la notion beaucoup plus naturelle de g-morphisme.

Un \emph{g-morphisme}\index{g-morphisme} $f$ est un triplet
$\left((\BBB_1;s_1;t_1;
\OP_1);(\BBB_2;s_2;t_2;\OP_2);f\right)$, souvent noté
$$f:\BBB_1\MorTo\BBB_2$$ où les deux premiers termes sont des groupements et $f$ est
une application d'ensembles $$f:\BBB_1\to\BBB_2$$ vérifiant l'axiome suivant
\begin{description}
\item[(GMOR)] si $x$ et $y$ sont deux éléments de $\BBB_1$ tels que
        $s_1(x)=t_1(y)$, alors $$s_2(f(x))=t_2(f(y))\quad\text{et}\quad f\left(x\OP_1y\right)=f(x)\OP_2f(y)$$
\end{description}
Il est évident que tout g-foncteur est un g-morphisme puisque $$s_2(f(x))=f(s_1(x))=f(t_1(y))=t_2(f(y))$$ Comme pour
les g-foncteurs, on peut parler de petits g-morphismes et montrer que ceux-ci forment un petit ensemble.
Notons $\GMOR$ l'ensemble des petits g-morphismes.
\begin{theo}\label{theo:gr:3.2}
Le quadruplet $\left(\GMOR;s;t;\OP\right)$ où $$s:\GMOR\to\GMOR,\quad t:\GMOR\to\GMOR\quad\text{et}
\quad\OP:\GMOR\times\GMOR\to\GMOR$$ sont les applications définies par
\begin{itemize}
\item $s\left(f:\BBB_1\to\BBB_2\right)=\Id_{\BBB_1}$,
\item $t\left(f:\BBB_1\to\BBB_2\right)=\Id_{\BBB_2}$,
\item $f_2\OP f_1=
\begin{cases}
f_2f_1\qquad\text{si }s\left(f_2\right)=t\left(f_1\right)\\
f_1\quad\qquad\text{sinon}
\end{cases}$
\end{itemize}
est un groupement et même une catégorie.
\end{theo}
La démonstration est identique à celle du théorème \ref{theo:gr:3.1}.
%%%%%%%%%%%%%%
%%%%  section 4  %%%%
%%%%%%%%%%%%%%
\section{Les g-transformations : Une première approche.}\label{sec:gr:4}
Alors que nous avons introduit la notion de g-mor\-phi\-sme pour satisfaire des exigences 
es\-thé\-ti\-ques, il n'en va pas du tout de la même manière pour généraliser les transformations 
naturelles. Le problème se situe au niveau même de leur définition.
Elles le sont traditionnellement comme une famille, indexée par les objets, de morphismes qui 
vérifient des conditions de commutativité de diagrammes. Or nous n'avons pas d'objets et l'étude 
qui suit va nous permettre d'observer que les axiomes (NAT 1) et (NAT 3) ne nous permettent pas de 
construire une théorie satisfaisante si nous ne rajoutons pas une condition sur les groupements 
qui en font pratiquement des catégories. Commençons par quelques remarques.

Premièrement, l'axiome (TRANS 1) et sa forme équivalente (TRANS 1')  ont pour unique but 
de permettre la création des diagrammes suivants
$$\xymatrix{
F_1(X)\ar[r]^{\eta_X}\ar[d]_{F_1(f)}\ar[dr] & F_2(X)\ar[d]^{F_2(f)}\\
F_1(Y)\ar[r]^{\eta_Y} & F_2(Y)
}$$
pour toute transformation naturelle $\eta:F_1\leadsto F_2$ et tout morphisme $f:X\to Y$.\\
Deuxièmement, l'axiome (TRANS 3) impose à ces différents carrés d'être commutatifs. Ce qui signifie 
que les deux compositions latérales possibles donnent le même résultat.\\
Troisièmement, si nous prenons deux morphismes composables $f$ et $g$, nous obtenons le diagramme commutatif
$$\xymatrix{
F_{1}(X)\ar[r]^{F_{1}(f)}\ar[d]^{\eta_{X}} & F_{1}(Y)\ar[r]^{F_{1}(g)}\ar[d]^{\eta_{Y}} & 
F_{1}(Z)\ar[d]^{\eta_{Z}}\\
F_{2}(X)\ar[r]_{F_{2}(f)} & F_{2}(Y)\ar[r]_{F_{2}(g)} & F_{2}(Z)}$$
Par conséquent, en notant $\eta^{1}(h)=\eta_{R}$ et $\eta^{2}(h)=\eta_{S}$ pour tout morphisme $h:R\to S$
nous avons les égalités
$$F_{2}(f)\OP_{2}\eta^{1}(f)=\eta^{2}(f)\OP_{2}F_{1}(f)$$
et
$$\eta^{1}(g)=\eta^{2}(f)$$
Finalement il semble raisonnable de définir les g-transformations de la manière suivante :\\
Une \emph{($\MFU$-)g-transformation}\index{g-transformation}
\index{g-transformation!$\MFU $-g-transformation} $\left(\eta^{1},\eta^{2};f_1,f_2\right)$ entre le
g-morphisme $f_1:\BBB_1\to\BBB_2$ et le g-morphisme $f_2:\BBB_1\to\BBB_2$ 
est un quadruplet où $\eta^{1},\eta^{2}:\BBB_1\to\BBB_2$ sont deux applications
d'ensembles satisfaisant les deux axiomes ci-dessous
\begin{description}
\item[(GTRANS 1)] $s_2\eta^{1}=s_{2}f_1$,\quad $t_2\eta^{1}=s_{2}f_2$,\quad $s_{2}\eta^{2}=t_{2}f_{1}$,
        \quad $t_{2}\eta_{2}=t_{2}f_{2}$;
\item[(GTRANS 2)] pour tous $x$ et $y$ dans $\BBB_1$ vérifiant $s_{1}(x)=t_{1}(y)$,
        on a
        $$f_{2}(x)\OP_{2}\eta^{1}(x)=\eta^{2}(x)\OP_2 f_1(x)$$
        et
        $$\eta^{1}(x)=\eta^{2}(y)$$
\end{description}
La notation $\left(\eta^1;\eta^2;f_1;f_2\right)$ n'étant pas très explicite, nous 
écrirons plus souvent $$(\eta^1;\eta^2):f_1\leadsto f_2$$
et même $$\eta:f_{1}\leadsto f_{2}$$
Dans l'axiome (GTRANS 2), la première égalité porte 
sur $x$. La proposition suivante nous montre que l'on peut indifféremment remlpacer le $x$ par $y$.
% proposition
\begin{prop}\label{prop:gr:4.1}
L'axiome (GTRANS 2) est équivalent à l'axiome
\begin{description}
\item[(GTRANS 2')] pour tous $x$ et $y$ dans $\BBB_1$ vérifiant $s_{1}(x)=t_{1}(y)$,
        on a
        $$f_{2}(y)\OP_{2}\eta^{1}(y)=\eta^{2}(y)\OP_2 f_1(y)$$
        et
        $$\eta^{1}(x)=\eta^{2}(y)$$
\end{description}
\end{prop}
\Dem\ Supposons que l'axiome (GTRANS 2) soit vérifié. Il nous suffit de montrer qu'alors (GTRANS 2') l'est aussi.
D'après l'axiome (GR 1), les éléments $y$ et $s_{1}(y)$ vérifient les conditions de l'axiome (GTRANS 2). D'où
$$f_{2}(y)\OP_{2}\eta^{1}(y)=\eta^{2}(y)\OP_2 f_1(y)$$
        et
$$\eta^{1}(y)=\eta^{2}(s_{1}(y))$$ En prenant $t_{1}(x)$ et $x$, on démontre la réciproque avec en plus
$$\eta^1(t_{1}(x))=\eta^2(x)$$\QEDb
Cette démonstration nous permet de voir que les applications $\eta^1$ et $\eta^2$ sont reliées par des relations 
très strictes.
% proposition
\begin{prop}\label{prop:gr:4.2}
Pour toute g-transformation $(\eta^1;\eta^2):f_{1}\leadsto f_{2}$ entre les g-morphismes $f_{1}$ et $f_{2}$ 
de $(\BBB_{1};s_{1};t_{1})$ dans $(\BBB_{2};s_{2};t_{2})$, on a les égalités suivantes :
$$\eta^1=\eta^2s_{1},\quad\eta^{2}=\eta^1t_{1},\quad\eta^1=\eta^1s_{1},\quad\eta^2=\eta^2t_{1}$$
\end{prop}
\Dem Les deux premières égalités ont été vues dans la démonstration de la proposition précédente. Les deux suivantes
découlent immédiatement d'elles et de l'axiome (GR 1).\QEDb
Cette définition semble à priori éloignée de celle donnée pour les transformations 
naturelles entre foncteurs de catégories. Les deux propositions suivantes vont nous montrer qu'en fait il 
n'en est rien.
% proposition
\begin{prop}\label{prop:gr:4.3}
Soient $f_{1}$ et $f_{2}$ deux g-morphismes entre les groupements $(\BBB_{1};s_{1};t_{1})$ et 
$(\BBB_{2};s_{2};t_{2})$.\\
Si $\left(\eta^{1};\eta^{2}\right):f_{1}\leadsto f_{2}$ est une g-transformation, alors 
l'application $\eta^{1}:\BBB_{1}\leadsto \BBB_{2}$ vérifie les conditions suivantes :
\begin{enumerate}
\item $s_{2}\eta^1 s_{1}=s_{2}f_{1}$\quad et\quad$t_{2}\eta^1 s_{1}=s_{2}f_{2}$
\item $\eta^1=\eta^1s_{1}$
\item Pour tout $x\in\BBB_{1}$, on a $f_{2}(x)\OP_{2}\eta^1(s_{1}(x))=
         \eta^1(t_{1}(x))\OP_{2}f_{1}(x)$.
\end{enumerate}
\end{prop}
\Dem Supposons que $(\eta^1;\eta^2):f_{1}\leadsto f_{2}$ soit une g-transformation.
La seconde propriété a été vue dans la proposition \ref{prop:gr:4.2}.\\
La première est la conséquence directe de l'axiome (GTRANS 1) et de la seconde propriété.
Pour la dernière, considérons un élément $x$ de $\BBB_{1}$. D'après l'axiome (GR 1), $x$ et 
$s_{1}(x)$ vérifient la condition de l'axiome (GTRANS 2). Doù
$$f_2(x)\OP_{2}\eta^1(x)=\eta^2(x)\OP_2 f_1(x)$$
La proposition \ref{prop:gr:4.2} nous dit que
$$\eta^1(x)=\eta^1(s_{1}(x))\quad\text{et}\quad\eta^2(x)=\eta^1(t_{1}(x))$$
D'où le résultat.\QEDb
Pour démontrer la réciproque dans le cas où $f_{1}$ et $f_{2}$ sont des g-foncteurs, nous aurons besoin
du
% lemme
\begin{lemm}\label{lemm:gr:4.1}
Si $f_{1}$ et $f_{2}$ sont des g-foncteurs de $\BBB_{1}$ vers $\BBB_{2}$, alors toute application 
$\eta:\BBB_{1}\to\BBB_{2}$ qui satisfait les conditions $$s_{2}\eta s_{1}=s_{2}f_{1}\quad\text{et}
\quad t_{2}\eta s_{1}=s_{2}f_{2}$$ vérifient aussi les conditions
 $$s_{2}\eta t_{1}=t_{2}f_{1}\quad\text{et}
\quad t_{2}\eta t_{1}=t_{2}f_{2}$$
\end{lemm}
\Dem\ C'est une conséquence immédiate des axiomes (GR 1) et (GFONC 1) :
$$s_{2}\eta t_{1}=s_{2}\eta s_{1}t_{1}=s_{2}f_{1}t_{1}=s_{2}t_{2}f_{1}=t_{2}f_{2}$$et
$$t_{2}\eta t_{1}=t_{2}\eta s_{1}t_{1}=s_{2}f_{2}t_{1}=s_{2}t_{2}f_{2}=t_{2}f_{2}$$\QED
% proposition
\begin{prop}\label{prop:gr:4.4}
Si $f_{1}$ et $f_{2}$ sont des g-foncteurs, alors toute fonction $\eta$ de $\BBB_{1}$ dans $\BBB_{2}$ 
qui vérifient les trois conditions  de la proposition \ref{prop:gr:4.3} définit une g-trans\-for\-ma\-tion 
$(\eta;\eta t_{1}):f_{1}\leadsto f_{2}$.
\end{prop}
\Dem\ Considérons une application $\eta:\BBB_{1}\to\BBB_{2}$ qui vérifient les 
trois propriétés précédentes. Posons $\eta^1=\eta$ et $\eta^2=\eta t_{1}$. Montrons que 
$(\eta^1;\eta^2)$ est une g-transformation de $f_{1}$ vers $f_{2}$.
\begin{description}
\item[(GTRANS 1)] Puisque $\eta^1=\eta^1s_{1}$, on a, d'après la première condition,
         $$s_{2}\eta^{1}=s_{2}f_{1}\quad\text{et}\quad t_{2}\eta^{1}=s_{2}f_{2}$$
         Comme $\eta^{2}=\eta t_{1}$, le lemme \ref{lemm:gr:4.1}, nous donne en plus
         $$s_{2}\eta^{2}=t_{2}f_{1}\quad\text{et}\quad t_{2}\eta^{2}=t_{2}f_{2}$$
\item[(GTRANS 2)] Soient $x$ et $y$ telles que $s_{1}(x)=t_{1}(y)$. Tout découle de $\eta=\eta s_{1}$ 
         et $\eta^2=\eta^1t_{1}$. En effet à la vue de la condition imposée à $x$ et $y$, on a         
         $$\eta^1(x)=\eta^1(s_{1}(x))=\eta^1(t_{1}(y))=\eta^2(y)$$Et d'après la troisième 
         propriété vérifiée par $\eta$, on a aussi
         $$f_{2}(x)\OP_{2}\eta^1(x)=\eta^2(x)\OP_{2}f_{1}(x)$$\QED
\end{description}
Avant de poursuivre, continuons à étudier de plus près ce qui se passe 
pour les transformations naturelles de la théorie des catégories.\\
Une propriété classique des transformations naturelles de catégories est que chaque foncteur 
$F:\BBC_{1}\to\BBC_{2}$ définit 
une transformation naturelle de $F$ dans lui-même. En effet, pour tout objet $X$ de $\BBC_{1}$, posons 
$\eta_{X}=F\left(\Id_{X}\right)$. Alors, pour tout morphisme $f:X\to Y$ de $\BBC_{1}$, le diagramme 
$$\xymatrix{
F(X)\ar[r]^{F(f)}\ar[d]_{\eta_{X}} & F(Y)\ar[d]^{\eta_{Y}}\\
F(X)\ar[r]_{F(f)} & F(Y)}$$
est commutatif car 
$$\eta_{Y}\OP_{2}F(f)=F(\Id_{Y})\OP_{2}F(f)=F(\Id_{Y}\OP_{1}f)
=F(f\OP_{1}\Id_{X})=F(f)\OP_{2}F(\Id_{X})=F(f)\OP_{2}\eta_{X}$$
Plus précisément, nous remarquons que nous y avons 
utilisé deux choses. La première est que, $F$ étant un foncteur, l'image d'une composition est égale à la 
composition des images ; et la seconde, que $\Id_{Y}\OP_{1}f=f\OP_{1}\Id_{X}$ dans la catégorie $\BBC_{1}$.
En d'autres termes, nous avons utilisé certaines propriétés des identités. Or le principe qui nous a guidé jusqu'à
maintenant est d'enlever les références à ces dernières. Cette impossibilité pour les g-foncteurs 
de devenir des g-transformations nous empêche de définir des applications source et but sur l'ensemble 
des g-transformations. C'est en cela qu'il ne semble pas y avoir de théorie satisfaisante  de g-transformation 
sans l'ajout de quelque chose de supplémentaire. Nous reviendrons là-dessus dans le chapitre \ref{chap:2-gr}.
Malgré tout, et afin de bien illustrer notre propos nous allons continuer la construction en supposant que les 
groupements que nous allons considérer dans le reste de cette section satisfont tous la 
condition supplémentaire
$$\mathbf{(\star)}\qquad t(x)\OP x=x\OP s(x),\quad\text{pour tout }x\in\BBB$$
Pour bien garder celà à l'esprit, nous parlerons de $\star$-groupements.

La proposition qui suit va nous permettre de définir sur l'ensemble des petites 
g-trans\-for\-ma\-tions entre g-morphismes de $\star$-groupements des applications source et but.
% proposition
\begin{prop}\label{prop:gr:4.5}
Si $f$ est un g-morphisme du $\star$-groupement
$\left(\BBB_1;s_1;t_1;\OP_1\right)$ vers le $\star$-groupement
$\left(\BBB_2;s_2;t_2;\OP_2\right)$, alors $(s_{2}f;t_{1}f)$ est une g-transformation de $f$ dans $f$.
\end{prop}
\Dem C'est une simple vérification des axiomes (GTRANS 1) et 
(GTRANS 2). Posons $\eta^1=s_{2}f$ et $\eta^{2}=t_{2}f$
\begin{description}
\item[(GTRANS 1)] \'Evident.
\item[(GTRANS 2)] Soient $x$ et $y$ deux éléments de $\BBB_{1}$ satisfaisant $s_{1}(x)=t_{1}(y)$.
 En utilisant la condition $\mathbf{\star}$, on trouve  
 \begin{align*}
 f(x)\OP_{2}\eta^1(x)&=f(x)\OP_{2}s_{2}\left(f(x)\right)\\
 &=t_{2}\left(f(x)\right)\OP_{1}f(x)\\
 &=\eta^2(x)\OP_{2}f(x)
 \end{align*}
De plus, grace à l'axiome (GMOR), on trouve$$\eta^1(x)=s_{2}f(x)=t_{2}f(y)=\eta^2(y)$$\QED
\end{description}
% remarque
\begin{rema}\label{rema:gr:4.1}
L'application identité $\Id_{\BBB}$ de l'ensemble $\BBB$ dans lui-même définit clairement un
g-morphisme, et même un g-foncteur, de $(\BBB;s;t;\OP)$ dans lui-même. Ainsi $(s_{1};t_{1})$ 
est une g-transformation de $\Id_{\BBB}$ vers lui-même.
\end{rema}
On dit qu'une g-transformation $\eta:f_1\leadsto f_2$ est \emph{petite}
\index{g-transformation!petite g-transformation} quand 
les g-morphismes $f_1$ et $f_2$ le sont.
% lemme
\begin{lemm}\label{lemm:gr:4.2}
Toute petite $\MFU$-g-transformation est un élément de l'univers $\MFU$.
\end{lemm}
\Dem Semblable à la démonstration du lemme \ref{lemm:gr:3.1}. \QEDb
Désignons par $\star\GTrans$ l'ensemble de toutes les petites g-transformations au dessus 
de $\star$-groupements. 
Comme le théorème suivant va le montrer, il est aisé de munir cet ensemble 
d'une structure de groupement.
% théorème
\begin{theo}\label{theo:gr:4.1} 
On peut définir deux applications $\sigma_1$ et $\tau_1$ de $\star\GTrans$ dans lui-même 
de la manière suivante :
$$\sigma_1\left(\eta:f_1\leadsto f_2\right)=\left((s_{2}f_1;t_{2}f_{1}):
f_1\leadsto f_1\right)$$et
$$\tau_1\left(\eta:f_1\leadsto f_2\right)=\left((s_{2}f_2;t_{2}f_{2}):f_2
\leadsto f_2\right)$$
Soient deux g-transformations
$\eta_1:f_1\leadsto f'_1$, $\eta_2:f_2\leadsto f'_2$.
\begin{itemize}
\item Si $\sigma_1(\eta_2)=\tau_1(\eta_1)$, notons 
        $\eta_2\circledast\eta_1:f_1\leadsto f'_2$
        la g-transformation définie, pour tout $x\in\BBB_{1}$, par $$\left(\eta_2
        \circledast\eta_1\right)(x)=\left(\eta^1_{2}(x)\OP_{2}\eta^1_{1}(x);
        \eta^2_{2}(x)\OP_{2}\eta^2_{2}(x)\right)$$
\item sinon, posons $\eta_2\circledast\eta_1=\eta_1:f_1\leadsto f'_1$.
\end{itemize} 
On obtient une application $\circledast$ de $\star\GTrans\times\star\GTrans$ 
dans $\star\GTrans$.\\ Le quadruplet $\left(\star\GTrans;\sigma_1;\tau_1;
\circledast\right)$ est un groupement.
\end{theo}
\Dem Avant de vérifier les axiomes des groupements, nous devons montrer 
que si $\sigma_1(\eta_2)=\tau_1(\eta_1)$ alors 
$\eta_2\circledast\eta_1:f_1\leadsto f'_2$ est bien 
une g-transformation. Pour plus de clareté, posons $\eta_2\circledast\eta_1=(\mu^1;\mu^2)$
\begin{description}
\item[(GTRANS 1)] Remarquons que $\sigma_1(\eta_2)=\tau_1(\eta_1)$ implique 
         $f_{2}=f'_{1}$. 
         En effet $\sigma_1(\eta_2)=\tau_1(\eta_1)$ signifie 
         $$(s_{2}f_{2};t_{2}f_{2};f_{2};f_{2})=(s_{2}f'_{1};t_{2}f'_{1};f'_{1};f'_{1})$$
        Ainsi, l'axiome (GR 2) implique
        $$s_{2}\mu^1=s_2\left(\eta^1_2\OP_{2}\eta^1_1\right)
        =s_2\eta^1_1=s_{2}f_1$$
        et
         $$t_{2}\mu^1=t_2\left(\eta^1_2\OP_{2}\eta^1_1\right)
        =t_2\eta^1_2=s_{2}f'_2$$
        car, d'après (GTRANS 1) et la remarque ci-dessus, 
        $s_{2}\eta^1_{2}=s_{2}f_{2}=s_{2}f'_{1}=t_{2}\eta^1_{1}$.\\
        De même, puisque  $s_{2}\eta^2_{2}=t_{2}f_{2}=t_{2}f'_{1}=t_{2}\eta^2_{1}$, on trouve
        $$s_{2}\mu^2=s_2\left(\eta^2_2\OP_{2}\eta^2_1\right)
        =s_2\eta^2_1=t_{2}f_1$$
        et
         $$t_{2}\mu^2=t_2\left(\eta^1_2\OP_{2}\eta^1_1\right)
        =t_2\eta^2_2=t_{2}f'_2$$
\item[(GTRANS 2)] Pour tous $x\in\BBB_1$ et $y\in\BBB_{1}$ tels que $s_{2}(x)=t_{2}(y)$, les 
        axiomes (GTRANS 2) vérifiés par $\eta_{1}$ et $\eta_{2}$ et le fait que $f_{2}=f'_{1}$, nous 
        permettent d'écrire
        \begin{align*}
                f'_2(x)\OP_2\mu^1(x)
                &=f'_2\OP_2\left(\eta^1_2(x)\OP_2\eta^1_1(x)\right)\\
                &=\left(f'_2(x)\OP_2\eta^1_2(x)\right)\OP_2\eta^1_1(x)\\
                &=\left(\eta^2_2(x)\OP_2f_2(x)\right)\OP_2\eta^1_1(x)\\
                &=\eta^2_2(x)\OP_2\left(f_2(x)\OP_2\eta^1_1(x)\right)\\
                &=\eta^2_2(x)\OP_2\left(f'_1(x)\OP_2\eta^1_1(x)\right)\\
                &=\eta^2_2(x)\OP_2\left(\eta^2_1(x)\OP_2f_1(x)\right)\\
                &=\left(\eta^2_2(x)\OP_2\eta^2_1(x)\right)\OP_2f_1(x)
                =\mu^2(x)\OP_2f_1(x)
        \end{align*}et
        $$\mu^1(x)=\eta^1_{2}(x)\OP_{2}\eta^1_{1}(x)
        =\eta^2_{2}(y)\OP_{2}\eta^2_{1}(y)=\mu^2(y)$$
\end{description} 
Passons maintenant à la seconde partie de notre démonstration.
\begin{description}
\item[(GR 1)] Pour tout $\eta:f_1\leadsto f_2$ appartenant à $\star\GTrans$, on a
        $$\left(\sigma_1\sigma_1\right)\left(\eta\right)
        =\sigma_1\left((s_{2}f_1;t_{2}f_{1}):f_1\leadsto f_1\right)
        =(s_{2}f_1;t_{2}f_{1}):f_1\leadsto f_1=\sigma_1\left(\eta\right)$$
        $$\left(\sigma_1\tau_1\right)\left(\eta\right)
        =\sigma_1\left((s_{2}f_2;t_{2}f_{2}):f_2\leadsto f_2\right)
        =(s_{2}f_2;t_{2}f_{2}):f_2\leadsto f_2=\tau_1\left(\eta\right)$$
        $$\left(\tau_1\tau_1\right)\left(\eta\right)
        =\tau_1\left((s_{2}f_2;t_{2}f_{2}):f_2\leadsto f_2\right)
        =(s_{2}f_{2};t_{2}f_2):f_2\leadsto f_2=\tau_1\left(\eta\right)$$
        $$\left(\tau_1\sigma_1\right)\left(\eta\right)
        =\tau_1\left((s_{2}f_{1};t_{1}f_1):f_1\leadsto f_1\right)
        =(s_{2}f_{1};t_{2}f_1):f_1\leadsto f_1=\sigma_1\left(\eta\right)$$
\item[(GR 2)] Soient $\eta_1:f_1\leadsto f'_1$ et $\eta_2:f_2\leadsto f'_2$ deux 
        éléments de $\star\GTrans$ tels que $\sigma_1\left(\eta_2\right)
=\tau_1\left(\eta_1\right)$. 
        Dans ce cas $\eta_2\circledast\eta_1$ est une g-transformation de $f_1$ dans 
        $f'_2$. Par conséquent
        $$\sigma_1\left(\eta_2\circledast\eta_1\right)
        =(s_{2}f_1;t_{2}f_{1})
        =\sigma_1\left(\eta_1\right)
        \quad\text{et}\quad
        \tau_1\left(\eta_2\circledast\eta_1\right)
        =(s_{2}f'_{2};t_{2}f'_2)
        =\tau_1\left(\eta_2\right)$$
\item[(GR 3)] Soient $\eta_1:f_1\leadsto f'_1$, $\eta_2:f_2\leadsto f'_2$ et 
        $\eta_3:f_3\leadsto f'_3$ trois éléments de $\star\GTrans$ tels que     
        $\sigma_1\left(\eta_2\right)=\tau_1\left(\eta_1\right)$ et 
        $\sigma_1\left(\eta_3\right)=\tau_1\left(\eta_2\right)$.\\
        Comme $\sigma_1\left(\eta_3\circledast\eta_2\right)=\sigma_1\left(\eta_2\right)
        =\tau_1\left(\eta_1\right)$ et $\tau_1\left(\eta_2\circledast\eta_1\right)=
        \tau_1\left(\eta_2\right)=\sigma_1\left(\eta_3\right)$, on a
        $$s_{2}(\eta^1_{3}\OP_{2}\eta^1_{2})=t_{2}\eta^1_{1}$$et
        $$s_{2}(\eta^2_{3}\OP_{2}\eta^2_{2})=t_{2}\eta^2_{1}$$
        En utilisant l'axiome (GR 3), on en déduit
        $$\left(\eta^1_3\OP_{2}\eta^1_2\right)\OP_{2}\eta^1_1
        =\eta^1_3\OP_2\left(\eta^1_2\OP_2\eta^1_1\right)$$
        et
        $$\left(\eta^2_3\OP_{2}\eta^2_2\right)\OP_{2}\eta^2_1
        =\eta^2_3\OP_2\left(\eta^2_2\OP_2\eta^2_1\right)$$
        D'où
        $$\left(\eta_{3}\circledast\eta_{2}\right)\circledast\eta_{1}
        =\eta_{3}\circledast\left(\eta_{2}\circledast\eta_{1}\right)$$
        \QED
\end{description}
% remarque
\begin{rema}\label{rema:gr:4.2}
\`A la différence de ce qui se passe pour les g-foncteurs, les g-transformations ne forment 
pas une catégorie. Ce qui n'est pas surprenant car c'est justement pour cette raison 
que nous avons introduit la notion de groupement. 
\end{rema}
On peut munir l'ensemble des (petites) g-transformations de deux autres structures de 
groupement. Pour cela nous aurons besoin du résultat intermédiaire suivant
% proposition
\begin{prop}\label{prop:gr:4.6}
Soient $f:\BBB'_1\to\BBB_1$, $g:\BBB_2\to\BBB'_2$ deux g-morphismes et 
$(\eta^1;\eta^2):(f_1:\BBB_1\to\BBB_2)\leadsto(f_2:\BBB_1\to\BBB_2)$ une g-transformation.
Les quadruplets $\left(\eta^1f;\eta^2f;f_1f;f_2f\right)$ et 
$\left(g\eta^1;g\eta^2;gf_1;gf_2\right)$ sont des g-transformations.
\end{prop}
\Dem Nous nous contenterons de faire la démonstration pour $(\eta^1f;\eta^2f)$ car celle de 
$(g\eta^1;g\eta^2)$ est similaire.
\begin{description}
\item[(GTRANS 1)] Comme $(\eta^1;\eta^2)$ est une g-transformation
         $$s_2(\eta^1f)=(s_{2}\eta^1)f=(s_{2}f_{1})f=s_{2}(f_{1}f)$$ 
         $$t_2(\eta^1f)=(t_2\eta^1)f=(s_{2}f_2)f=s_{2}(f_2f)$$
         et
         $$s_2(\eta^2f)=(s_{2}\eta^2)f=(t_{2}f_{1})f=t_{2}(f_{1}f)$$ 
         $$t_2(\eta^2f)=(t_2\eta^2)f=(t_{2}f_2)f=t_{2}(f_2f)$$ 
\item[(GTRANS 2)] Pour tout $x\in\BBB'_1$ et tout $y\in\BBB'_{1}$ 
         tels que $s'_{1}(x)=t'_{1}(x)$, 
         \begin{align*}
         (f_2f)(x)\OP_2(\eta^1 f)(x)
         &=f_2(f(x))\OP_2\eta^1(f(x))\\
         &=\eta^2(f(x))\OP_2 f_1(f(x))
         =(\eta^2 f)(x)\OP_2(f_1f)(x)
         \end{align*}
         car $\eta:f_1\leadsto f_2$ est une g-transformation et $f(x)\in\BBB_1$,
         et
         $$(\eta^1f)(x)=\eta^1(f(x))=\eta^2(f(x))=(\eta^2f)(x)$$
         d'après l'axiome (GMOR), $s_{2}f(x)=t_{2}f(y)$.\QED
\end{description}
Pour simplifier, si $\eta=(\eta^1;\eta^2)$ est la g-transformation de la proposition, 
alors on note $\eta f$ et $g\eta$ les nouvelles g-transformations.
% lemme
\begin{lemm}\label{lemm:gr:4.3}
Soient $\eta:f_{1}\leadsto f_{2}$ et $\eta':f'_{1}\leadsto f'_{2}$ deux 
g-transformations pour lesquelles le groupement but de $f_{1}$ et $f_{2}$ est aussi le 
groupement source de $f'_{1}$ et $f'_{2}$. Alors les g-transformations
$$(f'_{1}\eta)\circledast(\eta'f_{2}):f'_{1}f_{1}\leadsto f'_{2}f_{2}$$ et 
$$(\eta'f_{1})\circledast(f'_{2}\eta):f'_{1}f_{1}\leadsto f'_{2}f_{2}$$
existent.
\end{lemm}
\Dem\ Puisque $f'_{1}\eta:f'_{1}f_{1}\leadsto f'_{1}f_{2}$ et 
$\eta'f_{2}:f'_{1}f_{2}\leadsto f'_{2}f_{2}$, la première g-transformation 
$(f'_{1}\eta)\circledast(\eta'f_{2})$ existe. L'existence de la seconde 
se montre de la même manière.\QEDb
Grâce à ce lemme, on peut définir de nouvelles applications :
% théorème
\begin{theo}\label{theo:gr:4.2}
Considérons les quatre applications d'ensembles $\sigma_0$, $\tau_0$,
$\boxtimes$ et $\boxdot$sont les applications ci-dessous.
\begin{itemize}
\item $\sigma_0:\star\GTrans\to\star\GTrans$ définie par $$\sigma_0\left(\eta:
(f_1:\BBB_1\to\BBB_2)\leadsto(f_2:\BBB_1\to\BBB_2)\right)
=(s_{1};t_{1}):\Id_{\BBB_1}\leadsto\Id_{\BBB_1}$$
\item $\tau_0:\star\GTrans\to\star\GTrans$ définie par $$\tau_0\left(\eta:
(f_1:\BBB_1\to\BBB_2)\leadsto(f_2:\BBB_1\to\BBB_2)\right)
=(s_{2};t_{2}):\Id_{\BBB_2}\leadsto\Id_{\BBB_2}$$
\item $\boxtimes,\boxdot:\star\GTrans\times\star\GTrans\to\star\GTrans$ 
définies par
$$\eta_2\boxtimes\eta_1=
\begin{cases}
(f'_2\eta_1)\circledast(\eta_2f_1)
&\text{si}\ \sigma_0(\eta_2)=\tau_0(\eta_1)\\
\eta_1&\text{sinon}
\end{cases}$$et
$$\eta_2\boxdot\eta_1=
\begin{cases}
(\eta_2f'_1)\circledast(f_2\eta_1)
&\text{si}\ \sigma_0(\eta_2)=\tau_0(\eta_1)\\
\eta_1&\text{sinon}
\end{cases}$$
pour $\eta_1:(f_1:\BBB_1\to\BBB'_1)\leadsto(f'_1:\BBB_1\to\BBB'_1)$ et
$\eta_2:(f_2:\BBB_2\to\BBB'_2)\leadsto(f'_2:\BBB_2\to\BBB'_2)$ deux 
g-trans\-for\-ma\-tions.
\end{itemize}
Les deux quadruplets $\left(\GTrans;\sigma_0;\tau_0;\boxtimes\right)$ 
et $\left(\GTrans;\sigma_0;\tau_0;\boxdot\right)$ sont des grou\-pe\-ments.
\end{theo}
\Dem
La condition $\sigma_0(\eta_2)=\tau_0(\eta_1)$ n'est autre que celle imposée 
dans le lemme \ref{lemm:gr:4.3}. Par conséquent $\eta_{2}\boxtimes\eta_{1}$ 
et $\eta_{2}\boxdot\eta_{1}$ existent et sont bien des g-groupements.
\begin{description}
\item[(GR 1)] Pour tout $\eta:(f_1:\BBB_1\to\BBB_2)\leadsto(f_2:\BBB_1\to\BBB_2)$,
        \begin{alignat*}{2}
                \sigma_0\sigma_0(\eta)&=\sigma_0(s_{1};t_{1})&=(s_{1};t_{1})&=\sigma_0(\eta)\\
                \sigma_0\tau_0(\eta)&=\sigma_0(s_{2};t_{2})&=(s_{2};t_{2})&=\tau_0(\eta)\\
                \tau_0\tau_0(\eta)&=\tau_0(s_{2};t_{2})&=(s_{2};t_{2})&=\tau_0(\eta)\\
                \tau_0\sigma_0(\eta)&=\tau_0(s_{1};t_{1})&=(s_{1};t_{1})&=\sigma_0(\eta)
        \end{alignat*}
\item[(GR 2)] Pour tous $$\eta_1:(f_1:\BBB_1\to\BBB'_1)\leadsto(f'_1:\BBB_1\to\BBB'_1)$$et
        $$\eta_2:(f_2:\BBB_2\to\BBB'_2)\leadsto(f'_2:\BBB_2\to\BBB'_2)$$ avec 
        $\sigma_0(\eta_2)=\tau_0(\eta_1)$,
        $$\sigma_0(\eta_2\boxtimes\eta_1)=(s_{1};t_{1})=\sigma_0(\eta_1)
        \qquad\text{et}\qquad
        \tau_0(\eta_2\boxtimes\eta_1)=(s'_{2};t'_{2})=\tau_0(\eta_2)$$
        car $\eta_2\boxtimes\eta_1$ va de $f'_{1}f_{1}:\BBB_{1}\to\BBB'_{2}$ 
        vers $f'_{2}f_{2}:\BBB_{1}\to\BBB'_{2}$.
\item[(GR 3)] Soient $\eta_1:f_1\leadsto f'_1$, $\eta_2:f_2\leadsto f'_2$ et 
    $\eta_3:f_3\leadsto f'_3$ trois éléments de $\star\GTrans$ tels que     
    $\sigma_0\left(\eta_2\right)=\tau_0\left(\eta_1\right)$ et 
    $\sigma_0\left(\eta_3\right)=\tau_0\left(\eta_2\right)$.
$$(\eta_3\boxtimes\eta_2)\boxtimes\eta_1
=\left((f'_3\eta_2)\circledast(\eta_3f_2)\right)\boxtimes\eta_1
=(f'_3f'_2\eta_1)\circledast\left((f'_3\eta_2f_1)\circledast(\eta_3f_2f_1)\right)$$
$$\eta_3\boxtimes(\eta_2\boxtimes\eta_1)
=\eta_3\boxtimes\left((f'_2\eta_1)\circledast(\eta_2f_1)\right)
=\left((f'_3f'_2\eta_1)\circledast(f'_3\eta_2f_1)\right)\circledast(\eta_3f_2f_1)$$
et
$$(f'_3f'_2\eta_1)\circledast\left((f'_3\eta_2f_1)\circledast(\eta_3f_2f_1)\right)
=\left((f'_3f'_2\eta_1)\circledast(f'_3\eta_2f_1)\right)\circledast(\eta_3f_2f_1)$$
d'après le théorème \ref{theo:gr:4.1}. 
\end{description}
On démontre de la même manière que $\left(\GTrans;\sigma_0;\tau_0;\boxdot\right)$ est 
un groupement.\QEDb
Ce théorème nous montre que l'ensemble des petites g-transformations est naturellement 
muni de trois structures de groupement. Dans les deux propositions qui viennent maintenant 
nous allons étudier de plus près les relations qui existent entre elles.
\begin{prop}\label{prop:gr:4.7}
Avec les notations introduites précédemment, on a :
$$\sigma_0\sigma_1=\sigma_0\qquad\sigma_0\tau_1=\sigma_0\qquad
\tau_0\sigma_1=\tau_0\qquad\tau_0\tau_1=\tau_0$$
$$\sigma_1\sigma_0=\sigma_0\qquad\sigma_1\tau_0=\tau_0\qquad
\tau_1\sigma_0=\sigma_0\qquad\tau_1\tau_0=\tau_0$$
\end{prop}
\Dem
Soit $\eta:(f_1:\BBB_1\to\BBB_2)\leadsto(f_2:\BBB_1\to\BBB_2)$ une g-transformation.
\begin{align*}
\sigma_0\sigma_1(\eta)&=\sigma_0(s_{2}f_{1};t_{2}f_{1})
=(s_{1};t_{1})=\sigma_0(\eta)&
\sigma_0\tau_1(\eta)&=\sigma_0(s_{2}f_{2};t_{2}f_{2})
=(s_{1};t_{1})=\sigma_0(\eta)\\
\tau_0\sigma_1(\eta)&=\tau_0(s_{2}f_1;t_{2}f_{1})
=(s_{2};t_{2})=\tau_0(\eta)&
\tau_0\tau_1(\eta)&=\tau_0(s_{2}f_2;t_{2}f_{2})
=(s_{2};t_{2})=\tau_0(\eta)\\
\sigma_1\sigma_0(\eta)&=\sigma_1(s_{1};t_{1})=(s_{1};t_{1})=\sigma_0(\eta)&
\sigma_1\tau_0(\eta)&=\sigma_1(s_{2};t_{2})=(s_{2};t_{2})=\tau_0(\eta)\\
\tau_1\sigma_0(\eta)&=\tau_1(s_{1};t_{1})=(s_{1};t_{1})=\sigma_0(\eta)&
\tau_1\tau_0(\eta)&=\tau_1(s_{2};t_{2})=(s_{2};t_{2})=\tau_0(\eta)
\QEDb
\end{align*}
Une conséquence immédiate de cette proposition est le corollaire suivant.
\begin{coro}\label{coro:gr:4.1}
On a
$$\sigma_0\sigma_1=\sigma_1\sigma_0\quad\text{et}\quad\tau_0\tau_1=\tau_1\tau_0$$
$$\tau_0\sigma_1=\sigma_1\tau_0\quad\text{et}\quad\sigma_0\tau_1=\tau_1\sigma_0$$
\end{coro}
Le problème dans le cas présent c'est que si les relations entre les applications
sources et buts sont extrèmement simples, celles entre les compositions le sont
beaucoup moins. Nous reviendrons dessus dans le chapitre \ref{chap:alex}.
%%%%%%%%%%%%%%%%%%%%%%%%%%%%%%%%%%%%%%%%%%%%%%%%%%
\chapter{Exemples de groupement :\\ Les chemins de Moore et leurs extensions}
\label{chap:moore}
En cherchant un exemple de groupement moins trivial que les catégories, on est 
de manière assez naturelle amené à chercher un exemple géométrique. Pourquoi 
ne pas regarder du côté de la topologie algébrique ? Le groupe fondamental qui 
fut introduit par Henri Poincaré à la fin du 19ème siècle et qui fut à l'origine 
de la topologie algébrique mérite quelques attentions. Il est basé sur les notions 
de chemins et d'homotopie. L'homotopie sert à identifier des chemins entre eux afin 
de pouvoir construire une composition qui vérifie les axiomes de groupe. 
Géométriquement cette composition est basée sur la juxtaposition des chemins.

Dans toute la suite de ce chapitre $X$ désignera un espace topologique.
%%%%%%%%%%%%%%%%%
%%% section 1 %%%
%%%%%%%%%%%%%%%%%
\section{Les chemins de Moore}
Un \emph{chemin de Moore}\index{chemin de Moore} est tout simplement une 
application continue $\gamma$ d'un intervalle $\left[0;d\right]$, où $d>0$, dans 
l'espace topologique $X$.

Notons $\BBM(X)$ l'ensemble des chemins de Moore de l'espace $X$. Dans la définition 
ci-dessus, le nombre $d$ correspond à la longueur de l'intervalle. Il peut 
aussi être vu comme la durée mise pour parcourir le chemin. Nous pouvons ainsi 
construire une application durée $\delta$ de $\BBM(X)$ dans $\left]0;+\infty\right[$
qui a tout chemin $\gamma$ associe sa durée de parcours $\delta_{\gamma}$.

Nous allons maintenant munir l'ensemble $\BBP(X)$ d'une structure de groupement.
Pour cela nous pouvons commencer par prendre pour applications source et but les applications $$\MFs:\BBP(X)\to\BBP(X)\quad\text{et}\quad\MFt:\BBP(X)\to\BBP(X)$$ définies de la 
manière suivante :
$$\MFs(\gamma)(u)=\gamma(0)\qquad\text{et}\qquad \MFt(\gamma)(t)=\gamma
(\delta_{\gamma})$$
pour tout $t\in\left[0;1\right]$. $\MFs(\gamma)$ et $\MFt(\gamma)$ sont des applications 
constantes, donc continues, de $\left[0;1\right]$ dans $X$. Ce sont des chemins 
de Moore de durée $1$.

Pour tout $t$ appartenant à l'intervalle $\left[0;1\right]$ et tout chemin 
$\gamma\in\BBP(X)$, on a
$$\MFs(\MFs(\gamma))(t)=\MFs(\gamma)(0)=\gamma(0)=\MFs(\gamma)(t)$$
$$\MFs(\MFt(\gamma))(t)=\MFt(\gamma)(0)=\gamma(\delta_{\gamma})=\MFt(\gamma)(t)$$
$$\MFt(\MFt(\gamma))(t)=\MFt(\gamma)(\delta_{\gamma})
=\gamma(\delta_{\gamma})=\MFt(\gamma)(t)$$
$$\MFt(\MFs(\gamma))(t)=\MFs(\gamma)(\delta_{\gamma})
=\gamma(0)=\MFs(\gamma)(t)$$
Par conséquent l'axiome (GR 1) est vérifié.

Ces deux seules applications ne sont pas suffisantes, il nous faut une composition.
Celle-ci correspond en fait à l'opération géométrique qui consiste à mettre 
bout-à-bout deux chemins à partir du moment où l'un commence là où fini l'autre.
Pour être plus précis, définissons la composition comme l'application $\OP$ de 
$\BBP(X)\times\BBP(X)$ dans $\BBP(X)$ donnée par la construction qui suit.
Soient $\gamma_1$ et $\gamma_2$ deux chemins de Moore.
\begin{itemize}
\item Si $\MFs(\gamma_2)=\MFt(\gamma_1)$, alors $\gamma_2\OP\gamma_1$ est la chemin 
de durée $\delta_{\gamma_1}+\delta_{\gamma_2}$ défini par
$$(\gamma_2\OP\gamma_1)(t)
\begin{cases}
\gamma_1(t)\quad&\text{pour }t\in\left[0;\delta_{\gamma_1}\right]\\
\gamma_2(t-\delta_{\gamma_1})\quad&\text{pour }t\in\left[\delta_{\gamma_1};
\delta_{\gamma_1}+\delta_{\gamma_2}\right]
\end{cases}$$
\item Sinon $\gamma_2\OP\gamma_1=\gamma_1$.
\end{itemize}  
Dans le premier cas, comme $\MFs(\gamma_2)=\MFt(\gamma_1)$, on a $\gamma_1(\delta
_{\gamma_1})=\gamma_2(0)=\gamma_2(\delta_{\gamma_1}-\delta_{\gamma_1})$. De plus les 
intervalles $\left[0;\delta_{\gamma_1}\right]$ et $\left[\delta_{\gamma_1};
\delta_{\gamma_1}+\delta_{\gamma_2}\right]$ sont fermés. Par conséquent 
$\gamma_2\OP\gamma_1$ est bien continue de $\left[0;\delta_{\gamma_1}+
\delta_{\gamma_2}\right]$ dans $X$. En d'autres mots c'est un chemin de Moore.
Le second cas donne trivialement un chemin de Moore.

Il ne nous reste plus qu'à vérifier les axiomes (GR 2) et (GR 3). Soient $\gamma_1$ 
et $\gamma_2$ deux chemins de Moore tels que $\MFs(\gamma_2)=\MFt(\gamma_1)$. Pour 
tout $t\in\left[0;\delta_{\gamma_1}+\delta_{\gamma_2}\right]$, les calculs
$$\MFs(\gamma_2\OP\gamma_1)(t)=(\gamma_2\OP\gamma_1)(0)=\gamma_1(0)=\MFs(\gamma_1)(t)$$
$$\MFt(\gamma_2\OP\gamma_1)(t)=(\gamma_2\OP\gamma_1)(\delta_{\gamma_1}+
\delta_{\gamma_2})=\gamma_2(\delta_{\gamma_1}+\delta_{\gamma_2}-\delta_{\gamma_1})
=\gamma_2(\delta_{\gamma_2})=\MFt(\gamma_2)(t)$$
montrent que l'axiome (GR 2) est satisfait.

Si $\gamma_1$, $\gamma_2$ et $\gamma_3$ sont trois chemins de Moore tels que 
$\MFs(\gamma_3)=\MFt(\gamma_2)$ et $\MFs(\gamma_2)=\MFt(\gamma_1)$, alors
$$(\gamma_3\OP\gamma_2)\OP\gamma_1(t)=
\begin{cases}
\gamma_1(t)\quad&\text{pour }0\leq t\leq
\delta_{\gamma_1}\\
(\gamma_3\OP\gamma_2)(t-\delta_{\gamma_1})\quad&\text{pour }
\delta_{\gamma_1}\leq t\leq
\delta_{\gamma_1}+\delta_{\gamma_2}+\delta_{\gamma_3}
\end{cases}$$
et
$$\gamma_3\OP(\gamma_2\OP\gamma_1)(t)=
\begin{cases}
\gamma_2\OP\gamma_1(t)\quad&\text{pour }0\leq t\leq
\delta_{\gamma_1}+\delta_{\gamma_2}\\
\gamma_3(t-\delta_{\gamma_1}-\delta_{\gamma_2})\quad&\text{pour }
\delta_{\gamma_1}+\delta_{\gamma_2}\leq t\leq
\delta_{\gamma_1}+\delta_{\gamma_2}+\delta_{\gamma_3}
\end{cases}$$
D'où
$$(\gamma_3\OP\gamma_2)\OP\gamma_1(t)=\gamma_3\OP(\gamma_2\OP\gamma_1)(t)=
\begin{cases}
\gamma_1(t)\quad&\text{pour }0\leq t\leq
\delta_{\gamma_1}\\
\gamma_2(t-\delta_{\gamma_1})\quad&\text{pour }
\delta_{\gamma_1}\leq t\leq
\delta_{\gamma_1}+\delta_{\gamma_2}\\
\gamma_3(t-\delta_{\gamma_1}-\delta_{\gamma_2})\quad&\text{pour }
\delta_{\gamma_1}+\delta_{\gamma_2}\leq t\leq
\delta_{\gamma_1}+\delta_{\gamma_2}+\delta_{\gamma_3}
\end{cases}$$
Ainsi l'axiome (GR 3) est vérifié. Finalement
\begin{theo}\label{theo:moore:1.1}
Le quadruplet $\left(\BBP(X);
\MFs;\MFt;\OP\right)$ est un groupement.
\end{theo}
%%%%%%%%%%%%%%%%%
%%% section 2 %%%
%%%%%%%%%%%%%%%%%
\section{Les surfaces de Moore}
Une \emph{surface de Moore}\index{surface de Moore} est une application continue
de $\left[0;d_1\right]\times\left[0;d_2\right]$, où $d_1$ et $d_2$ sont des réels 
strictement positifs, dans l'espace $X$.

Il faut bien remarquer que le mot surface est ici utilisé de manière abusive. S'il 
est entendu que le produit $\left[0;d_1\right]\times\left[0;d_2\right]$ est une surface, 
il n'en va pas nécessairement de même pour son image dans $X$.

En se référant à ce que nous avons dit dans la section précédente, nous pouvons 
voir $d_1$ comme la durée de parcours dans la première direction et $d_2$ comme celle 
dans la seconde direction. Ainsi si nous notons $\BBS(X)$ l'ensemble des surfaces 
de Moore de l'espace $X$, alors nous avons naturellement deux applications de durée
$\delta_1$ et $\delta_2$ de $\BBS(X)$ dans $\left]0;+\infty\right[$ :
\begin{align*}
\delta_1(\gamma)&=d_1=\text{durée dans la première direction}\\
\delta_2(\gamma)&=d_2=\text{durée dans la deuxième direction}
\end{align*}
pour $\gamma\in\BBS(X)$.

Nous pouvons munir l'ensemble $\BBS(X)$ de deux structures de groupement. Considérons 
$\MFs_1$ et $\MFs_2$ les deux applications définies, sur $\BBS(X)$ à valeur dans 
lui-même, par les formules suivantes :
$$\MFs_1(\gamma)(u;v)=\gamma(0;v)\qquad\text{pour }(u;v)\in\left[0;1\right]\times
\left[0;\delta_2(\gamma)\right]$$
$$\MFs_2(\gamma)(u;v)=\gamma(u;0)\qquad\text{pour }(u;v)\in\left[0;\delta_1
(\gamma)\right]\times\left[0;1\right]$$
De même on peut définir deux applications $\MFt_i:\BBS(X)\to\BBS(X)$, $i=1,2$, de la 
manière suivante :
$$\MFt_1(\gamma)(u;v)=\gamma(\delta_1(\gamma);v)
\qquad\text{pour }(u;v)\in\left[0;1\right]\times\left[0;\delta_2(\gamma)\right]$$
$$\MFt_2(\gamma)(u;v)=\gamma(u;\delta_2(\gamma))
\qquad\text{pour }(u;v)\in\left[0;\delta_1(\gamma)\right]\times\left[0;1\right]$$
Les applications sources et buts étant données, nous pouvons maintenant nous 
intéresser aux compositions. Il y aura une composition $\OP_1$  suivant la première 
direction et une autre $\OP_2$ suivant la seconde. Pour deux surfaces de Moore 
$\gamma_1$ et $\gamma_2$,
\begin{itemize}
\item si $\MFs_1(\gamma_2)=\MFt_1(\gamma_1)$, on pose
$$\gamma_2\OP_1\gamma_1(u;v)=
\begin{cases}
\gamma_1(u;v)\quad&\text{pour }(u;v)\in\left[0;\delta_1(\gamma_1)\right]
\times\left[0;\delta_2(\gamma_1)\right]\\
\gamma_2(u-\delta_1(\gamma_1);v)\quad&\text{pour }
(u;v)\in\left[\delta_1(\gamma_1);\delta_1(\gamma_1)+\delta_1(\gamma_2)
\right]\times\left[0;\delta_2(\gamma_1)\right]
\end{cases}$$
(Cette définition est valide. En effet,
comme $$\MFs_1(\gamma_2):\left[0;1\right]\times\left[0;\delta_2(\gamma_2)\right]
\to\BBS(X)\quad\text{et}\quad\MFt_1(\gamma_1):\left[0;1\right]\times\left[0;\delta_2
(\gamma_1)\right]\to\BBS(X)$$ on a $\delta_2(\gamma_1)=\delta_2(\gamma_2)$.)
\item si $\MFs_2(\gamma_2)=\MFt_2(\gamma_1)$, on pose
$$\gamma_2\OP_2\gamma_1(u;v)=
\begin{cases}
\gamma_1(u;v)\quad&\text{pour }(u;v)\in\left[0;\delta_1(\gamma_1)\right]
\times\left[0;\delta_2(\gamma_1)\right]\\
\gamma_2(u;v-\delta_2(\gamma_1))\quad&\text{pour }
(u;v)\in\left[0;\delta_1(\gamma_1)\right]\times\left[\delta_2
(\gamma_1);\delta_2(\gamma_1)+\delta_2(\gamma_2)\right]
\end{cases}$$
(Comme précédemment on montre $\delta_2(\gamma_1)=\delta_2(\gamma_2)$.)
\item sinon, on pose $\gamma_2\OP_1\gamma_1=\gamma_2\OP_2\gamma_1=\gamma_1$
\end{itemize} 
\begin{theo}\label{theo:moore:2.1}
Les quadruplets $\left(\BBS(X);\MFs_1;\MFt_1;\OP_1\right)$ et 
$\left(\BBS(X);\MFs_2;\MFt_2;\OP_2\right)$ sont des groupements.
\end{theo}
\Dem 
\begin{description}
\item[(GR 1)] Soient $\gamma\in\BBS(X)$ et $(u;v)\in\left[0;1\right]\times
\left[0;\delta_2(\gamma)\right]$.
$$\MFs_1\MFs_1(\gamma)(u;v)=\MFs_1(\gamma(0;v))
=\gamma(0;v)=\MFs_1(\gamma)(u;v)$$
$$\MFs_1\MFt_1(\gamma)(u;v)=\MFs_1(\gamma(\delta_1(\gamma);v))
=\gamma(\delta_1(\gamma);v)=\MFt_1(\gamma)(u;v)$$
$$\MFt_1\MFt_1(\gamma)(u;v)=\MFt_1(\gamma(\delta_1(\gamma);v))
=\gamma(\delta_1(\gamma);v)=\MFt_1(\gamma)(u;v)$$
$$\MFt_1\MFs_1(\gamma)(u;v)=\MFs_1(\gamma(0;v))
=\gamma(0;v)=\MFs_1(\gamma)(u;v)$$
\item[(GR 2)] Soient $\gamma_1$ et $\gamma_2$ deux surfaces de Moore telles que 
$\MFs_1(\gamma_2)=\MFt_1(\gamma_1)$. Pour $(u;v)\in\left[0;1\right]\times
\left[0;\delta_2(\gamma_1)\right]$,
$$\MFs_1(\gamma_2\OP_1\gamma_1)(u;v)
=\gamma_2\OP_1\gamma_1(0;v)=\gamma_1(0;v)=\MFs_1(\gamma_1)(u;v)$$
$$\MFt_1(\gamma_2\OP_1\gamma_1)(u;v)
=\gamma_2\OP_1\gamma_1(\delta_1(\gamma_1)+\delta_1(\gamma_2);v)
=\gamma_2(\delta_1(\gamma_2);v)=\MFt_1(\gamma_2)(u;v)$$
\item[(GR 3)] Soient $\gamma_1$, $\gamma_2$ et $\gamma_3$ deux surfaces de Moore 
telles que $\MFs_1(\gamma_3)=\MFt_1(\gamma_2)$ et
$\MFs_1(\gamma_2)=\MFt_1(\gamma_1)$. Pour $(u;v)\in\left[0;\delta_1(\gamma_1)+
\delta_1(\gamma_2)+\delta_1(\gamma_3)\right]\times
\left[0;\delta_2(\gamma_1)\right]$,
$$(\gamma_3\OP\gamma_2)\OP\gamma_1(u;v)=
\begin{cases}
\gamma_1(t)\quad&\text{pour }0\leq u\leq
\delta_{\gamma_1}\\
(\gamma_3\OP\gamma_2)(u-\delta_{\gamma_1};v)\quad&\text{pour }
\delta_{\gamma_1}\leq u\leq
\delta_{\gamma_1}+\delta_{\gamma_2}+\delta_{\gamma_3}
\end{cases}$$
et
$$\gamma_3\OP(\gamma_2\OP\gamma_1)(u;v)=
\begin{cases}
\gamma_2\OP\gamma_1(u;v)\quad&\text{pour }0\leq u\leq
\delta_{\gamma_1}+\delta_{\gamma_2}\\
\gamma_3(u-\delta_{\gamma_1}-\delta_{\gamma_2};v)\quad&\text{pour }
\delta_{\gamma_1}+\delta_{\gamma_2}\leq u\leq
\delta_{\gamma_1}+\delta_{\gamma_2}+\delta_{\gamma_3}
\end{cases}$$
D'où
\begin{align*}
(\gamma_3\OP\gamma_2)\OP\gamma_1(u;v)\\
&=\gamma_3\OP(\gamma_2\OP\gamma_1)(u;v)\\
&=
\begin{cases}
\gamma_1(u;v)\quad&\text{pour }0\leq u\leq
\delta_{\gamma_1}\\
\gamma_2(u-\delta_{\gamma_1};v)\quad&\text{pour }
\delta_{\gamma_1}\leq u\leq
\delta_{\gamma_1}+\delta_{\gamma_2}\\
\gamma_3(u-\delta_{\gamma_1}-\delta_{\gamma_2};v)\quad&\text{pour }
\delta_{\gamma_1}+\delta_{\gamma_2}\leq u\leq
\delta_{\gamma_1}+\delta_{\gamma_2}+\delta_{\gamma_3}
\end{cases}
\end{align*}
Cela montre que $\left(\BBS(X);\MFs_1;\MFt_1;\OP_1\right)$ est un groupement. La 
démonstration est semblable pour $\left(\BBS(X);\MFs_2;\MFt_2;\OP_2\right)$.\QED
\end{description} 
\'Etudions de plus près les liens qui unissent ces deux structures de groupements.
\begin{prop}\label{prop:moore:2.1}
Soient $\left(\BBS(X);\MFs_1;\MFt_1;\OP_1\right)$ et $\left(\BBS(X);\MFs_2;\MFt_2;\OP_2\right)$
les deux structures de groupement sur $\BBS(X)$ introduites dans le théorème
\ref{theo:moore:2.1}. On a les égalités
$$\MFs_1\MFs_2=\MFs_2\MFs_1\qquad\MFt_1\MFt_2=\MFt_2\MFt1$$
$$\MFs_1\MFt_2=\MFt_2\MFs_1\qquad\MFs_2\MFt_1=\MFt_1\MFs_2$$
\end{prop}
\Dem
Soient $\gamma\in\BBS(X)$ et $(u;v)\in\left[0;1\right]\times\left[0;1\right]$.
$$\MFs_1\MFs_2\left(\gamma\right)(u;v)=\MFs_1\left(\gamma\right)(u;0)=\gamma(0;0)
\quad\text{et}\quad\MFs_2\MFs_1\left(\gamma\right)(u;v)=\MFs_2\left(\gamma\right)(0;v)
=\gamma(0;0)$$
$$\MFs_1\MFt_2\left(\gamma\right)(u;v)=\MFs_1\left(\gamma\right)(u;1)=\gamma(0;1)
\quad\text{et}\quad\MFt_2\MFs_1\left(\gamma\right)(u;v)=\MFt_2\left(\gamma\right)(0;v)
=\gamma(0;1)$$
$$\MFt_1\MFs_2\left(\gamma\right)(u;v)=\MFt_1\left(\gamma\right)(u;0)=\gamma(1;0)
\quad\text{et}\quad\MFs_2\MFt_1\left(\gamma\right)(u;v)=\MFs_2\left(\gamma\right)(1;v)
=\gamma(1;0)$$
$$\MFt_1\MFt_2\left(\gamma\right)(u;v)=\MFt_1\left(\gamma\right)(u;1)=\gamma(1;1)
\quad\text{et}\quad\MFt_2\MFt_1\left(\gamma\right)(u;v)=\MFt_2\left(\gamma\right)(1;v)
=\gamma(1;1)$$
\QEDb
Cette proposition est à rapprocher du corollaire \ref{coro:gr:4.1}.
\begin{prop}\label{prop:moore:2.2}
Soient $\gamma_i$, pour
$i=1,\ldots,4$, quatres éléments de $\BBS(X)$.
Si $\MFs_2(\gamma_2)=\MFt_2(\gamma_1)$, $\MFs_2(\gamma_4)=\MFt_2(\gamma_3)$,
$\MFs_1(\gamma_3)=\MFt_1(\gamma_1)$ et $\MFs_1(\gamma_4)=\MFt_1(\gamma_2)$
alors les surfaces de Moore
$$(\gamma_4\OP_2\gamma_3)\OP_1(\gamma_2\OP_2\gamma_1)\quad\text{et}\quad
(\gamma_4\OP_1\gamma_2)\OP_2(\gamma_3\OP_1\gamma_1)$$
existent et sont égales :
$$(\gamma_4\OP_2\gamma_3)\OP_1(\gamma_2\OP_2\gamma_1)=
(\gamma_4\OP_1\gamma_2)\OP_2(\gamma_3\OP_1\gamma_1)$$
\end{prop}
\Dem
Commençons par remarquer que
 les quatre conditions imposées impliquent $$\delta_1(\gamma_1)=
\delta_1(\gamma_2),\qquad\delta_1(\gamma_3)=\delta_1(\gamma_4)$$et
 $$\delta_2(\gamma_1)=
\delta_2(\gamma_3),\qquad\delta_2(\gamma_2)=\delta_2(\gamma_4)$$

Il est clair que la surface de Moore $(\gamma_4\OP_2\gamma_3)\OP_1(\gamma_2\OP_2\gamma_1)$
existe si les conditions $$\MFs_2(\gamma_2)=\MFt_2(\gamma_1),\quad
\MFs_2(\gamma_4)=\MFt_2(\gamma_3)\quad\text{et}\quad\MFs_1(\gamma_4\OP_2\gamma_3)=
\MFt_1(\gamma_2\OP_2\gamma_1)$$ sont satisfaites. Par hypothèse, seule la dernière
mérite une justification.\\
Pour tout $(u;v)\in\left[0;\delta_1(\gamma_1)\right]\times\left[0;
\delta_2(\gamma_1)+\delta_2(\gamma_2)\right]$,
$$\gamma_2\OP_2\gamma_1(u;v)=
\begin{cases}
\gamma_1(u;v) &\text{si }0\leq v\leq\delta_2(\gamma_1)\\
\gamma_2(u;v-\delta_2(\gamma_1)) &\text{si }
\delta_2(\gamma_1)\leq v\leq\delta_2(\gamma_1)+\delta_2(\gamma_2)
\end{cases}$$
et pour tout $(u;v)\in\left[0;\delta_1(\gamma_3)\right]\times\left[0;
\delta_2(\gamma_1)+\delta_2(\gamma_2)\right]$,
$$\gamma_4\OP_2\gamma_3(u;v)=
\begin{cases}
\gamma_3(u;v) &\text{si }0\leq v\leq\delta_2(\gamma_1)\\
\gamma_4(u;v-\delta_2(\gamma_1)) &\text{si }
\delta_2(\gamma_1)\leq v\leq\delta_2(\gamma_1)+\delta_2(\gamma_2)
\end{cases}$$
Ainsi, pour tout $(u;v)\in\left[0;1\right]\times
\left[0;\delta_2(\gamma_1)+\delta_2(\gamma_2)\right]$, on trouve
$$\MFt_1\left(\gamma_2\OP_2\gamma_1\right)(u;v)=
\begin{cases}
\gamma_1(\delta_1(\gamma_1);v) &\text{si }0\leq v\leq\delta_2(\gamma_1)\\
\gamma_2(\delta_1(\gamma_1);v-\delta_2(\gamma_1)) &\text{si }
\delta_2(\gamma_1)\leq v\leq\delta_2(\gamma_1)+\delta_2(\gamma_2)
\end{cases}$$
et
$$\MFs_1\left(\gamma_4\OP_2\gamma_3\right)(u;v)=
\begin{cases}
\gamma_3(0;v) &\text{si }0\leq v\leq\delta_2(\gamma_3)\\
\gamma_2(0;v-\delta_2(\gamma_1)) &\text{si }
\delta_2(\gamma_3)\leq v\leq\delta_2(\gamma_3)+\delta_2(\gamma_4)
\end{cases}$$
Or les conditions $\MFs_1(\gamma_3)=\MFt_1(\gamma_1)$ et $\MFs_1(\gamma_4)=\MFt_1(\gamma_2)$
s'écrivent de la manière suivante
$$\begin{cases}
\gamma_3(0;v)=\gamma_1(\delta_1(\gamma_1);v) &\text{pour }0\leq v\leq\delta_2(\gamma_1)\\
\gamma_4(0;v)=\gamma_2(\delta_1(\gamma_1);v) &\text{pour }0\leq v\leq\delta_2(\gamma_2)
\end{cases}$$
Par conséquent, l'égalité $\MFt_1\left(\gamma_2\OP_2\gamma_1\right)=
\MFs_1\left(\gamma_4\OP_2\gamma_3\right)$ est bien vérifiée et, pour tout $$(u;v)\in
\left[0;\delta_1(\gamma_1)+\delta_1(\gamma_3)\right]\times
\left[0;\delta_2(\gamma_1)+\delta_2(\gamma_2)\right]$$
$(\gamma_4\OP_2\gamma_3)\OP_1(\gamma_2\OP_1\gamma_1)(u;v)$ est égal à
$$\begin{cases}
\gamma_1(u;v) &\text{pour }(u;v)\in\left[0;\delta_1(\gamma_1)\right]\times
\left[0;\delta_2(\gamma_1)\right]\\
\gamma_3(u-\delta_1(\gamma_1);v) &\text{pour }
(u;v)\in\left[\delta_1(\gamma_1);\delta_1(\gamma_1)+
\delta_1(\gamma_3)\right]\times\left[0;\delta_2(\gamma_1)\right]\\
\gamma_2(u;v-\delta_2(\gamma_1)) &\text{pour }
(u;v)\in\left[0;\delta_1(\gamma_1)\right]\times
\left[\delta_2(\gamma_1);\delta_2(\gamma_1)+\delta_2(\gamma_2)\right]\\
\gamma_4(u-\delta_1(\gamma_1);v-\delta_2(\gamma_2)) &\text{pour }
(u;v)\in\left[\delta_1(\gamma_1);\delta_1(\gamma_1)+
\delta_1(\gamma_3)\right]\times
\left[\delta_2(\gamma_1);\delta_2(\gamma_1)+\delta_2(\gamma_2)\right]
\end{cases}$$

De même, il est clair que la surface de Moore $(\gamma_4\OP_1\gamma_2)\OP_2(\gamma_3\OP_1\gamma_1)$
existe si les conditions $\MFs_1(\gamma_4)=\MFt_1(\gamma_2)$,
$\MFs_1(\gamma_3)=\MFt_1(\gamma_1)$ et $\MFs_2(\gamma_4\OP_1\gamma_2)=
\MFt_2(\gamma_3\OP_1\gamma_1)$ sont satisfaites. Les deux premières étant évidentes,
intéressons nous à la troisième.\\
Pour tout $(u;v)\in\left[0;\delta_1(\gamma_1)+\delta_1(\gamma_3)\right]\times\left[0;
\delta_2(\gamma_1)\right]$,
$$\gamma_3\OP_1\gamma_1(u;v)=
\begin{cases}
\gamma_1(u;v) &\text{si }0\leq u\leq\delta_1(\gamma_1)\\
\gamma_3(u-\delta_1(\gamma_1);v) &\text{si }
\delta_1(\gamma_1)\leq u\leq\delta_1(\gamma_1)+\delta_1(\gamma_3)
\end{cases}$$
et pour tout $(u;v)\in\left[0;\delta_1(\gamma_1)+\delta_1(\gamma_3)\right]\times\left[0;
\delta_2(\gamma_2)\right]$,
$$\gamma_4\OP_1\gamma_2(u;v)=
\begin{cases}
\gamma_2(u;v) &\text{si }0\leq u\leq\delta_1(\gamma_1)\\
\gamma_4(u-\delta_1(\gamma_1);v) &\text{si }
\delta_1(\gamma_1)\leq u\leq\delta_1(\gamma_1)+\delta_1(\gamma_3)
\end{cases}$$
Pour tout $(u;v)\in\left[0;\delta_1(\gamma_1)+\delta_1(\gamma_3)\right]
\times\left[0;1\right]$, on a ainsi
$$\MFt_2\left(\gamma_3\OP_1\gamma_1\right)(u;v)=
\begin{cases}
\gamma_1(u;\delta_2(\gamma_1))
&\text{si }0\leq u\leq\delta_1(\gamma_1)\\
\gamma_3(u-\delta_1(\gamma_1);\delta_2(\gamma_1)) &\text{si }
\delta_1(\gamma_1)\leq u\leq\delta_1(\gamma_1)+\delta_1(\gamma_3)
\end{cases}$$
et
$$\MFs_2\left(\gamma_4\OP_1\gamma_1\right)(u;v)=
\begin{cases}
\gamma_3(u;0) &\text{si }0\leq u\leq\delta_1(\gamma_1)\\
\gamma_2(u-\delta_1(\gamma_1);0) &\text{si }
\delta_1(\gamma_1)\leq u\leq\delta_1(\gamma_1)+\delta_1(\gamma_3)
\end{cases}$$
Comme les conditions $\MFs_2(\gamma_2)=\MFt_2(\gamma_1)$ et
$\MFs_2(\gamma_4)=\MFt_2(\gamma_3)$ s'écrivent de la manière suivante
$$\begin{cases}
\gamma_2(u;O)=\gamma_1(u;\delta_2(\gamma_1)) &\text{pour }0\leq u\leq\delta_1(\gamma_1)\\
\gamma_4(u;0)=\gamma_3(u;\delta_2(\gamma_1)) &\text{pour }0\leq u\leq\delta_1(\gamma_3)
\end{cases}$$
Par conséquent, l'égalité $$\MFt_2\left(\gamma_3\OP_1\gamma_1\right)=
\MFs_2\left(\gamma_4\OP_1\gamma_2\right)$$ est satisfaite et, pour tout $$(u;v)\in
\left[0;\delta_1(\gamma_1)+\delta_1(\gamma_3)\right]\times
\left[0;\delta_2(\gamma_1)+\delta_2(\gamma_2)\right]$$
$(\gamma_4\OP_1\gamma_2)\OP_2(\gamma_3\OP_1\gamma_1)(u;v)$ est égal à
$$\begin{cases}
\gamma_1(u;v) &\text{pour }(u;v)\in\left[0;\delta_1(\gamma_1)\right]\times
\left[0;\delta_2(\gamma_1)\right]\\
\gamma_3(u-\delta_1(\gamma_1);v) &\text{pour }
(u;v)\in\left[\delta_1(\gamma_1);\delta_1(\gamma_1)+
\delta_1(\gamma_3)\right]\times\left[0;\delta_2(\gamma_1)\right]\\
\gamma_2(u;v-\delta_2(\gamma_1)) &\text{pour }
(u;v)\in\left[0;\delta_1(\gamma_1)\right]\times
\left[\delta_2(\gamma_1);\delta_2(\gamma_1)+\delta_2(\gamma_2)\right]\\
\gamma_4(u-\delta_1(\gamma_1);v-\delta_2(\gamma_2)) &\text{pour }
(u;v)\in\left[\delta_1(\gamma_1);\delta_1(\gamma_1)+
\delta_1(\gamma_3)\right]\times
\left[\delta_2(\gamma_1);\delta_2(\gamma_1)+\delta_2(\gamma_2)\right]
\end{cases}$$
D'où le résultat.\QED
%%%%%%%%%%%%%%%%%
%%% section 3 %%%
%%%%%%%%%%%%%%%%%
\section{$I$-espaces de Moore}
Soit $I$ un ensemble non vide. On appelle \emph{$I$-espace de Moore}
\index{$I$-espace de Moore} de $X$ toute application continue
$\eta:\prod_{i\in I}\left[0;d_i\right]\to X$ où $\left(d_i\right)_{i\in I}$ est une famille
de nombres réels strictement positifs. Ici l'espace $\prod_{i\in I}\left[0;d_i\right]$ est
muni de la topologie produit et tous les $\left[0;d_i\right]$ sont munis de celle induite par
la topologie usuelle de $\BBR$.

Pour $I=\{1\}$, nous retouvons les chemins de Moore et, pour $I=\{1;2\}$, nous obtenons
les surfaces de Moore.

Notons $\BBM_I(X)$ l'ensemble des $I$-espaces de Moore de $X$. Comme nous l'avons fait pour
les chemins et les surfaces nous pouvons définir une
application de durée $\delta_i:\BBE_I(X)\to\left]0;+\infty\right[$, pour chaque
$i\in I$, en posant $$\delta_i(\eta)=d_i$$
pour tout $\eta:\prod_{i\in I}\left[0;d_i\right]\to X$ appartenant à $\BBM_I(X)$. Nous
dirons que $\delta_i(\eta)$ est la durée de $\eta$ dans la $i$-ième direction.

Maintenant nous allons construire $\Card(I)$ structures de groupements sur l'ensemble
$\BBE_I(X)$. Mais avant, définissons une nouvelle notation. Si $(u_i)_I$ est une famille
de nombres, nous écrirons $(u_j\lhd_i x)_I$  pour indiquer que l'on s'intéresse à la
famille $(u_i)_I$ dans laquelle le $i$-ième terme $u_i$ est remplacé par le nombre $x$.
De manière plus générale, nous noterons $(u_k\lhd_i x \lhd_j y)_I$, pour $i\neq j$, la
famille $(u_i)_I$ dont les termes $u_i$ et $u_j$ sont remplacés respectivement par
$x$ et $y$. Maintenant, pour chaque $i\in I$, posons
\begin{itemize}
\item $\MFs_i:\BBE_I(X)\to\BBE_I(X)$ avec $\MFs_i(\eta)(u_j)_I=\eta(u_i\lhd_i 0)_I$
pour tout $(u_j)_I\in\prod_{j\in I}\left[0;\delta_j(\eta)\lhd_i 0\right]$
\item $\MFt_i:\BBE_I(X)\to\BBE_I(X)$ avec $\MFt_i(\eta)(u_j)_I=\eta(u_i\lhd_i
\delta_i(\eta))_I$
pour tout $(u_j)_I\in\prod_{j\in I}\left[0;\delta_j(\eta)\lhd_i 0\right]$
\item $\OP_i:\BBE_I(X)\times\BBE_I(X)\to\BBE_I(X)$ avec
\begin{itemize}
\item si $\MFt_i(\eta_2)=\MFs_i(\eta_1)$, alors
$\eta_2\OP_i\eta_1(u_j)_I$ est égal à
$$\begin{cases}
\eta_1(u_j)_I, &(u_j)\in\prod_{j\in I}
\left[0;\delta_j(\eta_1)\right]\\
\eta_2(u_j\lhd_i u_i-\delta_i(\eta_1))_I,&(u_j)\in\prod_{j\in I}
\left[0\lhd_i\delta_i(\eta_1);\delta_j(\eta_2)\lhd_i
\delta_i(\eta_1)+\delta_i(\eta_2)\right]
\end{cases}$$
\item sinon $\eta_2\OP_i\eta_1=\eta_1$
\end{itemize}
\end{itemize}
Des démonstrations semblables à celle du théorème \ref{theo:moore:2.1} et des propositions
\ref{prop:moore:2.1} et \ref{prop:moore:2.2}, nous donne
\begin{theo}\label{theo:moore:3.1}
Pour chaque $i\in I$, le quadruplet $\left(\BBE_I(X);\MFs_i;\MFt_i;\OP_i\right)$ est
un groupement.
\end{theo}
\begin{prop}\label{prop:moore:3.1}
Pour tous $i$ et $j$ de $I$, $i\neq j$, on a
$\qquad\MFs_i\MFs_j=\MFs_j\MFs_i,\qquad\MFt_i\MFt_j=\MFt_j\MFt_i,\qquad
\MFs_j\MFt_i=\MFt_i\MFs_j$
\end{prop}
\begin{prop}\label{prop:moore:3.2}
Soient $\eta_i$, pour
$i=1,\ldots,4$, quatres éléments de $\BBE_I(X)$ et $i$, $j$ deux éléments distincts de $I$.
Si $\MFs_i(\eta_2)=\MFt_i(\eta_1)$, $\MFs_i(\eta_4)=\MFt_i(\eta_3)$,
$\MFs_j(\eta_3)=\MFt_j(\eta_1)$ et $\MFs_i(\eta_4)=\MFt_j(\eta_2)$
alors les $I$-espaces de Moore
$$(\eta_4\OP_i\eta_3)\OP_j(\eta_2\OP_j\eta_1)\quad\text{et}\quad
(\eta_4\OP_j\eta_2)\OP_i(\eta_3\OP_j\eta_1)$$
existent et sont égaux :
$$(\eta_4\OP_i\eta_3)\OP_j(\eta_2\OP_i\eta_1)=
(\eta_4\OP_j\eta_2)\OP_i(\eta_3\OP_j\eta_1)$$
\end{prop}

%%%%%%%%%%%%%%%%%%%%%%%%%%%%%%%%%%%%%%%%%%%%%%%%%%
\chapter{Les groupements d'Alexandroff}\label{chap:alex}

Jusqu'à maintenant, nous avons cherché à calquer la théorie des groupements sur celle des 
catégories. Mais nous nous sommes alors trouvé dans l'impossibilité de définir un 
équivalent satisfaisant à la notion de transformation naturelle. Dans ce chapitre, objectif 
est de spécialiser très légèrement la notion de groupement afin de pouvoir construire 
des transformations \Dire{naturelles}.
%%%%%%%%%%%
%%% section 1 %%%
%%%%%%%%%%%
\section{Définition d'un groupement d'Alexandroff}\label{sec:alex:1}
Un \emph{groupement d'Alexandroff}\index{groupement!groupement d'Alexandroff} 
est un groupement $(\BBB;s;t;\OP)$ qui possède  un élément particulier $\alpha$ 
vérifiant les conditions suivantes :
\begin{description}
\item[(GALEX 1)] Pour tout $x\in\BBB$ tel que $x\neq\alpha$, alors
         $s(x)\neq\alpha$ et $t(x)\neq\alpha$.
\item[(GALEX 2)] Pour tout $x\in\BBB$, on a
         $x\OP\alpha=\alpha\OP x=x$.
\end{description}
$\alpha$ sera parfois appelé l'\emph{alexis}\index{alexis} du groupement 
d'Alexandroff.
% remarque
\begin{rema}\label{rema:alex:1.1}
Bien qu'il semble en vérifier toutes les propriétés, il ne faut pas croire que l'alexis 
soit une simple identité. La condition (GALEX 2) est plus forte que celle 
imposée aux identités. En effet, dans ce cas, la composée de l'alexis avec un autre 
élément quelconque laisse toujours ce dernier inchangé, alors que l'on est sûr que 
la composée d'une identité avec un autre élément ne le laisse invariant que si 
l'identité est la source ou le but de celui-ci. En résumé, un alexis est une super 
identité. De plus on a
\end{rema}
% lemme
\begin{lemm}\label{lemm:alex:1.1}
Un groupement d'Alexandroff possède un unique alexis.
\end{lemm}
\Dem\ Supposons que $\alpha$ et $\alpha'$ soient deux alexis. Comme $\alpha$ 
satisfait la condition (GALEX 2), on a
$$\alpha\OP\alpha'=\alpha'$$
De même, $\alpha'$ vérifiant aussi cette condition, on a
$$\alpha\OP\alpha'=\alpha$$
Donc $\alpha=\alpha'$.\QEDb
C'est pour cette raison que nous écrirons $(\BBB;\alpha)$ pour désigner un 
groupement d'Alexandroff d'alexis $\alpha$.
% remarque
\begin{rema}\label{rema:alex:1.2}
Bien sûr  nous aurions pu parler de groupements pointés, mais en topologie un 
espace pointé est simplement un espace topologique dans lequel on a choisi un 
élément sans imposé de condition sur celui-ci. En fait nous allons voir plus tard,
qu'à tout groupement, on peut associer, de manière très naturelle, un groupement
d'Alexandroff. La ressemblance avec le compactifié d'Alexandroff d'un espace 
topologique localement compact nous a semblé suffisamment forte pour utiliser 
ce nom. De plus les premiers exemples simples sont construit à partir d'un espace 
topologique.
\end{rema}
% exemple
\begin{exem}\label{exem:alex:1.1}
Soit $X$ un espace topologique. Nous savons que l'ensemble des ouverts 
de $X$ est stable par réunion finie et intersection finie.
De plus il possède deux éléments particuliers : l'ensemble vide $\emptyset$ et 
$X$. Ce qui nous permet de définir deux structures de groupements d'Alexandroff:
\begin{enumerate}
\item $\left((\Ouv(X);\stackrel{\cup}{s};\stackrel{\cup}{t};
         \cup);\emptyset\right)$ avec
         \begin{enumerate}
         \item $\Ouv(X)$ l'ensemble des ouverts de $X$.
         \item Pour tout $x\in\Ouv(X)$, 
                  $$\stackrel{\cup}{s}(x)=X
                  \quad\text{et}\quad
                  \stackrel{\cup}{t}(x)=X$$
         \item La composition est simplement la réunion.
         \item l'alexis est l'ensemble vide.
         \end{enumerate}
\item $\left((\Ouv(X);\stackrel{\cap}{s};\stackrel{\cap}{t};
         \cap);X\right)$ avec
         \begin{enumerate}
         \item $\Ouv(X)$ l'ensemble des ouverts de $X$.
         \item Pour tout $x\in\Ouv(X)$,  
                  $$\stackrel{\cap}{s}(x)=x
                  \quad\text{et}\quad 
                  \stackrel{\cap}{t}(x)=x$$
         \item La composition est simplement l'intersection.
         \item L'alexis est l'ensemble $X$.
         \end{enumerate}
\end{enumerate}
La vérification des axiomes (GR 1), (GR 2), (GR 3), (GALEX 1) et (GALEX 2) est 
tellement simple que nous nous permettons de ne pas l'écrire.
\end{exem}
% remarque
\begin{rema}\label{rema:alex:1.3}
On peut munir l'ensemble des ouverts d'autres structures de groupement d'Alexandroff, 
mais celles-ci s'avéreront, plus tard, un peu meilleures.\\
Nous aurions aussi pu faire le même genre de constructions avec l'ensemble des 
fermés.
\end{rema}
% exemple
\begin{exem}\label{exem:alex:1.2}
Nous avons vu, dans le chapitre \ref{chap:gr}, exemple \ref{exem:gr:2.2}, qu'à partir 
d'un monoïde non vide $(M;\bullet)$ nous pouvions définir des groupements. 
Le problème qui se présentait alors était que nous étions obligé
de choisir un élément quelconque dans $M$. Il n'y avait pas de construction canonique.
En fait il est tout à fait simple d'associer canoniquement un groupement d'Alexandroff
à un monoïde. En effet il suffit de réunir l'ensemble $M$ avec les singletons 
$\{M\}$ et $\{\{M\}\}$ pour obtenir l'ensemble $\hat{M}$ et de poser
\begin{itemize}
\item pour tout $x\in\hat{M}$, si $x\neq M$, alors $$\hat{s}(x)=\{M\}
          \quad\text{et}\quad\hat{t}(x)=\{M\}$$
          sinon $$\hat{s}(M)=\hat{t}(M)=M$$
\item pour tous éléments $x$ et $y$ de $\hat{M}$,
          $$x\hat{\bullet}y=
          \begin{cases}
          x\bullet y &\text{si }x\neq\{M\},M \text{ et }
          y\neq\{M\},M\\
          y &\text{si }x=M\\
          x &\text{si }y=M\\
          \{M\} & \text{sinon}
          \end{cases}$$
\end{itemize}
$((\hat{M};\hat{s};\hat{t};\hat{\bullet});M)$ est un groupement d'Alexandroff :
\begin{description}
\item[(GR 1)] La vérification est très simple. Il suffit de prendre $x\in\hat{M}$,
         et de distinguer dans les calculs le cas où $x\neq M$ et le cas où $x=M$.
\item[(GR 2)] Soient $x$ et $y$ deux éléments de $\hat{M}$ tels que 
         $\hat{s}(x)=\hat{t}(y)$. Il y a deux possibilités :
         \begin{itemize}
         \item $\hat{s}(x)=\hat{t}(y)=\{M\}$ et alors $x$ et $y$ sont différents de $M$.\\                  Dans ce cas $x\hat{\bullet}y$ est égal à $x\bullet y\neq M$ ou à $\{M\}$.
                  Ainsi $$\hat{s}(x\hat{\bullet}y)=\{M\}=\hat{s}(y)$$et
                   $$\hat{t}(x\hat{\bullet}y)=\{M\}=\hat{t}(x)$$
         \item $\hat{s}(x)=\hat{t}(y)=M$ et alors $x$ et $y$ sont égaux à $M$. D'où
                   $$\hat{s}(x\hat{\bullet}y)=\hat{s}(M)=M=\hat{s}(y)$$et 
                   $$\hat{t}(x\hat{\bullet}y)=\hat{t}(M)=M=\hat{t}(x)$$
         \end{itemize}
\item[(GR 3)] Soient $x$, $y$ et $z$ trois éléments de $\hat{M}$ tels que
         $\hat{s}(x)=\hat{t}(y)$ et $\hat{s}(y)=\hat{t}(z)$. Comme précédemment,
         deux cas se présentent
         \begin{itemize}
         \item $\hat{s}(x)=\hat{t}(y)=\hat{s}(y)=\hat{t}(z)=\{M\}$ 
                  et alors $x$, $y$ et $z$ sont différents de $M$.
                  Dans ce cas $(x\hat{\bullet}y)\hat{\bullet}z$ et
                  $x\hat{\bullet}(y\hat{\bullet}z)$ sont égaux, respectivement, à 
                  $(x\bullet y)\bullet z$ et 
                  $x\bullet (y\bullet z)$ si $x$, $y$ et $z$ sont tous différents de $\{M\}$.
                  Comme $\bullet$ est associative, on obtient l'égalité dans ce cas.\\
                  Si l'un des trois éléments est égal à $\{M\}$, alors les deux compositions 
                  sont égales à $\{M\}$.
         \item $\hat{s}(x)=\hat{t}(y)=\hat{s}(y)=\hat{t}(z)=M$ 
                   et alors $x$, $y$ et $z$ sont égaux à $M$. D'où
                   $$(x\hat{\bullet}y)\hat{\bullet}z=M= 
                   x\hat{\bullet}(y\hat{\bullet}z)$$
         \end{itemize}
\item[(GALEX 1)] Immédiat.
\item[(GALEX 2)] Ce déduit directement de la définition.
\end{description}
\end{exem}
% remarque
\begin{rema}\label{rema:alex:1.4}
Dans cet exemple, on a utilisé le fait bien connu qu'un ensemble ne peut pas 
appartenir à lui-même.
\end{rema}
%%%%%%%%%%%
%%% section 2 %%%
%%%%%%%%%%%
\section{Groupement d'Alexandroff associé à un groupement}\label{sec:alex:2}
Nous allons nous inspirer de la construction faite pour 
les monoïdes.\\
Soit $\left(\BBB;s;t;\OP\right)$ un groupement quelconque. Posons 
\begin{itemize}
\item $\tilde{\BBB}=\BBB\cup\{\BBB\}$ ;
\item $\tilde{s}:\tilde{\BBB}\to\tilde{\BBB}$, 
         $$\tilde{s}(x)=
         \begin{cases}
                  s(x) & \text{si }x\neq\BBB\\
                  \BBB & \text{si }x=\BBB
         \end{cases}$$
\item $\tilde{t}:\tilde{\BBB}\to\tilde{\BBB}$, 
         $$\tilde{t}(x)=
         \begin{cases}
                  t(x) & \text{si }x\neq\BBB\\
                  \BBB & \text{si }x=\BBB
         \end{cases}$$
\item $\tilde{\bullet}:\tilde{\BBB}\times\tilde{\BBB}\to\tilde{\BBB}$,
         $$x\tilde{\bullet}y=
         \begin{cases}
                  x\bullet y & \text{si }x\neq\BBB\text{ et }y\neq\BBB\\
                  x & \text{si }y=\BBB\\
                  y & \text{si }x=\BBB
         \end{cases}$$
\end{itemize}
Vérifions que $\left((\tilde{\BBB};\tilde{s};\tilde{t};\tilde{\bullet});\BBB\right)$ 
est un groupement d'Alexandroff.

Commençons par vérifier que 
$\left(\tilde{\BBB};\tilde{s};\tilde{t};\tilde{\bullet}\right)$ est un  groupement.
\begin{description}
\item[(GR 1)] Soient $x\in\tilde{\BBB}$. Deux cas se présentent. Si $x\neq\BBB$, 
         on a
         \begin{alignat*}{3}
         \tilde{s}\tilde{s}(x)&=\tilde{s}(s(x))&=s(s(x))&=s(x)&=\tilde{s}(x)\\
         \tilde{s}\tilde{t}(x)&=\tilde{s}(t(x))&=s(t(x))&=t(x)&=\tilde{t}(x)\\
         \tilde{t}\tilde{t}(x)&=\tilde{t}(t(x))&=t(t(x))&=t(x)&=\tilde{t}(x)\\
         \tilde{t}\tilde{s}(x)&=\tilde{t}(s(x))&=t(s(x))&=s(x)&=\tilde{s}(x)
         \end{alignat*}
         car $s(x)\neq\BBB$ et $t(x)\neq\BBB$. Si $x=\BBB$, alors
         \begin{alignat*}{2}
         \tilde{s}\tilde{s}(\BBB)&=\tilde{s}(\BBB)&=\BBB&=\tilde{s}(\BBB)\\
         \tilde{s}\tilde{t}(\BBB)&=\tilde{s}(\BBB)&=\BBB&=\tilde{t}(\BBB)\\
         \tilde{t}\tilde{t}(\BBB)&=\tilde{t}(\BBB)&=\BBB&=\tilde{t}(\BBB)\\
         \tilde{t}\tilde{s}(\BBB)&=\tilde{t}(\BBB)&=\BBB&=\tilde{s}(\BBB)
         \end{alignat*}
\item[(GR 2)] Soient $x$ et $y$ deux éléments de $\tilde{\BBB}$ tels que 
         $\tilde{s}(x)=\tilde{t}(y)$. D'après la construction, $\tilde{s}(x)=\BBB$ et 
         $\tilde{t}(y)=\BBB$ si, et seulement si $x$ et $y$ sont égaux à $\BBB$. Par 
         conséquent, nous n'avons que deux cas à étudier :
         \begin{itemize}
         \item soit $x$ et $y$ appartiennent tous les deux à $\BBB$ et alors
                  $$s(x)=\tilde{s}(x)=\tilde{t}(y)=t(y)$$
                  et
                  $$\tilde{s}(x\tilde{\OP}y)=s(x\OP y)=s(y)=\tilde{s}(y)$$
                  $$\tilde{t}(x\tilde{\OP}y)=t(x\OP y)=t(x)=\tilde{t}(x)$$
         \item soit $x=y=\BBB$ et alors
                  $$\tilde{s}(x\tilde{\OP}y)=\tilde{s}(\BBB)=\tilde{s}(y)$$
                  $$\tilde{s}(x\tilde{\OP}y)=\tilde{t}(\BBB)=\tilde{t}(x)$$
         \end{itemize}
\item[(GR 3)] Soient $x$, $y$ et $z$ trois éléments de $\tilde{\BBB}$ tels que 
         $\tilde{s}(x)=\tilde{t}(y)$ et $\tilde{s}(y)=\tilde{t}(z)$. Comme 
         précédemment, on remarque que l'on a deux cas : soit $x$, $y$ et $z$ 
         appartiennent tous les trois à $\BBB$, soit ils sont tous les trois égaux à
         $\BBB$. Dans le premier cas, on a
         $$(x\tilde{\OP}y)\tilde{\OP}z=(x\OP y)\OP z=x\OP(y\OP z)
         =x\tilde{\OP}(y\tilde{\OP}z)$$
         car $(\BBB;s;t;\OP)$ est un groupement. Dans le second cas, on a
         $$(x\tilde{\OP}y)\tilde{\OP}z=\BBB\quad\text{et}\quad
         \BBB=x\tilde{\OP}(y\tilde{\OP}z)$$
\end{description}
Il ne nous reste plus que les deux conditions (GALEX 1)  et (GALEX 2) qui en font 
un groupement d'Alexandroff. Mais celles-ci sont absolument évidentes par construction.
%%%%%%%%%%%
%%% section 3 %%%
%%%%%%%%%%%
\section{Les g-morphismes d'Alexandroff}\label{sec:alex:3}
Bien entendu, les groupements d'Alexandroff étant avant tout des groupements, les 
g-mor\-phi\-smes habituels sont toujours à envisager. Mais il paraît intéressant, vu 
la présence des alexis, de s'intéresser plus particulièrement au g-mor\-phi\-smes qui 
transforment un alexis en un alexis.

Ainsi, un g-morphisme $f$ du groupement d'Alexandroff $(\BBB_{1},\alpha_{1})$ vers
le groupement d'Alexandroff $(\BBB_{2};\alpha_{2})$ est dit \emph{d'Alexandroff}
\index{g-morphisme!g-morphisme d'Alexandroff} s'il vérifie
\begin{description}
\item[(MALEX)] $f(\alpha_{1})=\alpha_{2}$.
\end{description}
En suivant la démonstration du théorème \ref{theo:gr:3.2}, on obtient
\begin{theo}\label{theo:alex:3.1}
Si nous notons $\GMORALEX$ l'ensemble des petits g-morphismes d'Alexandroff, 
alors $$\left(\GMORALEX;s;t;\OP\right)$$ est une sous-catégorie de la catégorie
$$\left(\GMOR;s;t;\OP\right)$$
\end{theo}
Reprenons les exemples de la section \ref{sec:alex:1}.
% exemple
\begin{exem}\label{exem:alex:3.1}
Soit $f:X_{1}\to X_{2}$ une application continue entre 
espaces topologiques. Considérons l'application $\stackrel{\cup}{f}:
\Ouv(X_{2})\to\Ouv(X_{1})$ définie par
$$\stackrel{\cup}{f}(x)=f^{-1}(x)$$
pour tout $x\in\Ouv(X_{2})$. Cette application est bien définie car 
l'image réciproque d'un ouvert par une application continue est un ouvert. Vérifions 
que $\stackrel{\cup}{f}$ est un g-morphisme d'Alexandroff.
\begin{description}
\item[(GMOR)] Soient $x$ et $y$ deux éléments de $\Ouv(X_{2})$ 
         tels que $\stackrel{\cup}{s_{2}}(x)=\stackrel{\cup}{t_{2}}(y)$. 
         Or par définition, cela implique que $x$ et $y$ sont quelconques. Pour la même 
         raison, la condition 
         $$\stackrel{\cup}{s_{1}}(\stackrel{\cup}{f}(x))=
         \stackrel{\cup}{t_{1}}(\stackrel{\cup}{f}(y))$$
         est trivialement vérifiée car les deux sont égales à $X_{1}$. De plus, on a
         $$\stackrel{\cup}{f}(x\cup y)=f^{-1}(x\cup y)
         =f^{-1}(x)\cup f^{-1}(y)=\stackrel{\cup}{f}(x)\cup
         \stackrel{\cup}{f}(y)$$
\item[(MALEX)] Ce n'est rien d'autre que la conséquence du calcul
         $$\stackrel{\cup}{f}(X_{2})=f^{-1}(X_{2})=X_{1}$$
\end{description}
Si $f_{1}:X_{1}\to X_{2}$ et $f_{2}:X_{2}\to X_{3}$ sont deux applications 
continues, alors, pour tout $x\in\Ouv(X_{3})$, on a
$$(f_{2}f_{1})^{\cup}(x)=(f_{2}f_{1})^{-1}(x)=f_{1}^{-1}f_{2}^{-1}(x)
=\stackrel{\cup}{f_{1}}\OP\stackrel{\cup}{f_{2}}(x)$$
De plus, si $X$ est un espace topologique, alors 
$$(\Id_{X})^{\cup}(x)=x=\Id_{\Ouv(X)}$$
En résumé, nous venons de montrer que $^{\cup}$ est un foncteur contravariant de 
la catégorie des petits espaces topologiques dans la catégorie $\GMORALEX$.\\
De la même manière, $^{\cap}$ est lui aussi un foncteur contravariant de
la catégorie des petits espaces topologiques dans la catégorie des petits 
groupements d'Alexandroff $\GMORALEX$.
\end{exem}
% exemple
\begin{exem}\label{exem:alex:3.2}
Considérons $f:(M_{1};\bullet_{1})\to(M_{2};\bullet_{2})$ un morphisme 
de monoïdes. On définit une application $\hat{f}$ de $\hat{M_{1}}$ dans 
$\hat{M_{2}}$ en posant
$$\hat{f}(x)=
\begin{cases}
f(x) &\text{si }x\neq\{M\},M\\
\{M\}&\text{si }x=\{M\}\\
M&\text{si }x=M
\end{cases}$$
$\hat{f}$ est un g-morphisme d'Alexandroff car
\begin{description}
\item[(GMOR)] Soient $x$ et $y$ deux éléments de $\hat{M_{1}}$ vérifiant
         $\hat{s}(x)=\hat{t}(y)$. On a deux cas :
         \begin{itemize}
         \item $\hat{s}(x)=\hat{t}(y)=\{M_{1}\}$ et donc $x$ et $y$ sont différents 
                  de $M_{1}$.\\
                  On a donc
                  $$\hat{s}(\hat{f}(x))=\hat{s}(f(x))=\{M_{2}\}=\hat{t}(f(y))
                  =\hat{t}(\hat{f}(y))$$
                  Si les deux sont différents de $\{M_{1}\}$, alors
                  $$\hat{f}(x\hat{\bullet_{1}}y)=f(x\bullet_{1}y)
                  =f(x)\bullet_{1}f(y)=\hat{f}(x)\hat{\bullet_{1}}\hat{f}(y)$$
                  Si $x$ ou $y$ est égal à $\{M_{1}\}$, alors 
                  $x\hat{\bullet_{1}y}=\{M_{1}\}$ et
                  $\hat{f}(x)=\{M_{2}\}$ ou $\hat{f}(y)=\{M_{2}\}$. D'où
                  $$\hat{f}(x\hat{\bullet_{1}}y)=\hat{f}(\{M_{1}\})=\{M_{2}\}
                  =\hat{f}(x)\hat{\bullet_{1}}\hat{f}(y)$$
         \item $\hat{s}(x)=\hat{t}(y)=M_{1}$ et donc $x$ et $y$ sont égaux à $M_{1}$.
                  Ainsi
                  $$\hat{s}\hat{f}(x)=M_{2}=\hat{t}\hat{f}(y)$$et
                  $$\hat{f}(x\hat{\bullet_{1}}y)=\hat{f}(M_{1})=M_{2}
                  =M_{2}\hat{\bullet_{2}}M_{2}=\hat{f}(x)\hat{\bullet_{2}}
                  \hat{f}(y)$$
         \end{itemize}
\item[(MALEX)] Par construction.
\end{description}
Il est clair que $\hat{\ }$ est un foncteur de la catégorie des petits monoïdes dans la 
catégorie des petits groupements d'Alexandroff $\GMORALEX$.
\end{exem}
En s'inspirant de cet exemple on démontre
% théorème
\begin{theo}\label{theo:alex:3.2}
Il existe un foncteur\ $\tilde{ }$ de la catégorie des petits groupements $\GMOR$ dans 
la catégorie des petits groupements d'Alexandroff $\GMORALEX$.\\
Nous l'appellerons \emph{foncteurs d'Alexandroff}
\index{foncteur!foncteur d'Alexandroff}.\\
Pour tout $(f:\BBB_{1}\to\BBB_{2})\in\GMOR$ et tout $x\in\tilde{\BBB_{1}}
=\BBB_{1}\cup\{\BBB_{1}\}$, l'application $\tilde{f}:\tilde{\BBB_{1}}\to
\tilde{\BBB_{2}}$ est définie par
$$\tilde{f}(x)=
\begin{cases}
f(x) & \text{si }x\neq\BBB_{1}\\
\BBB_{2} & \text{si }y=\BBB_{1}
\end{cases}$$
\end{theo}
\Dem Commençons par vérifier que $\tilde{f}$ est bien un g-morphisme d'Alexandroff.
Il est évident, d'après la construction, que la condition (MALEX) est vérifiée. Considérons 
maintenant deux éléments $x$ et $y$ de $\tilde{\BBB_{1}}$ tels que 
$\tilde{s}(x)=\tilde{t}(y)$. On a deux cas
\begin{itemize}
\item $\tilde{s}(x)=\tilde{t}(y)=\BBB_{1}$ et donc $x=y=\BBB_{1}$. Ainsi
         $$\tilde{s}\tilde{f}(x)=\BBB_{2}=\tilde{t}\tilde{f}(y)$$
         et
         $$\tilde{f}(x\tilde{\OP_{1}}y)=\tilde{f}(\BBB_{1})=\BBB_{2}
         =\tilde{f}(x)\tilde{\OP_{2}}\tilde{f}(y)$$
\item $\tilde{s}(x)=\tilde{t}(y)\neq\BBB_{1}$ et donc $s(x)=t(y)$. Par conséquent, 
         $$\tilde{s}\tilde{f}(x)=sf(x)=tf(y)=\tilde{t}\tilde{f}(y)$$
         et
         $$\tilde{f}(x\tilde{\OP_{1}}y)=f(x\OP_{1}y)=f(x)\OP_{2}f(y)
         =\tilde{f}(x)\tilde{\OP_{2}}\tilde{f}(y)$$
\end{itemize}
La condition (GMOR) est donc satisfaite. Comme $\tilde{}$ est une application bien
définie de $\GMOR$ dans $\GMORALEX$, il nous suffit maintenant de prouver que c'est un 
foncteur de la catégorie $(\GMOR;s;t\OP)$ dans la catégorie $(\GMORALEX;s;t;\OP)$.
\begin{description}
\item[(FONC 1)] Pour tout g-morphisme $f:\BBB_{1}\to\BBB_{2}$, on a
         $$\tilde{}s(f)=\tilde{\Id_{\BBB_{1}}}\quad\text{et}\quad
         s\tilde{}(f)=s(\tilde{f})=\Id_{\tilde{\BBB_{1}}}$$
         $$\tilde{}t(f)=\tilde{\Id_{\BBB_{2}}}\quad\text{et}\quad
         t\tilde{}(f)=t(\tilde{f})=\Id_{\tilde{\BBB_{2}}}$$
         or il est clair au vue de la construction que $\tilde{\Id_{\BBB}}=
         \Id_{\tilde{\BBB}}$ pour tout groupement $\BBB$.
\item[(FONC 2)] Considérons maintenant deux g-morphismes $f:\BBB_{1}\to
         \BBB_{2}$ et $g:\BBB_{2}\to\BBB_{3}$ tels que $s(g)=t(f)=\Id_{\BBB_{2}}$.
         $$\tilde{}(g\OP f)(x)=(g\OP f)(x)=g(f(x))$$et
         $$(\tilde{g}\OP\tilde{f})(x)=\tilde{g}(\tilde{f}(x))=\tilde{g}(f(x))
         =g(f(x))$$
         pour $x\in\BBB_{1}$  car $f(x)\neq\alpha_{2}$.\\
         De plus
         $$\tilde{}(g\OP f)(\alpha_{1})=\alpha_{3}$$et
         $$(\tilde{g}\OP\tilde{f})(\alpha_{1})=\tilde{g}(\tilde{f}(x))
         =\tilde{g}(\alpha_{2})=\alpha_{3}$$Par conséquent
         $$\tilde{}(g\OP f)=\tilde{g}\OP\tilde{f}$$\QED
\end{description}
%%%%%%%%%%%
%%% section 4 %%%
%%%%%%%%%%%
\section{Les g-transformations d'Alexandroff}
\index{g-transformation!g-transformation d'Alexandroff}
\label{sec:alex:4}
En s'inspirant de la définition des g-transformations donnée dans la section \ref{sec:gr:3},
on définit une \emph{g-transformation d'Alexandroff} comme étant un quadruplet
$\left(\eta^{1};\eta^{2};f_{1};f_{2}\right)$ où $$f_{1},f_{2}:
(\BBB_{1};\alpha_{1})\to(\BBB_{2};\alpha_{2})$$ sont des g-morphismes d'Alexandroff
et où $$\eta^{1},\eta^{2}:\BBB_{1}\to\BBB_{2}$$ sont des applications d'ensembles qui 
vérifient les conditions suivantes :
\begin{description}
\item[(GTRALEX 1)] Pour tout $x\in\BBB_{1}$, on a
         \begin{itemize}
         \item soit $\eta^1(x)=\alpha_{2}$, soit\quad
                   $s_{2}\eta^{1}(x)=s_{2}f_{1}(x)\quad\text{et}\quad
                   t_{2}\eta^{1}(x)=s_{2}f_{2}(x)$ ;
          \item soit $\eta^2(x)=\alpha_{2}$, soit\quad
                   $s_{2}\eta^{2}(x)=t_{2}f_{1}(x)\quad\text{et}\quad
                   t_{2}\eta^{2}(x)=t_{2}f_{2}(x)$.
           \end{itemize}
\item[(GTRALEX 2)] Pour tous $x$ et $y$ dans $\BBB_{1}$ satisfaisant 
           $s_{1}(x)=t_{1}(y)$, on a
           $$f_{2}(x)\OP_{2}\eta^{1}(x)=\eta^{2}(x)\OP_{2}f_{1}(x)$$
           et
           $$\eta^{1}(x)=\eta^{2}(y)$$
\end{description}
Comme nous l'avons déjà fait dans la section \ref{sec:gr:4}, nous noterons
$$(\eta^{1};\eta^{2}):f_{1}\leadsto f_{2}$$
et même parfois $$\eta:f_{1}\leadsto f_{2}$$
% proposition
\begin{prop}\label{prop:alex:4.1}
Si $(\eta^{1};\eta^{2}):f_{1}\to f_{2}$ est une g-transformation 
d'Alexandroff au dessus de $\BBB_{1}$ et $\BBB_{2}$, alors, pour tout $x\in\BBB_{1}$,
\begin{enumerate}
\item  $f_{2}(x)\OP_{2}\eta^{1}(x)=\eta^{2}(x)\OP_{2}f_{1}(x)$
\item $\eta^{1}(x)=\eta^{2}(s_{1}(x))$ et $\eta^{1}(t_{1}(x))=\eta^{2}(x)$
\end{enumerate}
\end{prop}
\Dem\ Ce sont des conséquences directes du fait que le couple $x$, $s_{1}(x)$ et le 
couple $t_{1}(x)$, $x$ vérifie la condition (GTRALEX).\QEDb
Le premier et le plus simple des exemples que l'on puisse donner est celui
construit dans la proposition ci-dessous.
% proposition
\begin{prop}\label{prop:alex:4.2}
Soit $f:(\BBB_{1};\alpha_{1})\to(\BBB_{2};\alpha_{2})$ un g-morphisme
d'Alexandroff. Par abus, notons $\alpha_{2}:\BBB_{1}\to\BBB_{2}$ l'application 
constante qui a tout $x\in\BBB_{1}$ associe l'alexis $\alpha_{2}$ de $\BBB_{2}$.

Le couple d'applications $(\alpha_{2};\alpha_{2})$ définit une 
g-transformation d'Alexandroff de $f$ vers lui-même.
$$(\alpha_{2};\alpha_{2}):f\leadsto f$$
\end{prop}
\Dem\ Par définition, la condition (GRALEX 1) est trivialement vérifiée. De plus, 
pour tout $x\in\BBB_{1}$, on a
$$f(x)\OP_{2}\alpha_{2}(x)=f(x)\OP_{2}\alpha_{2}=f(x)$$
et
$$\alpha_{2}(x)\OP_{2}f(x)=\alpha_{2}\OP_{2}f(x)=f(x)$$
Par conséquent pour tous $x$ et $y$ de $\BBB_{1}$ tel que $s_{1}(x)=t_{1}(y)$, on a
$$f(x)\OP_{2}\alpha_{2}(x)=\alpha_{2}(x)\OP_{2}f(x)$$
et
$$\alpha_{2}(x)=\alpha_{2}(y)=\alpha_{2}$$
D'où la condition (GRALEX 2).\QEDb
Comme nous l'avons fait dans le chapitre \ref{chap:gr}, notons
$$\GTRANSALEX$$
l'ensemble des \emph{petites g-transformations d'Alexandroff}
\index{g-transformation!petite g-transformation d'Alexandroff} 
(les groupements d'Alexandroff de base sont petits). Avant de munir cet ensemble
de structures de groupement, commençons par les lemmes suivants :
% lemme
\begin{lemm}\label{lemm:alex:4.1}
Soient $\eta:f_{1}\leadsto f_{2}$ et $\eta':f_{2}\leadsto f_{3}$ deux 
g-transformations d'alexandroff au dessus des groupements d'Alexandroff
$(\BBB_{1};\alpha_{1})$ et $(\BBB_{2};\alpha_{2})$.

 Si nous écrivons
$$\eta'\OP\eta=(\eta'^1\OP_{2}\eta^1;\eta'^2\OP_{2}\eta^2)$$
où pour tout $x\in\BBB_{1}$ et $i\in\{1;2\}$, on a
$$(\eta'^{i}\OP_{2}\eta^{i})(x)=\eta'^{2}(x)\OP_{2}\eta^{2}(x)$$
alors $\eta'\OP\eta$ est une g-transformation d'Alexandroff de $f_{1}$ vers $f_{3}$.
\end{lemm}
\Dem\ Bien qu'elle soit simple, la démonstration est particulièrement pénible.
\begin{description}
\item[(GTRALEX 1)] Soit $x$ un élément de $\BBB_{1}$. Nous allons nous contenter
         de prouver la condition sur $\eta'^{1}\OP_{2}\eta^{1}$ car la justification
         pour $\eta'^{2}\OP_{2}\eta^{2}$ est identique. Plusieurs cas sont à prendre
         en considération :
         \begin{itemize}
         \item $\eta'^{1}(x)=\alpha_{2}$ et $\eta^{1}(x)=\alpha_{2}$. Alors
                   $(\eta'^{1}\OP_{2}\eta^{1})(x)=\alpha_{2}$.
         \item $\eta'^{1}(x)=\alpha_{2}$ et $\eta^{1}(x)\neq\alpha_{2}$. Alors
                   $(\eta'^{1}\OP_{2}\eta^{1})(x)=\eta^{1}(x)$ et, d'après (GTRALEX 1),
                   $$s_{2}(\eta'^{1}\OP_{2}\eta^{1})(x)=s_{2}\eta^{1}(x)=s_{2}f_{1}(x)$$
                   $$t_{2}(\eta'^{1}\OP_{2}\eta^{1})(x)=t_{2}\eta^{1}(x)=s_{2}f_{2}(x)$$
                   La proposition \ref{prop:alex:4.1} appliquée à $\eta'$ et 
                   $\eta'^{1}(x)=\alpha_{2}$, nous donne
                   $$f_{3}(x)=\eta'^{2}(x)\OP_{2}f_{2}(x)$$
                   Soit $\eta'^{2}(x)=\alpha_{2}$ et alors $$f_{3}(x)=f_{2}(x)$$
                   $$s_{2}f_{3}(x)=s_{2}f_{2}(x)$$
                   soit $s_{2}\eta'^{2}(x)=t_{2}f_{2}(x)$ et alors, d'après (GR 2), 
                   $$s_{2}f_{3}(x)=s_{2}f_{2}(x)$$
                   Ainsi on a toujours 
                   $$s_{2}(\eta'^{1}\OP_{2}\eta^{1})(x)=s_{2}f_{1}(x)$$
                   $$t_{2}(\eta'^{1}\OP_{2}\eta^{1})(x)=s_{2}f_{3}(x)$$
         \item $\eta'^{1}(x)\neq\alpha_{2}$ et $\eta^{1}(x)=\alpha_{2}$. Ainsi, 
                   en utilisant (ALEX 2), on trouve
                   $$(\eta'^{1}\OP_{2}\eta^{1})(x)=\eta'^{1}(x)$$
                   D'aprés (GTRALEX 1)
                   $$s_{2}(\eta'^{1}\OP_{2}\eta^{1})(x)=s_{2}\eta'^{1}(x)=s_{2}f_{2}(x)$$
                   $$t_{2}(\eta'^{1}\OP_{2}\eta^{1})(x)=t_{2}\eta'^{1}(x)=s_{2}f_{3}(x)$$ 
                   $\eta$ est une g-transformation d'Alexandroff et $\eta^{1}(x)=\alpha_{2}$.
                   D'où
                   $$f_{2}(x)=\eta^{2}(x)\OP_{2}f_{1}(x)$$
                   Comme précédemment, soit $\eta^{2}=\alpha_{2}$, soit 
                   $s_{2}\eta^{2}(x)=t_{2}f_{1}(x)$. Dans les deux cas, on obtient
                   $$s_{2}f_{2}(x)=s_{2}f_{1}(x)$$
                   Finalement,
                   $$s_{2}(\eta'^{1}\OP_{2}\eta^{1})(x)=s_{2}f_{1}(x)$$
                   $$t_{2}(\eta'^{1}\OP_{2}\eta^{1})(x)=s_{2}f_{3}(x)$$
          \item  $\eta'^{1}(x)=\alpha_{2}$ et $\eta^{1}(x)=\alpha_{2}$. C'est plus 
                   simple car on a
                   $$(\eta'^{1}\OP_{2}\eta^{1})(x)=\alpha_{2}$$
          \end{itemize}
\item[(GTRALEX 2)] Pour tous $x$ et $y$ dans $\BBB_{1}$ vérifiant $s_{2}(x)=t_{2}(y)$,
          on a
          \begin{align*}
          f_{3}(x)\OP_{2}(\eta'^{1}\OP_{2}\eta^{1})(x)
          &=f_{3}(x)\OP_{2}(\eta'^{1}(x)\OP_{2}\eta^{1}(x))\\
          &=(f_{3}(x)\OP_{2}\eta'^{1}(x))\OP_{2}\eta^{1}(x)\\
          &=(\eta'^{2}(x)\OP_{2}f_{2}(x))\OP_{2}\eta^{1}(x)\\
          &=\eta'^{2}(x)\OP_{2}(f_{2}(x)\OP_{2}\eta^{1}(x))\\
          &=\eta'^{2}(x)\OP_{2}(\eta^{2}(x)\OP_{2}f_{1}(x))\\
          &=(\eta'^{2}(x)\OP_{2}\eta^{2}(x))\OP_{2}f_{1}(x)
          =(\eta'^{2}\OP_{2}\eta^{2})(x)\OP_{2}f_{1}(x)
          \end{align*}
           et
           $$(\eta'^{1}\OP_{2}\eta^{1})(x)=\eta'^{1}(x)\OP_{2}\eta^{1}(x)
           =\eta'^{2}(y)\OP_{2}\eta^{2}(y)=(\eta'^{2}\OP_{2}\eta^{2})(y)$$
           \QED
\end{description}
% lemme
\begin{lemm}\label{lemm:alex:4.2}
Soient $f:\BBB'_1\to\BBB_1$, $g:\BBB_2\to\BBB'_2$ deux g-morphismes d'Alexandroff et 
$(\eta^1;\eta^2):(f_1:\BBB_1\to\BBB_2)\leadsto(f_2:\BBB_1\to\BBB_2)$ une 
g-transformation.

Les quadruplets $\left(\eta^1f;\eta^2f;f_1f;f_2f\right)$ et 
$\left(g\eta^1;g\eta^2;gf_1;gf_2\right)$
sont des g-transformations d'Alexandroff.
\end{lemm}
\Dem\ Nous ne ferons la démonstration que pour $(\eta^{1}f;\eta^{2}f):f_{1}f
\leadsto f_{2}f$. Nous avons déjà vu que $f_{1}f$ et $f_{2}f$ sont des g-morphismes 
d'Alexandroff. 
\begin{description}
\item[(GTRALEX 1)] Pour $x\in\BBB'_{1}$, on a 
         soit $$(\eta^{1}f)(x)=\eta^{1}(f(x))=\alpha_{2}$$soit
         $$s_{2}(\eta^{1}f)(x)=s_{2}\eta^{1}(f(x))=s_{2}f_{1}(f(x))
         =s_{2}(f_{1}f)(x)$$et
         $$t_{2}(\eta^{1}f)(x)=t_{2}\eta^{1}(f(x))=s_{2}f_{2}(f(x))
         =s_{2}(f_{2}f)(x)$$
         On prouve de la même façon les égalités concernant $\eta^{2}f$.
\item[(GTRALEX 2)] Soient $x$ et $y$ deux éléments de $\BBB'_{1}$ tels que 
         $s'_{1}(x)=t'_{1}(y)$. Comme $f$ est un g-morphisme, on a
         $$s_{1}f(x)=t_{1}f(y)$$
         Or $(\eta^{1};\eta^{2})$ est une g-transformation d'Alexandroff, on a donc
         $$f_{2}(f(x))\OP_{2}\eta^{1}(f(x))=\eta^{2}(f(x))\OP_{2}f_{1}(f(x))$$
         $$(f_{2}f)(x)\OP_{2}(\eta^{1}f)(x)=(\eta^{2}f)(x)\OP_{2}(f_{1}f)(x)$$
         et
         $$(\eta^{1}f)(x)=\eta^{1}(f(x))=\eta^{2}(f(y))=(\eta^{2}f)(y)$$
         \QED         
\end{description}
Pour simplifier, si au lieu de noter $(\eta^{1};\eta^{2})$ nous notons $\eta$ la
g-transformation d'Alexandroff du lemme, alors nous écrirons
$\eta f$ pour $(\eta^{1}f;\eta^{2}f)$ et $g\eta$ pour $(g\eta^{1};g\eta^{2})$.
% lemme
\begin{lemm}\label{lemm:alex:4.3}
Soient $\eta:f_{1}\leadsto f_{2}$ et $\eta':f'_{1}\leadsto f'_{2}$ deux 
g-transformations d'Alexandroff. Supposons que $f_{1},f_{2}:(\BBB_{1};\alpha_{1})
\to(\BBB_{2};\alpha_{2})$ et $f'_{1},f'_{2}:(\BBB_{2};\alpha_{2})\to
(\BBB_{3};\alpha_{3})$. Les g-transformations d'Alexandroff 
$$(\eta' f_{2})\OP(f'_{1}\eta)\quad\text{et}\quad
(f'_{2}\eta)\OP(\eta'f_{1})$$de $f'_{1}f_{2}$
vers $f'_{2}f_{2}$ existent.
\end{lemm}
\Dem\ D'après les deux lemmes précédents, il suffit de constater que 
$$\eta'f_{2}:f'_{1}f_{2}\leadsto f'_{2}f_{2}
\quad;\quad f'_{1}\eta:f'^{1}f_{1}\leadsto f'^{1}f_{2}$$et 
$$f'_{2}\eta:f'_{2}f_{1}\leadsto f'_{2}f_{2}
\quad;\quad\eta'f_{1}:f'^{1}f_{1}\leadsto f'^{2}f_{1}$$\QEDb
Nous pouvons maintenir établir le théorème principal de ce chapitre
\begin{theo}\label{theo:alex:4.1}
Notons $\GTRANSALEX$ l'ensemble des petites g-transformations d'Alexandroff. Ce sont 
des g-transformations d'Alexandroff au dessus de petits groupements d'Alexandroff.

D'après les lemmes ci-dessus, nous pouvons définir quatre applications
$$\sigma_{0},\tau_{0},\sigma_{1},\tau_{1}:\GTRANSALEX\to\GTRANSALEX$$
de la manière suivante : Pour tout $\eta:f_{1}\leadsto f_{2}$ avec
$f_{1},f_{2}:(\BBB_{1};\alpha_{1})\to(\BBB_{2};\alpha_{2})$
\begin{align*}
\sigma_{0}(\eta)&=(\alpha_{1};\alpha_{1}):\Id_{\BBB_{1}}\leadsto
\Id_{\BBB_{1}}\\
\tau_{0}(\eta)&=(\alpha_{2};\alpha_{2}):\Id_{\BBB_{2}}\leadsto
\Id_{\BBB_{2}}\\
\sigma_{1}(\eta)&=(\alpha_{2};\alpha_{2}):f_{1}\leadsto f_{1}\\   
\tau_{1}(\eta)&=(\alpha_{2};\alpha_{2}):f_{2}\leadsto f_{2}
\end{align*}
Et toujours à l'aide des mêmes lemmes, nous pouvons construire trois nouvelles
applications
$$\OP,\boxtimes,\boxdot:\GTRANSALEX\times\GTRANSALEX\to\GTRANSALEX$$
en posant pour tous $\eta$ et $\eta'$ de $\GTRANSALEX$,
\begin{align*}
\eta'\OP\eta&=
        \begin{cases}
         \eta'\OP\eta &\text{si }\sigma_{1}(\eta')=\tau_{1}(\eta)\text{ (lemme 
         \ref{lemm:alex:4.1})}\\
         \eta &\text{sinon}
         \end{cases} \\
\eta'\boxtimes\eta&=
        \begin{cases}
         (f'_{2}\eta)\OP(\eta'f_{1}) &\text{si }\sigma_{0}(\eta')=\tau_{0}(\eta)
         \text{ (lemme \ref{lemm:alex:4.3})}\\
         \eta &\text{sinon}
         \end{cases}\\
\eta'\boxdot\eta&=
        \begin{cases}
         (\eta'f_{2})\OP(f'_{1}\eta) &\text{si }\sigma_{0}(\eta')=\tau_{0}(\eta)
         \text{ (lemme \ref{lemm:alex:4.3})}\\
         \eta &\text{sinon}
         \end{cases}
\end{align*}
Les quadruplets $\left(\GTRANSALEX;\sigma_{1};\tau_{1};\OP\right)$, 
$\left(\GTRANSALEX;\sigma_{1};\tau_{1};\OP\right)$, 
$\left(\GTRANSALEX;\sigma_{1};\tau_{1};\OP\right)$
 sont des groupements. 
\end{theo}
\Dem\ Les parties techniques les plus difficiles ont été prouvées dans les différents 
lemmes ci-dessus. Les vérifications des axiomes (GR 1), (GR 2) et (GR 3) sont 
simples mais longues à écrire. Pour cette raison, nous nous permettons
de les laisser au lecteur.\QEDb
En fait, une autre raison nous incite à ne pas écrire la démonstration de ce théorème. Il
n'est pas tout à fait satisfaisant car il ne semble pas possible de démontrer que l'on ait
toujours, quand les calculs sont possibles, les égalités
$$(\eta_{4}\OP\eta_{3})\boxtimes(\eta_{2}\OP\eta_{1})
=(\eta_{4}\boxtimes\eta_{2})\OP(\eta_{3}\boxtimes\eta_{1})$$
$$(\eta_{4}\OP\eta_{3})\boxdot(\eta_{2}\OP\eta_{1})
=(\eta_{4}\boxdot\eta_{2})\OP(\eta_{3}\boxdot\eta_{1})$$
Et puis il ne semble pas non plus que l'on ait
$$\eta_{1}\boxtimes\eta_{2}=\eta_{1}\boxdot\eta_{2}$$
Commençons par étudier une condition pour laquelle cette égalité soit vérifiée.
% proposition
\begin{prop}\label{prop:alex:4.3}

\end{prop}
\Dem\ $(\eta_{1}\boxtimes\eta_{2})^{1}
=(f'_{2}\eta_{2}^{1})\OP_{2}(\eta_{1}^{1}f_{1})
=$
La raison est que nous nous sommes trop laissé influencer par ce qui se passe pour les 
transformations naturelles de catégories. Comme pour les surfaces de Moore, si nous 
devions faire un schéma représentatif d'une g-transformation d'Alexandroff 
$(\eta^{1};\eta^{2}):f_{1}\leadsto f_{2}$, nous 
serions sans doute tous amené à dessiner un carré du genre
$$\xymatrix{
(\BBB_{1};\alpha_{1})\ar[rr]^{f_{1}}\ar[dd]_{\eta^{1}} & &
(\BBB_{2};\alpha_{2})\ar[dd]^{\eta^{2}}\\
\\
(\BBB_{1};\alpha_{1})\ar[rr]_{f_{2}} & & (\BBB_{2};\alpha_{2})
}$$
Ainsi, si nous avons $(\eta^{2};\eta^{3}):f_{3}\leadsto f_{4}$ une g-transformation 
d'Alexandroff, représentée par
$$\xymatrix{
(\BBB_{2};\alpha_{2})\ar[rr]^{f_{3}}\ar[dd]_{\eta^{2}} & &
(\BBB_{3};\alpha_{3})\ar[dd]^{\eta^{3}}\\
\\
(\BBB_{2};\alpha_{2})\ar[rr]_{f_{4}} & & (\BBB_{3};\alpha_{3})
}$$
alors nous avons, comme pour les surfaces de Moore, envie d'écrire que 
$(\eta^{2};\eta^{3})\square(\eta^{1};\eta^{2})$ est une g-transformation 
d'Alexandroff représentée par la justaposition des deux carrés ci-dessus.
$$\xymatrix{
(\BBB_{1};\alpha_{1})\ar[rr]^{f_{3}f_{1}}\ar[dd]_{\eta^{1}} & &
(\BBB_{3};\alpha_{3})\ar[dd]^{\eta^{3}}\\
\\
(\BBB_{1};\alpha_{1})\ar[rr]_{f_{4}f_{2}} & & (\BBB_{3};\alpha_{3})
}$$
Commençons par prouver que c'est bien une g-transformation d'Alexandroff.
% lemme
\begin{lemm}\label{lemm:alex:4.4}
Soient $(\eta^1;\eta^{2}):f_{1}\leadsto f_{2}$ et 
$(\eta^{2};\eta^{3}):f_{3}\leadsto f_{4}$ deux g-transformations 
d'Alexandroff. Alors le couple $(\eta^{1};\eta^{3})$ détermine lui aussi 
une g-transformation d'Alexandroff du g-morphisme d'Al\-ex\-an\-droff 
$f_{3}\OP f_{1}$ vers le g-morphisme d'Alexandroff $f_{4}\OP f_{2}$.
\end{lemm}
\Dem\ Remarquons que par définition, les g-morphismes d'Alexandroff
$f_{3}\OP f_{1}$ et  $f_{4}\OP f_{2}$ ne sont à la base que les applications
d'ensembles $f_{3}f_{1}$ et $f_{4}f_{2}$.

Supposons de plus que l'on a $f_{1},f_{2}:(\BBB_{1};\alpha_{1})\to
(\BBB_{2};\alpha_{2})$ et $f_{3},f_{4}:(\BBB_{2};\alpha_{2})\to
(\BBB_{3};\alpha_{3})$.

Pour tout $y$ appartenant à $\BBB_{2}$, on  sait que
$$\text{soit }\eta^{2}(y)=s_{3}\eta^{2}(y)$$
%%%%%%%%%%%%%%%%%%%%%%%%%%%%%%%%%%%%%%%%%%%%%%%%%%
\chapter{Les $2$-groupements stricts}\label{chap:2-gr}
Dans les deux chapitres précédents, nous avons vu qu'il apparaît très vite et de manière assez
naturelle des ensembles munis de plusieurs structures de groupement. Dans ce chapitre nous 
allons nous contenter de déduire de quelques exemples la définition d'une notion intéressante
construite à partir de deux structures de groupements.
\section{Les surfaces de Moore}\label{sec:2-gr:1}
Nous avons vue dans le chapitre \ref{chap:moore} que l'ensemble $\BBS(X)$ des surfaces de 
Moore d'un espace topologique $X$ peut être muni de deux structures de groupement (théorème 
\ref{theo:moore:2.1})
$$(\BBS(X);\MFs_{1};\MFt_{1};\OP_{1})\quad\text{et}\quad(\BBS(X);\MFs_{2};\MFt_{2};\OP_{2})$$
qui sont reliées entre elles par les deux propriétés suivantes :
\begin{enumerate}
\item (proposition \ref{prop:moore:2.1}) On a les égalités
         \begin{align*}
         \MFs_{1}\MFs_{2}&=\MFs_{2}\MFs_{1}& \MFt_{1}\MFt_{2}=\MFt_{2}\MFt_{1}\\
         \MFs_{1}\MFt_{2}&=\MFt_{2}\MFs_{1}& \MFs_{2}\MFt_{1}=\MFt_{1}\MFs_{2}
         \end{align*}
\item (proposition \ref{prop:moore:2.2}) Pour tous chemins de Moore $\gamma_{1}$, $\gamma_{2}$, 
         $\gamma_{3}$ et $\gamma_{4}$ tels que
         $$\MFs_{2}(\gamma_{2})=\MFt_{2}(\gamma_{1}),\quad
         \MFs_{2}(\gamma_{4})=\MFt_{2}(\gamma_{3}),\quad
         \MFs_{1}(\gamma_{3})=\MFt_{1}(\gamma_{1}),\quad
         \MFs_{1}(\gamma_{4})=\MFt_{1}(\gamma_{2})$$on a
         $$(\gamma_{4}\OP_{2}\gamma_{3})\OP_{1}(\gamma_{2}\OP_{2}\gamma_{1})
         =(\gamma_{4}\OP_{1}\gamma_{2})\OP_{2}(\gamma_{3}\OP_{1}\gamma_{1})$$
\end{enumerate}
\'Etudions maintenant le cas des espaces topologiques.
\section{Espaces topologiques}\label{sec:2-gr:2}
Si $X$ est un espace topologique, alors, d'après le chapitre \ref{chap:alex}, l'ensemble
des ouverts $\Ouv(X)$ peut être muni de deux structures de groupements
$$(\Ouv(X);\stackrel{\cup}{s};\stackrel{\cup}{t};\cup)\quad
\text{et}\quad(\Ouv(X);\stackrel{\cap}{s};\stackrel{\cap}{t};\cap)$$
Puisque pour tout ouvert $x$ de $X$, on a, par définition,
$$\stackrel{\cup}{s}(x)=\stackrel{\cup}{t}(x)=X
\quad\text{et}\quad
\stackrel{\cap}{s}(x)=\stackrel{\cap}{t}(x)=x$$
il est très facile de vérifier les égalités
\begin{align*}
\stackrel{\cup}{s}\stackrel{\cap}{s}
&=\stackrel{\cap}{s}\stackrel{\cup}{s}
&\stackrel{\cup}{t}\stackrel{\cap}{t}
=\stackrel{\cap}{t}\stackrel{\cup}{t}\\
\stackrel{\cup}{s}\stackrel{\cap}{t}
&=\stackrel{\cap}{t}\stackrel{\cup}{s}
&\stackrel{\cap}{s}\stackrel{\cup}{t}
=\stackrel{\cup}{t}\stackrel{\cap}{s}
\end{align*}
Supposons que $x_{1}$, $x_{2}$, $x_{3}$ et $x_{4}$ soient quatre ouverts de $X$ qui
vérifient les égalités
$$\stackrel{\cup}{s}(x_{2})=\stackrel{\cup}{t}(x_{1}),\quad
\stackrel{\cup}{t}(x_{4})=\stackrel{\cup}{t}(x_{3}),\quad
\stackrel{\cap}{s}(x_{3})=\stackrel{\cap}{t}(x_{1}),\quad
\stackrel{\cap}{s}(x_{4})=\stackrel{\cap}{t}(x_{2})$$
Ces conditions impliquent
$$x_{1}=x_{3}\quad\text{et}\quad x_{2}=x_{4}$$
Cela implique
$$(x_{4}\cup x_{3})\cap(x_{2}\cup x_{1})=
(x_{2}\cup x_{1})\cap(x_{2}\cup x_{1})=x_{2}\cup x_{1}$$
et
$$(x_{4}\cap x_{2})\cup(x_{3}\cap x_{1})
=x_{2}\cup x_{1}$$
D'où
$$(x_{4}\cup x_{3})\cap(x_{2}\cup x_{1})=
(x_{4}\cap x_{2})\cup(x_{3}\cap x_{1})$$
\section{Les g-carrés}\label{sec:2-gr:3}
En utilisant le fait que $\GMOR$ est une catégorie, nous allons montrer 
que sa catégorie des flèches est naturellement muni de deux structures
de catégorie qui satisfont les deux conditions vues dans les deux sections
précédentes. Pour rester cohérent avec nos dénomination, nous parlerons de
g-carrés. L'ensemble $\GSQUARE$ des g-carrés est l'ensemble des quadruplets
$(x_{1};x_{2};y_{1};y_{2})$ où $x_{1}$, $x_{2}$, $y_{1}$ et $y_{2}$
sont des petits g-morphismes tels que
$$x_{2}\OP y_{1}=y_{2}\OP x_{1}$$
Cela signifie que l'on a aussi
$$s(y_{1})=s(x_{1})\quad t(y_{1})=s(x_{2})\quad 
s(y_{2})=t(x_{1})\quad t(y_{2})=t(x_{2})$$
Notons $\MFs_{1},\MFs_{2},\MFt_{1},\MFt_{2}:\GSQUARE\to\GSQUARE$, les
applications définies par
\begin{align*}
\MFs_{1}(x_{1};x_{2};y_{1};y_{2})&=(x_{1};x_{1};\Id;\Id)\\
\MFs_{2}(x_{1};x_{2};y_{1};y_{2})&=(\Id;\Id;y_{1};y_{1})\\
\MFt_{1}(x_{1};x_{2};y_{1};y_{2})&=(x_{2};x_{2};\Id;\Id)\\
\MFt_{2}(x_{1};x_{2};y_{1};y_{2})&=(\Id;\Id;y_{2};y_{2})
\end{align*}
et $\OP_{1};\OP_{2}:\GSQUARE\times\GSQUARE\to\GSQUARE$ les applications
définies par
\begin{align*}
(x_{1};x_{2};y_{1};y_{2})&\OP_{1}(x'_{1};x'_{2};y'_{1};y'_{2})\\
&=
\begin{cases}
(x'_{1};x_{2};y_{1}\OP y'_{1};y_{2}\OP y'_{2})  &\text{si }\MFs_{1}(x_{1};x_{2};y_{1};y_{2})=
\MFt_{1}(x'_{1};x'_{2};y'_{1};y'_{2})\\
(x'_{1};x'_{2};y'_{1};y'_{2})&\text{sinon}
\end{cases}
\end{align*}
et
\begin{align*}
(x_{1};x_{2};y_{1};y_{2})&\OP_{2}(x'_{1};x'_{2};y'_{1};y'_{2})\\
&=
\begin{cases}
(x_{1}\OP x'_{1};x_{2}\OP x'_{2};y'_{1};y_{2})  &\text{si }\MFs_{2}(x_{1};x_{2};y_{1};y_{2})=
\MFt_{2}(x'_{1};x'_{2};y'_{1};y'_{2})\\
(x'_{1};x'_{2};y'_{1};y'_{2})&\text{sinon}
\end{cases}
\end{align*}
Il est très facile de montrer que $(\GSQUARE;\MFs_{1};\MFt_{1};\OP_{1})$ et 
$(\GSQUARE;\MFs_{2};\MFt_{2};\OP_{2})$ sont des groupements et même des 
catégories.

De plus, pour tout $(x_{1};x_{2};y_{1};y_{2})\in\GSQUARE$, on a
$$\MFs_{1}\MFs_{2}(x_{1};x_{2};y_{1};y_{2})=\MFs_{1}(\Id;\Id;y_{1};y_{1})
=(\Id;\Id;\Id;\Id)$$
et
$$\MFs_{2}\MFs_{1}(x_{1};x_{2};y_{1};y_{2})=\MFs_{2}(x_{1};x_{1};\Id;\Id)
=(\Id;\Id;\Id;\Id)$$
Pour être plus précis, on peut facilement voir que dans les deux cas les identités finales 
sont $\Id_{\BBB_{1}}$ où $\BBB_{1}$ est le groupement source de $x_{1}$.
Par conséquent
$$\MFs_{1}\MFs_{2}=\MFs_{2}\MFs_{1}$$
et de même
$$\MFt_{1}\MFt_{2}=\MFt_{2}\MFt_{1}\quad
\MFs_{1}\MFt_{2}=\MFt_{2}\MFs_{1}\quad
\MFs_{2}\MFt_{1}=\MFt_{1}\MFt_{2}$$
Prenons $\gamma^{i}=(x_{1}^{i};x_{2}^{i};y_{1}^{i};y_{2}^{i})$, 
$i=1,\ldots,4$, quatre éléments de $\GSQUARE$ vérifiant
$$ \MFs_2(\gamma^{2})=\MFt_2(\gamma^{1}),
\ \MFs_2(\gamma^{4})=\MFt_2(\gamma^{3}),
\ \MFs_1(\gamma^{3})=\MFt_1(\gamma^{1})
\text{ et }\MFs_1(\gamma^{4})=\MFt_1(\gamma^{2})$$
On a alors
\begin{align*}
(\gamma^{4}\OP_{2}\gamma^{3})\OP_{1}
(\gamma^{2}\OP_{2}\gamma^{1})
&=(x_{1}^{4}\OP x_{1}^{3};x_{2}^{4}\OP x_{2}^{3};y_{1}^{3};y_{2}^{4})
\OP_{1}(x_{1}^{2}\OP x_{1}^{1};x_{2}^{2}\OP x_{2}^{1};y_{1}^{1};y_{2}^{2})\\
&=(x_{1}^{2}\OP x_{1}^{1};x_{2}^{4}\OP x_{2}^{3};y_{1}^{3}\OP y_{1}^{1};
     y_{2}^{4}\OP y_{2}^{2})
\end{align*}
et
\begin{align*}
(\gamma^{4}\OP_{1}\gamma^{2})\OP_{2}
(\gamma^{3}\OP_{1}\gamma^{1})
&=(x_{1}^{2};x_{2}^{4};y_{1}^{4}\OP y_{1}^{2};y_{2}^{4}\OP y_{2}^{2})
\OP_{2}(x_{1}^{1};x_{2}^{3};y_{1}^{3}\OP y_{1}^{1};y_{2}^{3}\OP y_{2}^{1})\\
&=(x_{1}^{2}\OP x_{1}^{1};x_{2}^{4}\OP x_{2}^{3};y_{1}^{3}\OP y_{1}^{1};
     y_{2}^{4}\OP y_{2}^{2})
\end{align*}
On en déduit
$$(\gamma^{4}\OP_{2}\gamma^{3})\OP_{1}
(\gamma^{2}\OP_{2}\gamma^{1})
=(\gamma^{4}\OP_{1}\gamma^{2})\OP_{2}
(\gamma^{3}\OP_{1}\gamma^{1})$$
%\section{Les anneaux}\label{sec:2-gr:3}
%Nous allons ici montrer que tout anneau quelconque $(\BBA;+;\cdot)$, unitaire ou pas, 
%commutatif ou pas, peut être muni de deux structures de groupement.
%
%Posons $$\stackrel{+}{s}(x)=x\quad\stackrel{+}{t}(x)=-x$$et
%$$\dot{s}(x)=\dot{t}(x)=x$$pour tout $x\in\BBA$. Puisque les opérations
%$+$ et $\cdot$ sont associatives, il est immédiat que 
%$$(\BBA;\stackrel{+}{s};\stackrel{+}{t};+)
%\quad\text{et}\quad(\BBA;\dot{s};\dot{t};\cdot)$$
%sont des groupements.
%
%De plus, pour tous éléments $x_{1}$, $x_{2}$, $x_{3}$ et $x_{4}$ de $\BBA$, 
%les égalités
%$$\stackrel{+}{s}(x_{2})=\stackrel{+}{t}(x_{1}),\quad
%\stackrel{+}{t}(x_{4})=\stackrel{+}{t}(x_{3}),\quad
%\dot{s}(x_{3})=\dot{t}(x_{1}),\quad
%\dot{s}(x_{4})=\dot{t}(x_{2})$$
%implique
%$$x_{1}=x_{3}\quad\text{et}\quad x_{2}=x_{4}$$.
%Par conséquent
%$$(x_{4}+ x_{3})\cdot(x_{2}+ x_{1})=
%(x_{2}\cup x_{1})\cap(x_{2}\cup x_{1})=x_{2}\cup x_{1}$$
%et
%$$(x_{4}\cdot x_{2})+(x_{3}\cdot x_{1})
%=x_{2}^{2}+ x_{1}^{2}$$
%D'où
%$$(x_{4}\cup x_{3})\cap(x_{2}\cup x_{1})=
%(x_{4}\cap x_{2})\cup(x_{3}\cap x_{1})$$
\section{Définition d'un $2$-groupement strict}\label{sec:2-gr:4}
Tous les exemples précédents, nous conduisent à la définition ci-dessous.

Nous dirons qu'un septuplet $\left(\BBB;s_1;t_1;\OP_1;s_2;t_2;\OP_2\right)$ est un
\emph{$2$-groupement strict}\index{$2$-groupement strict} quand
\begin{description}
\item[(2-GR 1)] les quadruplets $\left(\BBB;s_1;t_1;\OP_1\right)$ et
$\left(\BBB;s_2;t_2;\OP_2\right)$ sont des groupements ;
\item[(2-GR 2)] les égalités ci-dessous sont satisfaites
$$\MFs_1\MFs_2=\MFs_2\MFs_1\qquad\MFt_1\MFt_2=\MFt_2\MFt1$$
$$\MFs_1\MFt_2=\MFt_2\MFs_1\qquad\MFs_2\MFt_1=\MFt_1\MFs_2$$
\item[(2-GR 3)] si $x_i$,
$i=1,\ldots,4$, sont quatre éléments de $\BBB$ et
si $$ s_2(x_2)=t_2(x_1),\ s_2(x_4)=t_2(x_3),\
s_1(x_3)=t_1(x_1)\text{ et }s_1(x_4)=t_1(x_2)$$
alors les compositions
$$(x_4\OP_2x_3)\OP_1(x_2\OP_2x_1)\quad\text{et}\quad
(x_4\OP_1x_2)\OP_2(x_3\OP_1x_1)$$
existent et sont égales
$$(x_4\OP_2x_3)\OP_1(x_2\OP_2x_1)=
(x_4\OP_1x_2)\OP_2(x_3\OP_1x_1)$$
\end{description}
Puisque nous l'avons voulu ainsi, les sections \ref{sec:2-gr:1}, \ref{sec:2-gr:2} et 
\ref{sec:2-gr:3} sont des exemples de 2-grou\-pe\-ments.

Le premier d'entre eux
est intéressant car il semble que ce soit l'une des premières fois que l'on arrive
à munir l'ensemble des surfaces de Moore de compositions. Ce fait était connu depuis 
très longtemps pour les chemins, mais restait méconnu pour les surfaces.

Il est bien connu que les espaces topologiques définissent une catégorie (ensemble des ouverts 
muni de la structure de préordre induite par l'inclusion). Le défaut de se point de vue est 
que la réunion et l'intersection y sont complètement oubliées. De plus ils est facile de voir
que les applications continues définissent par images réciproques des g-morphismes pour 
les deux structures de groupements envisagées précédemment.
%%%%%%%%%%%%%%%%%%%%%%%%%%%%%%%%%%%%%%%%%%%%%%%%%%

%%%%%%%%%%%%%%%%%%%%%%%%%%%%%%%%%%%%%%%%%%%%%%%%%%
\printindex

\begin{thebibliography}{6}
\bibitem{Ehr65}
            C. Ehresmann
            \textit{Catégories et structures}
            Dunod, 1965
\bibitem{KaVo91}
             M. Kapranov et V. Voevodsky
             \textit{Infinity-groupoids and homotopy types}
             Cahiers de Topologie et Géométrie Différentielle Catégoriques \textbf{32}, 1991
\bibitem{ML71}
             S. MacLane
             \textit{Categories for the Working Mathematician}
             Graduate Texts in Mathematics \textbf{5}, Springer-Verlag, 1971
\bibitem{MSS99}
             P. Mateus, A. Sernadas et C. Sernadas
             \textit{Precategories for Combining Probabilistic Automata}
             Electronic Notes in Theoretical Computer Science \textbf{29}, 1999
\bibitem{ScHe00}
            L. Schröder et H. Herrlich
            \textit{Free Adjunction of Morphisms}
            Applied Categorical Strutures \textbf{8}(4), 2000
\end{thebibliography}
\end{document}